\documentclass[pdflatex,sn-mathphys-num]{sn-jnl}

\usepackage{graphicx}%
\usepackage{multirow}%
\usepackage{amsmath,amssymb,amsfonts}%
\usepackage{amsthm}%
\usepackage{rsfso}%
\usepackage[title]{appendix}%
\usepackage{xcolor}%
\usepackage{textcomp}%
\usepackage{manyfoot}%
\usepackage{booktabs}%
\usepackage{algorithm}%
\usepackage{algorithmic}
\usepackage{caption}%
\usepackage{subcaption}
\usepackage{enumerate}
\usepackage{listings}%

\allowdisplaybreaks
\hypersetup{hypertexnames=false}

\def\bbR{\mathbb{R}}

\DeclareMathOperator{\F}{F}
\DeclareMathOperator{\T}{T}

\DeclareMathOperator{\GRAD}{grad}

\DeclareMathOperator{\sym}{sym}

\DeclareMathOperator{\tr}{tr}

\newcommand\STM[2]{{\rm St}(#1,#2)}

\def\caR{\mathcal{R}}

\def\caF{\mathcal{F}}

\def\caL{\mathcal{L}}

\def\caW{\mathcal{W}}

\def\caZ{\mathcal{Z}}
\def\caB{\mathcal{B}}
\def\caU{\mathcal{U}}
\def\caV{\mathcal{V}}
\def\caK{\mathcal{K}}

\def\lam{\Lambda}

\def\l{\left}
\def\r{\right}

\def\cT{\mathcal{T}}

\theoremstyle{thmstyleone}%
\newtheorem{theorem}{Theorem}
\newtheorem{lemma}[theorem]{Lemma}

\theoremstyle{thmstyletwo}%

\theoremstyle{thmstylethree}%
\newtheorem{definition}{Definition}%

\begin{document}

\title[Splitting and Alternating Minimization Methods for OrthDFNs]{Splitting and Alternating Minimization Methods for Training Orthogonal Deep Feedforward Networks}


\author*[1]{\fnm{Jiawei} \sur{Wang}}
\email{22483004@life.hkbu.edu.hk}

\author[1]{\fnm{Richard Yi Da} \sur{Xu}}
\email{xuyida@hkbu.edu.hk}

\affil*[1]{%
    \orgdiv{Department of Mathematics},
    \orgname{Hong Kong Baptist University},
    \orgaddress{%
        \street{Kowloon Tong},
        \city{Hong Kong},
        \country{China}%
    }%
}




\abstract{Training deep feedforward networks with weight matrices constrained to layer-wise Stiefel manifolds serves as an effective regularization technique to improve generalization. However, the simultaneous presence of orthogonality constraints and nonsmooth activation functions poses severe theoretical and computational challenges for standard stochastic gradient and Riemannian optimization methods. 
In this paper, we propose two splitting and alternating-based algorithms to train orthogonal deep feedforward networks (OrthDFNs). By introducing auxiliary variables, we decouple the nested nonsmooth compositions and manifold constraints into a structured equality-constrained optimization problem. 
To solve this model, we develop OrthALM, an inexact augmented Lagrangian method (ALM) that employs a proximal alternating minimization scheme to solve its subproblems. Crucially, every block coordinate subproblem admits an exact, closed-form solution.  
We establish a rigorous convergence theory for OrthALM, proving that every accumulation point of the generated sequence is a Karush--Kuhn--Tucker (KKT) point of the regularized, reformulated problem.
Furthermore, we introduce OrthADMM as a practical, proximal-free variant that executes a single inner alternating sweep per multiplier update to accelerate computation. 
Extensive numerical experiments show that OrthALM attains lower training error and constraint residuals than OrthADMM in most configurations. OrthADMM achieves competitive test accuracy and the lowest reported runtime in most comparisons with the selected baselines.}

\keywords{Stiefel manifold, Orthogonal deep feedforward networks, Augmented Lagrangian method, ADMM, Nonsmooth nonconvex optimization}

\pacs[MSC Classification]{90C26, 90C30, 65K05, 68T07, 15B10}

\maketitle

\makeatletter
\let\@oddhead\@empty
\let\@evenhead\@empty
\makeatother

\section{Introduction}\label{sec.intro}

Deep feedforward networks (DFNs) have demonstrated significant advantages across a wide range of machine learning tasks. The deep architectures of DFNs endow them with high model capacities, allowing them to represent highly complex relationships between inputs and outputs. However, DFNs are usually over-parameterized; this massive number of parameters, coupled with their complex architectures, introduces a potential risk of overfitting. 

To alleviate this issue, \cite{li2019orthogonal, huang2018orthogonal} propose training DFNs with weight matrices constrained to their respective Stiefel manifolds. Enforcing this orthogonality acts as a regularization technique, which can improve the generalization ability of networks.
In this paper, we focus on the parameter estimation problem of DFNs whose weight matrices are constrained to the Stiefel manifold:
\begin{equation}\label{eq:STM}
    \STM{p}{n} := \left\{ X \in \bbR^{n\times p} \;\middle|\; X^{\T}X = I_p \right\} \quad (p \le n).
\end{equation}
Given a set of input data $x_i \in \mathbb{R}^{n_0}$ and corresponding output data $y_i \in \mathbb{R}^{n_L}$, for $i = 1, \dots, m$, the problem of training an $L$-layer OrthDFN can be formulated as the following optimization problem: 
\begin{equation}\label{orthdnn.ori.prob}
    \begin{aligned}
        \min_{W_l \in \mathcal{W}, b_l \in \mathcal{B}} \quad & \frac{1}{m} \sum_{i = 1}^{m} \left\| \sigma \Big( W_L \sigma \big( \dots \sigma ( W_1 x_i + b_1 ) \dots \big) + b_L \Big) - y_i \right\|^2 \\
        \text{s.t.} \quad & 
        \begin{cases} 
            W_l^{\T}W_l = I_{n_{l-1}}, & \text{if } n_{l} \geq n_{l-1}, \\ 
            W_l W_l^{\T} = I_{n_l}, & \text{if } n_{l-1} > n_{l}, 
        \end{cases} \\
        & \text{for } l = 1, \dots, L,
    \end{aligned}
\end{equation}
where $\mathcal{W} := \{W_l \in \mathbb{R}^{n_l \times n_{l-1}}\}_{l = 1}^L$ and $\mathcal{B} := \{b_l \in \mathbb{R}^{n_l}\}_{l = 1}^L$ are the weights and biases, respectively. For the nonsmooth activation function $\sigma(\cdot)$, we employ either the Rectified Linear Unit (ReLU), defined element-wise as
\[
\sigma_{\mathrm{relu}}(u) := \max \{ u, 0 \}, \quad u \in \mathbb{R}, 
\]
or the leaky ReLU, defined as
\[
\sigma_{\mathrm{lrelu}}(u) := \max \{ u, \alpha u \}, \quad \alpha \in (0,1), \; u \in \mathbb{R}.
\]

\subsection{Motivation and Related Work}\label{intro.moti}

While stochastic gradient descent (SGD)-based methods are most commonly adopted for training modern deep neural networks, directly applying them to solve problem \eqref{orthdnn.ori.prob} presents two major difficulties. The first one is due to the use of nonsmooth activation functions (such as ReLU and leaky ReLU). Second, the Stiefel manifold constraints pose a significant optimization challenge. Although SGD methods combined with Riemannian optimization techniques can be adopted to address \eqref{orthdnn.ori.prob}, they incur high computational costs \cite{li2019orthogonal} and frequently suffer instability in convergence, leading to inferior empirical performance \cite{huang2018orthogonal}. 

To circumvent these manifold-related issues, alternative strategies have been proposed. Assuming $n_l \leq n_{l-1}$ to ensure that the rows of $W_l$ are strictly orthogonal, Huang et al. \cite{huang2018orthogonal} propose an orthogonal weight normalization approach. To handle scenarios where $n_l > n_{l-1}$, they introduce a group-based orthogonalization method; we refer interested readers to \cite{huang2018orthogonal} for further details. Fundamentally, this orthogonal weight normalization ensures the orthogonality of a weight matrix by constructing an orthogonal transformation over unconstrained proxy parameters. Leveraging this framework, the authors construct the Orthogonal Linear Module (OLM) for practical implementation, explicitly deriving both its forward and backward passes. Instead of constraining the weight matrices strictly on the Stiefel manifold, the authors of \cite{li2019orthogonal} develop a Singular Value Bounding (SVB) algorithm, which permits slight deviations from the Stiefel manifold to lower computational expense. Despite these practical innovations, the aforementioned methods for solving \eqref{orthdnn.ori.prob} are still SGD-based, lack rigorous convergence analysis, and are usually sensitive to their algorithmic hyperparameters. 

In the absence of Stiefel manifold constraints, several splitting and alternating-based methods have been proposed as alternatives to SGD for deep learning optimization. For instance, the method of auxiliary coordinates (MAC), introduced by Carreira-Perpi\~{n}\'{a}n and Wang \cite{carreira2014distributed} in 2014, transforms standard neural network training into an equality-constrained optimization problem via auxiliary variables. The constraints are then enforced using quadratic penalties.
The ADMM framework has also been adopted for deep learning applications. To avoid the saturation problem inherent to DFNs with sigmoid-type activation functions, the authors of \cite{zeng2021admm} developed an ADMM-based algorithm. While they provide convergence analysis, the theoretical guarantees rely on highly restrictive assumptions: the activation function must be bounded, twice-differentiable with bounded first- and second-order derivatives, and either real analytic or semialgebraic.
Furthermore, augmented Lagrangian-based methods have been explored for network training. By introducing auxiliary variables and constraints, \cite{liu2023inexact} reformulates the problem of leaky ReLU network with the $\ell_{2,1}$ regularizer. This is solved via an inexact ALM, utilizing an alternating minimization method for the subproblems. Similarly, for Elman recurrent neural networks \cite{elman1990finding} with ReLU activations, Wang, Zhang, and Chen \cite{wang2025augmented} reformulate the network training as a constrained optimization problem, solving it with an ALM whose subproblems are tackled by the block coordinate descent (BCD) method.

Crucially, none of the aforementioned splitting and alternating-based methods consider the scenario where the network's weight matrices are constrained to the Stiefel manifold, which introduces fundamental modeling and theoretical differences. Motivated by these limitations, we aim to design a numerically stable approach to efficiently solve OrthDFNs with rigorous convergence guarantees.

\subsection{Problem Reformulation}\label{intro.prob.refo}

A natural strategy for handling this multi-layer nonsmooth optimization problem is to split the difficulties caused by the orthogonality constraints and the composite objective of \eqref{orthdnn.ori.prob}. Denote $X := [x_1, x_2,\cdots,x_m] \in \mathbb{R}^{n_0 \times m}$ and $Y := [y_1,y_2,\cdots,y_m] \in \mathbb{R}^{n_L \times m}$ as the aggregated input and output data matrices, respectively. By introducing auxiliary variables, problem \eqref{orthdnn.ori.prob} can be equivalently reformulated as follows: 
\begin{equation}\label{ori.prob.orth.three.sets.vari}
    \begin{aligned}
        \min_{\substack{
    W_l \in \mathcal{W},\; Z_l \in \mathcal{Z},\; b_l \in \mathcal{B},\\
    V_l \in \mathcal{V},\; U_l \in \mathcal{U}
}} \quad & \frac{1}{m} \left\Vert V_{L} - Y \right\Vert_{\F}^2 \\
        \text{s.t.} \quad & W_l V_{l-1} + b_l 1_m^{\T} = U_l, \quad V_l = \sigma ( U_{l} ), \\
        & \begin{cases} 
            W_l = Z_l, \quad Z_l \in \STM{n_{l-1}}{n_l}, & \text{if } n_{l} \geq n_{l-1}, \\ 
            W_l^{\T} = Z_l, \quad Z_l \in \STM{n_l}{n_{l-1}}, & \text{if } n_{l-1} > n_{l}, 
        \end{cases} \\       
        & \text{for } l = 1, \dots, L.
    \end{aligned}
\end{equation}
Here, $\mathcal{U} := \{U_l \in \mathbb{R}^{n_l \times m}\}_{l = 1}^L$ and $\mathcal{V} := \{V_l \in \mathbb{R}^{n_l \times m}\}_{l = 1}^L$ denote the collections of pre- and post-activation auxiliary matrices, while $\mathcal{Z} := \{Z_l\}_{l=1}^L$ contains the auxiliary weight variables. The initial state is given by $V_0 = X \in \mathbb{R}^{n_0 \times m}$. The use of auxiliary variables allows us to decouple the composite structure of the original objective function and reduce the objective in problem \eqref{ori.prob.orth.three.sets.vari} to a simple smooth function. This separation allows us to effectively employ splitting and alternating-based algorithms to deal with the reformulated model.  

To mitigate overfitting and improve the robustness of the model, we further introduce regularization terms for the biases $\mathcal{B}$ and the pre-activation variables $\mathcal{U}$. This yields the following regularized minimization problem:
\begin{align}\label{prob.orth.three.sets.vari}
    \min_{\substack{
    W_l \in \mathcal{W},\; Z_l \in \mathcal{Z},\; b_l \in \mathcal{B},\\
    V_l \in \mathcal{V},\; U_l \in \mathcal{U}
}} \quad & \frac{1}{m} \left\Vert V_{L} - Y \right\Vert_{\F}^{2} + \sum_{l = 1}^L \left( \frac{\lambda_1}{2}\|b_l\|^2 + \frac{\lambda_2}{2} \|U_l\|_{\F}^2 \right) \notag\\
    \text{s.t.} \quad & W_l V_{l-1} + b_l1_m^{\T} = U_l, \quad V_l = \sigma\left( U_{l} \right), \notag\\
     & \begin{cases} 
        W_l = Z_l, \quad Z_l \in \STM{n_{l-1}}{n_l}, & \text{if } n_{l} \geq n_{l-1},  \\ 
        W_l^{\T} = Z_l, \quad Z_l \in \STM{n_l}{n_{l-1}}, & \text{if } n_{l-1} > n_l, 
    \end{cases} \notag\\
    & \text{for } l = 1, \dots, L,
\end{align}
where $\lambda_1, \lambda_2 > 0$ are the regularization parameters.

\subsection{Contributions}\label{intro.cont}

To the best of our knowledge, this is the first study to train OrthDFNs using splitting and alternating-based methods---OrthALM and OrthADMM---supported by the convergence analysis of OrthALM for the regularized problem~\eqref{prob.orth.three.sets.vari}. 
Our main contributions are summarized below:

\begin{itemize}

    \item \textbf{Algorithm Design:} To solve the regularized problem~\eqref{prob.orth.three.sets.vari}, we propose OrthALM, an augmented Lagrangian method that employs an alternating minimization scheme with exact, closed-form updates for every variable block. We additionally develop OrthADMM as an efficient, tailored variant to solve the unregularized problem~\eqref{ori.prob.orth.three.sets.vari}.

    \item \textbf{Rigorous Convergence Analysis:} We establish theoretical convergence guarantees for our proposed OrthALM approach. Specifically, we prove that in the $k$-th iteration of the augmented Lagrangian method, the inner alternating minimization algorithm is guaranteed to meet its stopping criterion within $\mathcal{O}(1/\epsilon^2)$ steps. Then we demonstrate that the sequence generated by our augmented Lagrangian method has at least one accumulation point, and that any accumulation point is a KKT point of the regularized problem \eqref{prob.orth.three.sets.vari}.

    \item \textbf{Numerical Evaluation:} We evaluate our proposed algorithms through extensive numerical experiments on several benchmark datasets. 
    The numerical results show that OrthALM achieves lower training error and constraint violations in most of the tested configurations. 
    Meanwhile, OrthADMM emerges as a highly computationally efficient framework, achieving competitive test accuracy while delivering substantial CPU time speedups over state-of-the-art SGD-based and splitting-based solvers.
    
\end{itemize}

\subsection{Organization}\label{intro.organ}

The rest of this paper is organized as follows. 
In Section~\ref{sec:framework}, we establish the KKT conditions for problem~\eqref{prob.orth.three.sets.vari} and introduce some results. Then, we present the detailed algorithmic frameworks and closed-form subproblem updates for OrthALM and OrthADMM. 
In Section~\ref{sec:convergence}, we provide the convergence analysis for our proposed OrthALM framework. Section~\ref{sec:experiments} presents comprehensive numerical results.
Finally, concluding remarks and are drawn in Section~\ref{sec:conclusion}.

\subsection{Notation and Terminology}\label{intro.nota}

Throughout this paper, $\bbR^{n \times m}$ denotes the set of $n \times m$ real matrices, $\bbR^n = \bbR^{n \times 1}$, and $\bbR = \bbR^1$. We denote by $1_m \in \bbR^m$ the all-ones vector and by $I_n \in \bbR^{n \times n}$ the $n \times n$ identity matrix. 
For a real number $x \in \bbR$, $\lceil x \rceil$ denotes the smallest integer greater than or equal to $x$. For a vector $v \in \bbR^n$, $\|v\|$ denotes its Euclidean $\ell_2$-norm. For a matrix $A = [A_{(i,j)}] \in \bbR^{n \times m}$, $A^{\T}$ denotes its transpose, $\|A\|_{\F} = \sqrt{\tr(A^{\T}A)}$ denotes its Frobenius norm, and $\|A\|_{\max} = \max_{i,j} |A_{(i,j)}|$ denotes its entrywise maximum norm. 
The Euclidean inner product of two matrices $A, B \in \bbR^{n \times m}$ is defined as $\langle A, B \rangle = \tr(A^{\T}B)$, and $A \odot B$ denotes their Hadamard product. For a square matrix $X$, $\sym(X) = (X + X^{\T})/2$.

For a point $Z \in \STM{p}{n}$ on the Stiefel manifold defined in \eqref{eq:STM}, the tangent space to $\STM{p}{n}$ at $Z$ is given by $\cT_Z\STM{p}{n} = \{ V \in \bbR^{n \times p} \mid Z^{\T}V + V^{\T}Z = 0 \}$. The orthogonal projection of a matrix $P \in \bbR^{n \times p}$ onto the tangent space $\cT_Z\STM{p}{n}$ is given by (see, e.g., \cite{absil2008optimization,li2024theory})
\[
    {\mit{\Pi}}_Z(P) := P - Z\sym(Z^{\T}P).
\]
For a smooth function $f$, $\nabla f$ denotes its Euclidean gradient, and its Riemannian gradient with respect to $\STM{p}{n}$ at $Z$ is $\GRAD f_{|\STM{p}{n}}(Z) = {\mit{\Pi}}_Z(\nabla f(Z))$. 
For a locally Lipschitz nonsmooth function $r$, $\partial r$ denotes its limiting subdifferential \cite{rockafellar2009variational}.

\section{KKT Conditions and Algorithm Framework}\label{sec:framework}

In this section, we establish the KKT conditions for problem~\eqref{prob.orth.three.sets.vari} and exploit the structure of the reformulated problems to design splitting and alternating-based methods for OrthDFNs training.

\subsection{KKT Conditions of Problem \texorpdfstring{\eqref{prob.orth.three.sets.vari}}{(3)}}\label{subsec:kkt}

In this subsection, we first establish the KKT conditions for the regularized problem \eqref{prob.orth.three.sets.vari}. Subsequently, we derive several foundational properties that will help the convergence analysis presented in subsequent sections. 
To streamline the notation in the subsequent analysis, we first group the primal variables and the upcoming Lagrangian multipliers into the tuples $\Theta$ and $\Lambda$, respectively, defined as:
\[
    \Theta := \left(\caW, \caZ, \caB, \caV, \caU \right), \quad \text{and} \quad \Lambda := \left(\{\Lambda_{1,l}\}_{l=1}^L, \{\Lambda_{2,l}\}_{l=1}^L, \{\Lambda_{3,l}\}_{l=1}^L\right).
\]
Additionally, we define the index set $\mathcal{I}$ to denote the layers where the number of neurons is non-decreasing:
\[
    \mathcal{I} := \{ l \in \{1, \dots, L\} \mid n_l \geq n_{l-1} \},
\]
and we let $\mathcal{I}^c := \{1, \dots, L\} \setminus \mathcal{I}$ denote its complement. 
We express the objective function as:
\[
\caF (\Theta) := \frac{1}{m}\left\Vert V_{L} - Y \right\Vert_{\F}^{2} + \sum_{l = 1}^L \left(\frac{\lambda_1}{2}\|b_l\|^2 + \frac{\lambda_2}{2} \|U_l\|_{\F}^2 \right).
\]
Next, to write the two weight orientations uniformly, we define the linear isometry
and its adjoint as
\[
\mathcal T_l(A)=
\begin{cases}A,&l\in\mathcal I,\\ A^{\T},&l\in\mathcal I^c,\end{cases}
\qquad
\mathcal T_l^*(B)=
\begin{cases}B,&l\in\mathcal I,\\ B^{\T},&l\in\mathcal I^c.\end{cases}
\]
Here $A$ has the dimensions of $W_l$ and $B$ those of $Z_l$.
We also abbreviate the three equality-constraint residuals as
\[
E_{1,l}(\Theta) = W_lV_{l-1}+b_l1_m^{\T}-U_l, \quad E_{2,l}(\Theta) = V_l-\sigma(U_l), \quad E_{3,l}(\Theta) = \mathcal T_l(W_l)-Z_l.
\]
We write $E_{a,l}^k=E_{a,l}(\Theta^k)$ and
$\|\Lambda\|^2=\sum_{a=1}^3\sum_{l=1}^L\|\Lambda_{a,l}\|_{\F}^2$.

The Lagrangian function of \eqref{prob.orth.three.sets.vari} is
\[
\caL(\Theta,\Lambda)
=\caF(\Theta)+\sum_{a=1}^3\sum_{l=1}^L
\langle\Lambda_{a,l},E_{a,l}(\Theta)\rangle,
\]
where $\{ \lam_{1,l} \in \mathbb{R}^{n_l \times m} \}_{l = 1}^L $, $ \{ \lam_{2,l} \in \mathbb{R}^{n_l \times m} \}_{l = 1}^L$, $ \{ \lam_{3,l} \in \mathbb{R}^{n_l \times n_{l-1}} \}_{l \in \mathcal{I}}$, and $ \{ \lam_{3,l} \in \mathbb{R}^{n_{l-1} \times n_l} \}_{l \in \mathcal{I}^c}$ are the multiplier matrices. 
The KKT system used in this paper consists of the following stationarity and feasibility conditions for a primal--multiplier pair $(\Theta^\ast,\Lambda^\ast)$. A feasible point admitting such multipliers is called a KKT point, as formalized below:
\begin{subequations}
\begin{align}
    0 & = \nabla_{W_l} \caL \left( \Theta^\ast, \Lambda^\ast \right) = \Lambda_{1,l}^\ast {V_{l-1}^\ast}^{\T} + \Lambda_{3,l}^\ast, \quad l \in \mathcal{I}, \label{orthdnn.kkt.w.i} \\
    0 & = \nabla_{W_l} \caL \left( \Theta^\ast, \Lambda^\ast \right) = \Lambda_{1,l}^\ast {V_{l-1}^\ast}^{\T} + {\Lambda_{3,l}^\ast}^{\T}, \quad l \in \mathcal{I}^c, \label{orthdnn.kkt.w.i.c} \\
    0 & = {\mit{\Pi}}_{Z_l^{\ast}} \l( \nabla_{Z_l} \caL \left( \Theta^\ast, \Lambda^\ast \right)  \r) = {\mit{\Pi}}_{Z_l^{\ast}} \l( - \lam_{3,l}^\ast \r), \quad l = 1,\dots,L, \label{orthdnn.kkt.z} \\
    0 & = \nabla_{b_l} \caL \left( \Theta^\ast, \Lambda^\ast \right) = \Lambda_{1,l}^\ast 1_m + \lambda_1 b_l^\ast, \quad l = 1,\dots,L, \label{orthdnn.kkt.b} \\
    0 & \in \partial_{U_l} \caL \left( \Theta^\ast, \Lambda^\ast \right) = \lambda_2 U_l^\ast -\Lambda_{1,l}^\ast + \partial_{U_l} \langle - \lam_{2,l}^\ast , \sigma(\cdot) \rangle \l( U_l^\ast \r), \quad l = 1,\dots,L, \label{orthdnn.kkt.u} \\
    0 & = \nabla_{V_l} \caL \left( \Theta^\ast, \Lambda^\ast \right) = {W_{l+1}^\ast}^{\T} \Lambda_{1,l+1}^\ast + \Lambda_{2,l}^\ast, \quad l = 1,\dots,L-1, \label{orthdnn.kkt.vl} \\
    0 & = \nabla_{V_L} \caL \left( \Theta^\ast, \Lambda^\ast \right) = \frac{2}{m}(V_L^\ast - Y) + \Lambda_{2,L}^\ast, \label{orthdnn.kkt.vlast} \\
    0 & = U_l^\ast - W_{l}^\ast V_{l-1}^\ast - b_{l}^\ast1_{m}^{\T}, \quad l = 1,\dots,L, \label{orthdnn.kkt.lam1} \\
    0 & = \sigma(U_l^\ast) - V_l^\ast, \quad l = 1,\dots,L, \label{orthdnn.kkt.lam2} \\
    0 & = Z_l^\ast - W_l^\ast, \quad l \in \mathcal{I}, \label{orthdnn.kkt.lam3.i} \\
    0 & = Z_l^\ast - {W_l^\ast}^{\T}, \quad l \in \mathcal{I}^c, \label{orthdnn.kkt.lam3.i.c} \\
    Z_l^\ast & \in \STM{n_{l-1}}{n_l}, \quad l \in \mathcal{I}, \label{orthdnn.kkt.stiefel.i} \\
    Z_l^\ast & \in \STM{n_l}{n_{l-1}}, \quad l \in \mathcal{I}^c, \label{orthdnn.kkt.stiefel.i.c}
\end{align}
\end{subequations}
where $V_0^\ast = X$.
We define the feasible set for problem \eqref{prob.orth.three.sets.vari} as:
\[
    \caR := \left\{ \Theta \;\middle|\; 
    \begin{aligned}
        & W_l V_{l-1} + b_l 1_m^{\T} = U_l, \quad V_l = \sigma(U_l), \\
        & 
        \begin{cases} 
            W_l = Z_l, \quad Z_l \in \mathrm{St}(n_{l-1}, n_l), & \text{if } n_l \geq n_{l-1}, \\ 
            W_l^{\T} = Z_l, \quad Z_l \in \mathrm{St}(n_l, n_{l-1}), & \text{if } n_{l-1} > n_l, 
        \end{cases} \\[1ex]
        & \text{for } l = 1,\dots,L
    \end{aligned}
    \right\}.
\]

Furthermore, $\partial_{\Theta} \caL$ is defined as:
\[
    \begin{aligned}
        \partial_{\Theta} \caL \left( \Theta, \Lambda \right) = \; & \bigg( \left\{ \nabla_{W_l} \caL \left( \Theta, \Lambda \right) \right\}_{l = 1}^L, \; \left\{ {\mit{\Pi}}_{Z_l} \l( \nabla_{Z_l} \caL \left( \Theta, \Lambda \right)  \r) \right\}_{l = 1}^L, \\
        & \;\;\; \left\{ \nabla_{b_l} \caL \left( \Theta, \Lambda \right) \right\}_{l = 1}^L, \; \left\{ \nabla_{V_l} \caL \left( \Theta, \Lambda \right) \right\}_{l = 1}^L, \; \left\{ \partial_{U_l} \caL \left( \Theta, \Lambda \right) \right\}_{l = 1}^L  \bigg).
    \end{aligned}
\]
Here $\partial_\Theta$ denotes the product-manifold stationarity mapping: it collects Euclidean gradients in the smooth unrestricted blocks, limiting subdifferentials in the $U$ blocks, and tangent projections of the gradients in the Stiefel blocks. 
Consequently, the equivalent KKT conditions for problem \eqref{prob.orth.three.sets.vari} are formalized in the following definition.

\begin{definition}

We say that $\Theta^\ast\in\caR$ is a KKT point of problem~\eqref{prob.orth.three.sets.vari} if there exists a multiplier tuple $\Lambda^\ast$, with the dimensions specified above, such that
\[
0 \in \partial_{\Theta} \caL \l( \Theta^\ast, \lam^\ast \r).
\]

\end{definition}

The following lemma excludes nonzero abnormal multipliers in the system used in our convergence analysis.

\begin{lemma}\label{le.partial.phat.0}

For every feasible point $\Theta\in\caR$ of problem~\eqref{prob.orth.three.sets.vari}, the only multiplier tuple $\Lambda$ satisfying the following system is the zero tuple:
    \begin{equation}\label{eq.nnamcq}
        \begin{aligned}
            0 & = \Lambda_{1,l} V_{l-1}^{\T} + \Lambda_{3,l}, \quad && l \in \mathcal{I}, \\
            0 & = \Lambda_{1,l} V_{l-1}^{\T} + \Lambda_{3,l}^{\T}, \quad && l \in \mathcal{I}^c, \\
            0 & = \Lambda_{1,l} 1_m, \quad && l = 1,\dots,L, \\
            0 & \in - \Lambda_{1,l} + \partial_{U_l} \langle - \lam_{2,l}, \sigma(\cdot) \rangle \l( U_l \r), \quad && l = 1,\dots,L, \\
            0 & = W_{l+1}^{\T} \Lambda_{1,l+1} + \Lambda_{2,l}, \quad && l = 1,\dots,L-1, \\
            0 & = \Lambda_{2,L}.
        \end{aligned}
    \end{equation}
\end{lemma}

\begin{proof}

The last equation of \eqref{eq.nnamcq} gives $\Lambda_{2,L}=0$. Suppose that $\Lambda_{2,l}=0$ for some $l\in\{1,\ldots,L\}$. Since the limiting subdifferential of the zero function is $\{0\}$, the $U_l$ inclusion in \eqref{eq.nnamcq} yields $\Lambda_{1,l}=0$. The corresponding $W_l$ equation then gives $\Lambda_{3,l}=0$, whether $l\in\mathcal I$ or $l\in\mathcal I^c$. If $l\ge2$, the $V_{l-1}$ equation further implies
\[
\Lambda_{2,l-1}=-W_l^{\T}\Lambda_{1,l}=0.
\]
Backward induction from $l=L$ to $l=1$ therefore proves
\[
\Lambda_{1,l}=0,\quad\Lambda_{2,l}=0,\quad\Lambda_{3,l}=0,
\qquad l=1,\ldots,L.
\]
For $L=1$, the first induction step already proves the conclusion. This completes the proof.

\end{proof}

\subsection{OrthALM for Problem \texorpdfstring{\eqref{prob.orth.three.sets.vari}}{(3)}}\label{almam.al}

\subsubsection{The Augmented Lagrangian Method}\label{almam.al.alm}
    
Penalizing the equality constraint of the regularized problem \eqref{prob.orth.three.sets.vari} yields the following augmented Lagrangian function:
\begin{align}\label{eq.aug.lag}
\caL_\beta(\Theta,\Lambda)
&=\caF(\Theta)+\sum_{a=1}^3\sum_{l=1}^L
\left[\langle\Lambda_{a,l},E_{a,l}(\Theta)\rangle
+\frac{\beta}{2}\|E_{a,l}(\Theta)\|_{\F}^2\right]\notag\\
&=\caF(\Theta)+\frac{\beta}{2}\sum_{a=1}^3\sum_{l=1}^L
\left\|E_{a,l}(\Theta)+\frac{\Lambda_{a,l}}{\beta}\right\|_{\F}^2
-\frac{\|\Lambda\|^2}{2\beta}.
\end{align}
where $\{ \Lambda_{1,l} \in \mathbb{R}^{n_l \times m} \}_{l = 1}^L$ and $\{ \Lambda_{2,l} \in \mathbb{R}^{n_l \times m} \}_{l = 1}^L$, $\{ \Lambda_{3,l} \in \mathbb{R}^{n_l \times n_{l-1}} \}_{l \in \mathcal{I}}$ and $\{ \Lambda_{3,l} \in \mathbb{R}^{n_{l-1} \times n_l} \}_{l \in \mathcal{I}^c}$ represent the multiplier matrices, and $\beta > 0$ is the penalty parameter for the quadratic penalty terms.
Let
\[
\mathcal D=\left\{\Theta:\ Z_l\in
\operatorname{St}\big(\min\{n_{l-1},n_l\},\max\{n_{l-1},n_l\}\big),
\quad l=1,\ldots,L\right\}.
\]
We derive the following properties of the augmented Lagrangian function $\caL_\beta$.

\begin{lemma}\label{le.bcd.basic}
    For any fixed $\beta$ and $\lam$, the following statements hold:
    \begin{enumerate}
        \item The augmented Lagrangian function $\caL_\beta$ is lower bounded and satisfies
        \[
        \caL_{\beta} \l( \Theta, \Lambda \r) \ge - \sum_{l = 1}^L \frac{\l\| \lam_{1,l} \r\|_{\F}^2 + \l\| \lam_{2,l} \r\|_{\F}^2 + \l\| \lam_{3,l} \r\|_{\F}^2}{2\beta} \;\;\;\; {\rm for \; all } \; \Theta.
        \]
        \item For any $\hat{\Theta} \in\mathcal D $ and $ \tilde{c} \ge \hat r := \caL_\beta \l(\hat{\Theta}, \Lambda \r)$, the level set
        \[
        \Omega_{\caL_{\beta}}(\tilde{c}) := \{ \Theta \in\mathcal D : \caL_{\beta} \l( \Theta, \Lambda \r) \le \tilde{c} \}
        \]
        is nonempty and compact. 
    \end{enumerate}
\end{lemma}

\begin{proof}
\begin{enumerate}
        \item The first statement follows directly from the definition of the augmented Lagrangian function $\caL_\beta$. Since the squared Frobenius norms and the sum of squared vector norms are inherently non-negative, we have:
        \[
        \caL_\beta(\Theta,\Lambda) = \caF(\Theta)+\frac{\beta}{2}\sum_{a=1}^3\sum_{l=1}^L
\left\|E_{a,l}(\Theta)+\frac{\Lambda_{a,l}}{\beta}\right\|_{\F}^2 - \frac{\|\Lambda\|^2}{2\beta} \ge 0-\frac{\|\Lambda\|^2}{2\beta}.
        \]

        \item By definition, the level set $\Omega_{\mathcal{L}_{\beta}}(\tilde{c})$ is nonempty, and it is closed because $\mathcal D$ is closed and the augmented Lagrangian function $\mathcal{L}_{\beta}$ is continuous. To establish the boundedness, we know that $\mathcal{L}_\beta (\Theta, \Lambda ) \le \tilde{c}$ for all $\Theta \in \Omega_{\mathcal{L}_{\beta}}(\tilde{c})$. This upper bound guarantees that the penalty terms $\left\| W_l - Z_l + \Lambda_{3,l}/\beta \right\|_{\F}$ or $\left\| W_l^{\T} - Z_l + \Lambda_{3,l}/\beta \right\|_{\F}$ are strictly bounded for all $l = 1, \dots, L$. 
        Because the variables $\{Z_l\}_{l=1}^L$ are constrained to the Stiefel manifold, they are inherently bounded. Consequently, it follows that the variables $\{W_l\}_{l=1}^L$ must also be bounded. Furthermore, the upper bound on the objective restricts the regularization terms $\|b_l\|^2$ and $\|U_l\|_{\F}^2$, as well as the penalty terms $\left\| V_l - \sigma(U_l) + \Lambda_{2,l}/\beta \right\|_{\F}^2$ for all $l = 1, \dots, L$. This directly implies that the sequences $\{b_l\}_{l=1}^L$, $\{U_l\}_{l=1}^L$, and subsequently $\{V_l\}_{l=1}^L$, are all bounded. Since every component of the tuple $\Theta$ is bounded, the level set $\Omega_{\mathcal{L}_\beta}(\tilde{c})$ is closed and bounded. Therefore, it is both nonempty and compact.

\end{enumerate}

\end{proof}

Within the proposed framework, problem \eqref{prob.orth.three.sets.vari} is solved using an augmented Lagrangian scheme, as outlined in Algorithm~\ref{alg:alm}. At the $k$-th iteration, we first update the primal variables by holding the Lagrange multipliers and the penalty parameter fixed and solving the following subproblem:
\begin{equation} \label{eq:3.5}
    \begin{aligned}
        \min_{\Theta} \quad & \caL_{\beta^{k-1}} ( \Theta, \Lambda^{k-1} ) \\
        \text{s.t.} \quad & Z_l \in 
        \begin{cases} 
            \mathrm{St}(n_{l-1}, n_{l}), & \text{if } l \in \mathcal{I},  \\ 
            \mathrm{St}(n_{l}, n_{l-1}), & \text{if } l \in \mathcal{I}^c,
        \end{cases} \quad \forall l \in \{1,\dots,L\},
    \end{aligned}
\end{equation}
where $\Lambda^{k-1} = \left( \{ \Lambda_{1,l}^{k-1} \}_{l = 1}^L, \{ \Lambda_{2,l}^{k-1} \}_{l = 1}^L, \{ \Lambda_{3,l}^{k-1} \}_{l = 1}^L \right)$. Rather than solving \eqref{eq:3.5} exactly, we compute an approximate primal update $\Theta^k$ that satisfies the following conditions within a specified tolerance level $\epsilon^{k-1}$:
\begin{equation} \label{eq:3.6}
    \begin{aligned}
        & \mathrm{dist}\left(0, \; \partial_{\Theta} \caL_{\beta^{k-1}}\left( \Theta^k, \Lambda^{k-1} \right) \right) \le \epsilon^{k-1}, \\
        & Z_l^k \in 
        \begin{cases} 
            \mathrm{St}\left(n_{l-1}, n_{l}\right), & \text{if } l \in \mathcal{I},  \\ 
            \mathrm{St}\left(n_{l}, n_{l-1}\right), & \text{if } l \in \mathcal{I}^c, 
        \end{cases} \quad \forall l \in \{1,\dots,L\}.
    \end{aligned} 
\end{equation}
Here $\partial_\Theta\caL_{\beta^{k-1}}$ is the same product-manifold
stationarity mapping defined above, applied to the augmented Lagrangian.
For clarity, write $G_l^W,G_l^b,G_l^V$ for its three Euclidean gradient
blocks, $G_l^Z=\Pi_{Z_l^k}(\nabla_{Z_l}\caL_{\beta^{k-1}})$, and
$G_l^U=\partial_{U_l}\caL_{\beta^{k-1}}$, all evaluated at
$(\Theta^k,\Lambda^{k-1})$. Thus the stopping residual is exactly
\[
\begin{aligned}
&\mathrm{dist}\bigl(0,\partial_\Theta\caL_{\beta^{k-1}}
(\Theta^k,\Lambda^{k-1})\bigr)\\
&\quad=\left[\sum_{l=1}^L\left(
\|G_l^W\|_{\F}^2+\|G_l^Z\|_{\F}^2+\|G_l^b\|^2
+\|G_l^V\|_{\F}^2+\min_{g\in G_l^U}\|g\|_{\F}^2
\right)\right]^{1/2}.
\end{aligned}
\]
Following this inexact primal update and a subsequent reduction of the tolerance parameter, we update the Lagrange multipliers. Finally, the penalty parameter $\beta^k$ is adjusted according to the degree of constraint violation.

\begin{algorithm}
\caption{The Augmented Lagrangian Method for Solving
\eqref{prob.orth.three.sets.vari}}
\label{alg:alm}
\begin{algorithmic}[1]

\REQUIRE data matrices $X$ and $Y$;
initial feasible primal variables $\Theta^0\in\caR$ with $V_0=X$;
initial multipliers $\Lambda^0$;
initial penalty parameter $\beta^0>0$;
parameters $\eta_1,\eta_2,\eta_4\in(0,1)$ and $\eta_3>1$;
initial tolerance $\epsilon^0>0$;
initial constraint residual $\caR^0 = 0$.

\ENSURE an approximate solution to
\eqref{prob.orth.three.sets.vari}.

\STATE $k=1$;

\WHILE{the stopping criterion is not satisfied}

    \STATE solve subproblem \eqref{eq:3.5} to obtain
    $\Theta^k$ satisfying conditions \eqref{eq:3.6};

    \STATE update the subproblem tolerance as $\epsilon^k=\eta_4\epsilon^{k-1}$;

    \STATE update the multipliers for
    $a=1,2,3$ and $l=1,\dots,L$:
    \[
        \Lambda_{a,l}^k
        =\Lambda_{a,l}^{k-1}
        +\beta^{k-1}E_{a,l}^k;
    \]

    \STATE compute the constraint residual:
    \begin{equation}\label{orthdnn.cons.vio.r}
        \caR^k
        =\max_{\substack{a=1,2,3\\l=1,\dots,L}}
        \|E_{a,l}^k\|_{\F};
    \end{equation}

    \STATE update the penalty parameter:
    \begin{equation}\label{eq:3.8}
        \beta^k=
        \begin{cases}
            \beta^{k-1},
            & \text{if }\caR^k\le\eta_1\caR^{k-1},\\[1ex]
            \displaystyle
            \max\left\{
                \frac{\beta^{k-1}}{\eta_2},
                \max_{\substack{a=1,2,3\\ 1 \le l \le L}}
                \|\Lambda_{a,l}^k\|_{\F}^{1+\eta_3}
            \right\},
            & \text{otherwise};
        \end{cases}
    \end{equation}

    \STATE $k=k+1$;

\ENDWHILE

\RETURN $\Theta^{k-1}$.

\end{algorithmic}
\end{algorithm}

\subsubsection{The Alternating Minimization Method}\label{sec:alm.inner}

In this subsection, we propose an alternating minimization algorithm for solving the subproblem \eqref{eq:3.5} within Algorithm~\ref{alg:alm}. Assuming we are at the $k$-th iteration of Algorithm~\ref{alg:alm}, we drop the outer iteration index $k$ and adopt the following simplified notation for clarity:
\[
\begin{split}
    & \Theta^{(j)} = \l( \caW^{(j)}, \caZ^{(j)}, \caB^{(j)}, \caV^{(j)}, \caU^{(j)} \r) = \l( \caW^{k-1,j}, \caZ^{k-1,j}, \caB^{k-1,j}, \caV^{k-1,j}, \caU^{k-1,j} \r), \\
    & \{\lam_{i,l}\}_{l = 1}^L = \{\lam_{i,l}^{k-1}\}_{l = 1}^L, \quad i = 1,2,3, \qquad \beta = \beta^{k-1}, \quad \epsilon = \epsilon^{k-1}.
\end{split}
\] 
Before detailing the inner algorithmic framework, we first establish its initialization strategy. Let $\hat c$ be a constant upper bound such that 
\begin{equation}\label{gamma.cons.0}
    \hat c \geq \caL_{\beta}^0,
\end{equation}
where $\caL^0_{\beta} = \caL_{\beta^0} \left( \Theta^0, \lam^0 \right)$. 
The initial point for this inner alternating minimization algorithm is then defined as:
\begin{equation} \label{eq:3.13}
\begin{split}
     \Theta^{(0)}
      = \begin{cases}
        \Theta^{k-1}, & \text{if } k > 1 \text{ and } \caL_{\beta}^{k-1} \le \hat c, \\ 
        \Theta^0, & \text{otherwise},
    \end{cases}
\end{split}
\end{equation}
where $\caL_{\beta}^{k-1} := \caL_{\beta^{k-1}}(\Theta^{k-1}, \Lambda^{k-1})$. Note that since $\Theta^0 \in \caR$ is feasible, we have $\caL_{\beta^{k-1}}(\Theta^0, \Lambda^{k-1}) = \caF(\Theta^0) = \caL_{\beta}^0 \le \hat{c}$ for all $k \ge 1$.
Next, we denote 
\[
\hat{r}^{k-1} = \caL_\beta \l( \Theta^{(0)}, \lam \r), 
\]
and define the bounding parameter 
\begin{equation}\label{constant.l.max}
\begin{aligned}
    L_{\max}^{k-1} = \max \bigg\{ & \max_{1 \le l \le L} \l\{ L_{1,l}^{k-1} \r\} + \mu_1^2, \; \max_{1 \le l < L} \l\{ L_{1,l+1}^{k-1} + L_{2,l+1}^{k-1} \r\}, \\
    & \max_{1 \le l \le L} \l\{ L^{k-1}_{1,l} + L^{k-1}_{2,l} + L^{k-1}_{3,l} + L^{k-1}_{3,l-1} \r\} +\mu_2^2 \bigg\},
\end{aligned}
\end{equation}
where the constituent sequences $\l\{ L_{1,l}^{k-1} \r\}_{l=1}^L$, $\l\{ L_{2,l}^{k-1} \r\}_{l=1}^L$, and $\l\{ L_{3,l}^{k-1} \r\}_{l=1}^L$ are defined in \eqref{upper.bound.l.constants} and $L_{3,0}^{k-1} = 0$. 
When $L=1$, the empty middle maximum is defined to be zero.
This $L_{\max}^{k-1}$ plays a pivotal role in formulating the stopping criterion for the alternating minimization method.

In the following, we focus on the $j$-th iteration of this inner algorithm, explaining the resolution of each subproblem in detail and providing their corresponding analytical solutions.

{\bf The \(W_l\)-subproblem:}
We obtain $W_l^{(j)}$ ($l = 1,\dots,L$) by solving:
\begin{equation}\label{eq:3.14}
    \begin{aligned}
        W_l^{(j)} = \; & \arg\min_{W_l} \; \caL_{\beta} \Big( \{W_i^{(j)}\}_{i = 1}^{l-1}, W_l, \{W_i^{(j-1)}\}_{i = l+1}^L, \caZ^{(j-1)}, \\
        & \qquad \qquad \qquad \caB^{(j-1)}, \caV^{(j-1)}, \caU^{(j-1)}, \Lambda \Big) \\
        = \; & \arg\min_{W_l} \; \frac{\beta}{2}\left\Vert W_{l}V_{l-1}^{(j-1)} - \left( U_l^{(j-1)} - b_l^{(j-1)} 1_m^{\T} - \frac{1}{\beta}\Lambda_{1,l} \right) \right\Vert_{\F}^2 \\
        & \qquad \qquad + \frac{\beta}{2}\left\| W_l - \mathcal T_l^*\left( Z_l^{(j-1)} - \frac{1}{\beta}\Lambda_{3,l} \right) \right\|_{\F}^2.
    \end{aligned}
\end{equation}
Therefore, the closed-form solution of problem \eqref{eq:3.14} is 
\begin{equation}\label{w.l.sub.solu}
    \begin{aligned}
        W_l^{(j)} = \Bigg[ & \left( U_l^{(j-1)} - b_l^{(j-1)} 1_m^{\T} - \frac{1}{\beta}\Lambda_{1,l} \right) {V_{l-1}^{(j-1)}}^{\T} \\
        & + \mathcal T_l^*\left( Z_l^{(j-1)} - \frac{1}{\beta}\Lambda_{3,l} \right)
        \Bigg] \left( V_{l-1}^{(j-1)}{V_{l-1}^{(j-1)}}^{\T} + I_{n_{l-1}} \right)^{-1}.
    \end{aligned}
\end{equation}

\textbf{The $Z_l$-subproblem:}
To update $Z_l^{(j)}$ ($l = 1, \dots, L$), we introduce a proximal term and solve the following optimization problem:
\begin{equation}\label{eq:3.15}
\begin{aligned}
    \min_{Z_l} & \;\; \caL_{\beta} \Big(\caW^{(j)}, \{Z_i^{(j)}\}_{i = 1}^{l - 1}, Z_l, \{Z_i^{(j-1)}\}_{i = l+1}^L, \caB^{(j-1)}, \\
    & \qquad \qquad \caV^{(j-1)}, \caU^{(j-1)}, \lam \Big) + \frac{\mu_1}{2} \left\|Z_l - Z_l^{(j-1)} \right\|_{\F}^2 \\
    \text{s.t.} \;\; & Z_l \in \begin{cases} 
        \mathrm{St}\left(n_{l-1}, n_{l}\right), & \text{if } l \in \mathcal{I},  \\ 
        \mathrm{St}\left(n_{l}, n_{l-1}\right), & \text{if } l \in \mathcal{I}^c. 
    \end{cases} 
\end{aligned}
\end{equation}
Dropping the constant terms, the objective function in \eqref{eq:3.15} reduces to:
\[
    \frac{\beta}{2}\left\|Z_l-\left(\mathcal T_l(W_l^{(j)})+
    \frac{\Lambda_{3,l}}{\beta}\right)\right\|_{\F}^2
    +\frac{\mu_1}{2}\left\|Z_l-Z_l^{(j-1)}\right\|_{\F}^2.
\]

This problem is equivalent to computing a Euclidean projection of a matrix onto the Stiefel manifold. Thus, $Z_l^{(j)}$ is updated via the SVD as follows:
\begin{equation}\label{z.l.sub.solu}
    Z_l^{(j)} = \hat{U}_l \hat{V}_l^{\T},
\end{equation}
where $\hat{U}_l$ and $\hat{V}_l$ are obtained from the thin SVD:
\[
    \hat U_l\hat\Sigma_l\hat V_l^{\T}
    =\frac{\beta\mathcal T_l(W_l^{(j)})+\Lambda_{3,l}
    +\mu_1Z_l^{(j-1)}}{\beta+\mu_1}.
\]
The thin SVD retains all $\min\{n_l,n_{l-1}\}$ singular directions, including those with zero singular values. If the matrix is rank deficient, this construction selects one of the possibly nonunique minimizers.

{\bf The \(b_l\)-subproblem:}
We update $b_l^{(j)}$ ($l = 1, \dots, L$) by solving:  
\begin{equation}\label{eq:3.16}
\begin{split}
    b_{l}^{(j)} = & \; \arg\min_{b_l} \; \caL_{\beta} \Big(\caW^{(j)}, \caZ^{(j)}, \{b_i^{(j)}\}_{i = 1}^{l-1}, b_l, \{b_i^{(j-1)}\}_{i = l+1}^L, \caV^{(j-1)}, \caU^{(j-1)}, \lam \Big) \\
    = & \; \arg\min_{b_l} \; \frac{\lambda_1}{2} \|b_l\|^2 + \frac{\beta}{2} \left\Vert b_l 1_m^{\T} - \l( U_l^{(j-1)} - W_l^{(j)}V_{l-1}^{(j-1)} - \frac{\lam_{1,l}}{\beta} \r) \right\Vert_{\F}^2.
\end{split}
\end{equation}
The closed-form solution for \eqref{eq:3.16} is 
\begin{equation}\label{b.l.sub.solu}
b_l^{(j)} = \frac{\beta}{\lambda_1 + m \beta} \l( U_l^{(j-1)} - W_l^{(j)}V_{l-1}^{(j-1)} - \frac{\lam_{1,l}}{\beta} \r) 1_m.
\end{equation}

\textbf{The $V_l$-subproblem:} 
We obtain $V_l^{(j)}$ ($l = 1,\dots,L$) by solving:
\begin{equation}\label{eq:3.18}
\begin{aligned}
    V_l^{(j)} = \; & \arg\min_{V_l} \; \caL_{\beta} \Big(\caW^{(j)}, \caZ^{(j)}, \caB^{(j)},  \{V_i^{(j)}\}_{i = 1}^{l-1}, V_l, \{V_i^{(j-1)}\}_{i = l+1}^L, \caU^{(j-1)}, \lam \Big).
\end{aligned}
\end{equation}
For the intermediate layers $l \in \{1, \dots, L-1\}$, the objective function of \eqref{eq:3.18} is
\[
    \frac{\beta}{2} \left\Vert V_l - \left( \sigma(U_l^{(j-1)}) - \frac{\lam_{2,l}}{\beta} \right) \right\Vert_{\F}^2 + \frac{\beta}{2} \left\Vert W_{l+1}^{(j)} V_l - \left( U_{l+1}^{(j-1)} - b_{l+1}^{(j)} 1_m^{\T} - \frac{\lam_{1,l+1}}{\beta} \right) \right\Vert_{\F}^2.
\] 
The optimal solution is 
\begin{align}
    V_l^{(j)}
    ={} & \left(
        I + {W_{l+1}^{(j)}}^{\T} W_{l+1}^{(j)}
    \right)^{-1}
    \Biggl[
        \left(
            \sigma(U_l^{(j-1)}) - \frac{\lam_{2,l}}{\beta}
        \right)
    \nonumber\displaybreak[1]\\
    & \quad
    + {W_{l+1}^{(j)}}^{\T}
    \left(
        U_{l+1}^{(j-1)}
        - b_{l+1}^{(j)}1_m^{\T}
        - \frac{\lam_{1,l+1}}{\beta}
    \right)
    \Biggr].
    \label{v.l.sub.solu}
\end{align}
For the final layer $l = L$, the objective function for updating $V_L^{(j)}$ simplifies to:
\[
    \frac{1}{m} \left\Vert V_L - Y \right\Vert_{\F}^2 + \frac{\beta}{2} \left\Vert V_L - \left( \sigma( U_L^{(j-1)} ) - \frac{\lam_{2,L}}{\beta} \right) \right\Vert_{\F}^2.
\]
Similarly, this yields a closed-form update for $V_L^{(j)}$ given by:
\begin{equation}\label{v.fL.sub.solu}
   V_L^{(j)} = \frac{1}{\frac{2}{m} + \beta} \left[ \frac{2}{m}Y + \beta\sigma( U_L^{(j-1)} ) - \lam_{2,L} \right].
\end{equation}

\textbf{The $U_l$-subproblem:} 
To update $U_l^{(j)}$ ($l = 1,\dots,L$), again we introduce a proximal term and solve the following problem:
\begin{equation}\label{eq:3.17}
\begin{aligned}
    U_{l}^{(j)} \in & \; \arg\min_{U_l} \; \caL_{\beta} \Big(\caW^{(j)}, \caZ^{(j)}, \caB^{(j)}, \caV^{(j)}, \{U_i^{(j)}\}_{i = 1}^{l-1}, U_l, \\
    & \qquad \qquad \qquad \{U_i^{(j-1)}\}_{i = l+1}^L, \lam \Big) + \frac{\mu_2}{2} \left\| U_l - U_l^{(j-1)} \right\|_{\F}^2 \\
    = & \; \arg\min_{U_l} \; \frac{\lambda_2}{2} \| U_l \|_{\F}^2 + \frac{\beta}{2} \left\Vert U_l - \left( W_l^{(j)} V_{l-1}^{(j)} + b_l^{(j)} 1_m^{\T} + \frac{\lam_{1,l}}{\beta} \right) \right\Vert_{\F}^2 \\
    & \qquad \qquad \quad + \frac{\beta}{2}\left\Vert \sigma( U_l ) - \left( V_l^{(j)} + \frac{\lam_{2,l}}{\beta} \right) \right\Vert_{\F}^2 + \frac{\mu_2}{2} \left\| U_l - U_l^{(j-1)} \right\|_{\F}^2.
\end{aligned}
\end{equation}
Because the objective function in \eqref{eq:3.17} is completely separable with respect to the elements of $U_l$, solving it is equivalent to solving $n_l m$ independent one-dimensional problems. For brevity, let $u$, $c_1$, $c_2$, and $u'$ denote the corresponding scalar elements of the matrices $U_l$, $C_1$, $C_2$, and $U_l^{(j-1)}$ respectively, where $C_1 = W_l^{(j)} V_{l-1}^{(j)} + b_l^{(j)} 1_m^{\T} + \lam_{1,l}/\beta$ and $C_2 = V_l^{(j)} + \lam_{2,l}/\beta$. The element-wise optimization problem is given by:
\begin{equation}\label{eq2.4.4}
    \min_{u \in \mathbb{R}} \; g(u) = \frac{\lambda_2}{2} u^2 + \frac{\beta}{2} (u - c_1)^2 + \frac{\beta}{2} (\sigma(u) - c_2)^2 + \frac{\mu_2}{2} (u - u')^2.
\end{equation}
We define the optimal solutions over the nonnegative and nonpositive half-lines as:
\[
    u^+ := \arg\min_{u \in [0, +\infty)} \; g(u) \quad \text{and} \quad u^- := \arg\min_{u \in (-\infty, 0] } \; g(u).
\]

On each half-line, $g$ is a strictly convex quadratic. Setting its derivative to zero and projecting onto the corresponding half-line gives the following formulas. 
For leaky ReLU, the minimizers on the nonnegative and nonpositive half-lines are, respectively,
\begin{equation}\label{eq2.4.5}
    u^+=\max\left\{0,\frac{\beta(c_1+c_2)+\mu_2u'}{2\beta+\mu_2+\lambda_2}\right\}, \quad     u^-=\min\left\{0,\frac{\beta c_1+\alpha\beta c_2+\mu_2u'}
    {\beta+\alpha^2\beta+\mu_2+\lambda_2}\right\}.
\end{equation}
For ReLU, setting $\alpha=0$ leaves $u^+$ unchanged and gives
\begin{equation}\label{relueq2.4.5}
    u^-=\min\left\{
        0,\frac{\beta c_1+\mu_2u'}
        {\beta+\mu_2+\lambda_2}
    \right\}.
\end{equation}
Finally, the optimal element $u^{(j)}$ for the current iteration is determined by comparing the objective values:
\begin{equation}\label{eq2.4.7}
    u^{(j)} = \begin{cases}
        u^+, & \text{if } g(u^+) \leq g(u^-), \\ 
        u^-, & \text{otherwise.}
    \end{cases}
\end{equation}

For the inner iterations, define
\[
    \Delta_j A_l
    :=A_l^{(j)}-A_l^{(j-1)},
    \qquad
    A\in\{W,Z,b,V,U\},\quad
    l=1,\dots,L,\quad j\ge1.
\]
Since $V_0=X$ is fixed, we set $\Delta_jV_0=0$.
The framework of this alternating minimization method for solving \eqref{eq:3.5} is outlined in Algorithm~\ref{alg:bcd}. Analogous to the stopping criteria analyzed in \cite{wang2025augmented}, we will formally prove in Theorem~\ref{th.8.3.1} that satisfying condition \eqref{eq:3.19} guarantees condition \eqref{eq:3.6} is met.

\begin{algorithm}
\caption{The Alternating Minimization Method for Solving \eqref{eq:3.5}}
\label{alg:bcd}
\begin{algorithmic}[1]

\REQUIRE initial primal variables $\Theta^{(0)}$
chosen according to \eqref{eq:3.13};
multipliers $\Lambda$;
data matrices $X$ and $Y$;
penalty parameter $\beta>0$;
proximal parameters $\mu_1,\mu_2>0$;
tolerance $\epsilon>0$;
bounding parameter $L_{\max}^{k-1}$
defined in \eqref{constant.l.max}.

\ENSURE an approximate solution to \eqref{eq:3.5}.

\STATE $j = 1$;

\WHILE{the stopping criterion is not satisfied}

\FOR{$ l = 1, \ldots, L $}

    \STATE update $W_l^{(j)}$ by solving problem \eqref{eq:3.14};

\ENDFOR

\FOR{$ l = 1, \ldots, L $}

    \STATE update $Z_l^{(j)}$ by solving problem \eqref{eq:3.15};

\ENDFOR

\FOR{$ l = 1, \ldots, L $}

    \STATE update $b_l^{(j)}$ by solving problem \eqref{eq:3.16};

\ENDFOR

\FOR{$ l = 1, \ldots, L $}

    \STATE update $V_l^{(j)}$ by solving problem \eqref{eq:3.18};

\ENDFOR

\FOR{$ l = 1, \ldots, L $}

    \STATE update $U_l^{(j)}$ by solving problem \eqref{eq:3.17};

\ENDFOR

    \STATE if the stop criterion
    \begin{equation}\label{eq:3.19}
        \sum_{l=1}^{L}
        \Bigl(
            \|\Delta_j Z_l\|_{\F}^{2}
            +\|\Delta_j b_l\|^{2}
            +\|\Delta_j V_l\|_{\F}^{2}
            +\|\Delta_j U_l\|_{\F}^{2}
        \Bigr)
        \le
        \frac{\epsilon^{2}}{L_{\max}^{k-1}}.
    \end{equation}
    is not satisfied, then set $j = j + 1$;

\ENDWHILE

\STATE \textbf{return} $\Theta^{(j)}$.

\end{algorithmic}
\end{algorithm}

\subsection{OrthADMM for Problem \texorpdfstring{\eqref{ori.prob.orth.three.sets.vari}}{(2)}}\label{admm.al}

In this subsection, we present OrthADMM to solve the unregularized problem \eqref{ori.prob.orth.three.sets.vari}. From a methodological standpoint, OrthADMM serves as a practical, proximal-free ($\mu_1 = \mu_2 = 0$) variant of the OrthALM framework that executes only a single inner alternating step. This design trades the KKT convergence guarantees of our OrthALM for enhanced computational efficiency and reduced CPU time. The augmented Lagrangian function associated with \eqref{ori.prob.orth.three.sets.vari} is defined as:
\begin{equation}\label{eq.aug.lag.ori}
\widehat{\caL}_{\hat\beta}(\Theta,\hat\Lambda)
=\frac{1}{m}\|V_L-Y\|_{\F}^2
+\frac{\hat\beta}{2}\sum_{a=1}^3\sum_{l=1}^L
\left\|E_{a,l}(\Theta)+\frac{\hat\Lambda_{a,l}}{\hat\beta}\right\|_{\F}^2
-\frac{\|\hat\Lambda\|^2}{2\hat\beta}.
\end{equation}
where $\{ \hat{\Lambda}_{1,l} \in \mathbb{R}^{n_l \times m} \}_{l = 1}^L $, $ \{ \hat{\Lambda}_{2,l} \in \mathbb{R}^{n_l \times m} \}_{l = 1}^L$, $ \{ \hat{\Lambda}_{3,l} \in \mathbb{R}^{n_l \times n_{l-1}} \}_{l \in \mathcal{I}}$ and $ \{ \hat{\Lambda}_{3,l} \in \mathbb{R}^{n_{l-1} \times n_l} \}_{l \in \mathcal{I}^c}$ denote the multiplier matrices, and $\hat{\beta} > 0$ is the penalty parameter. In the $k$-th iteration of OrthADMM, we sequentially update the primal variables for $l = 1,\dots,L$ by solving the following subproblems: 
\begin{subequations}\label{nn.admm.prime.sub}
\begin{align}
    W_l^k = \; & \arg\min_{W_l} \; \hat{\caL}_{\hat{\beta}} \Big( \{W_i^k\}_{i = 1}^{l-1}, W_l, \{W_i^{k-1}\}_{i = l+1}^L, \caZ^{k-1}, \caB^{k-1}, \caV^{k-1}, \caU^{k-1}, \hat{\lam}^{k-1} \Big), \label{nn.admm.prime.sub.w} \\
    Z_l^k \in \; & \arg\min_{Z_l} \; \hat{\caL}_{\hat{\beta}} \Big( \caW^k, \{Z_i^k\}_{i = 1}^{l-1}, Z_l, \{Z_i^{k-1}\}_{i = l+1}^L, \caB^{k-1}, \caV^{k-1}, \caU^{k-1}, \hat{\lam}^{k-1} \Big), \label{nn.admm.prime.sub.z} \\
    b_l^k = \; & \arg\min_{b_l} \; \hat{\caL}_{\hat{\beta}} \Big( \caW^k, \caZ^k, \{b_i^k\}_{i = 1}^{l-1}, b_l, \{b_i^{k-1}\}_{i = l+1}^L, \caV^{k-1}, \caU^{k-1}, \hat{\lam}^{k-1} \Big), \label{nn.admm.prime.sub.b} \\
    V_l^k = \; & \arg\min_{V_l} \; \hat{\caL}_{\hat{\beta}} \Big( \caW^k, \caZ^k, \caB^k, \{V_i^k\}_{i = 1}^{l-1}, V_l, \{V_i^{k-1}\}_{i = l+1}^L, \caU^{k-1}, \hat{\lam}^{k-1} \Big), \label{nn.admm.prime.sub.v} \\
    U_l^k \in \; & \arg\min_{U_l} \; \hat{\caL}_{\hat{\beta}} \Big( \caW^k, \caZ^k, \caB^k, \caV^{k}, \{U_i^k\}_{i = 1}^{l-1}, U_l, \{U_i^{k-1}\}_{i = l+1}^L, \hat{\lam}^{k-1} \Big), \label{nn.admm.prime.sub.u}
\end{align}
\end{subequations}
where the minimization with respect to $Z_l$ in \eqref{nn.admm.prime.sub.z} is subject to $Z_l \in \STM{n_{l-1}}{n_l}$ if $l \in \mathcal{I}$ and $Z_l \in \STM{n_l}{n_{l-1}}$ if $l \in \mathcal{I}^c$.
Because problem \eqref{ori.prob.orth.three.sets.vari} is the unregularized counterpart to problem \eqref{prob.orth.three.sets.vari}, the explicit closed-form solutions for these subproblems are similar to those derived in the last subsection, naturally obtained by setting the regularization and proximal parameters to zero ($\lambda_1 = \lambda_2 = \mu_1 = \mu_2 = 0$). 

The complete OrthADMM procedure, including the dual variable updates, is outlined in Algorithm~\ref{alg:admm.nn}.

\begin{algorithm}[htbp]
\caption{OrthADMM for Solving
\eqref{ori.prob.orth.three.sets.vari}}
\label{alg:admm.nn}
\begin{algorithmic}[1]

\REQUIRE data matrices $X$ and $Y$;
initial feasible primal variables $\Theta^0\in\caR$ with $V_0=X$;
initial multipliers $\hat{\Lambda}^0$;
penalty parameter $\hat{\beta}>0$.

\ENSURE an approximate solution to
\eqref{ori.prob.orth.three.sets.vari}.

\STATE $k=1$;

\WHILE{the stopping criterion is not satisfied}

    \FOR{$l=1,\dots,L$}
        \STATE update $W_l^k$ by solving problem 
        \eqref{nn.admm.prime.sub.w};
    \ENDFOR

    \FOR{$l=1,\dots,L$}
        \STATE update $Z_l^k$ by solving problem 
        \eqref{nn.admm.prime.sub.z};
    \ENDFOR

    \FOR{$l=1,\dots,L$}
        \STATE update $b_l^k$ by solving problem 
        \eqref{nn.admm.prime.sub.b};
    \ENDFOR

    \FOR{$l=1,\dots,L$}
        \STATE update $V_l^k$ by solving problem 
        \eqref{nn.admm.prime.sub.v};
    \ENDFOR

    \FOR{$l=1,\dots,L$}
        \STATE update $U_l^k$ by solving problem 
        \eqref{nn.admm.prime.sub.u};
    \ENDFOR

    \STATE update the multipliers for
    $a=1,2,3$ and $l=1,\dots,L$:
    \[
        \hat{\Lambda}_{a,l}^k
        =\hat{\Lambda}_{a,l}^{k-1}
        +\hat{\beta}E_{a,l}^k;
    \]

    \STATE $k=k+1$;

\ENDWHILE

\RETURN $\Theta^{k-1}$.

\end{algorithmic}
\end{algorithm}

\section{Convergence Analysis}\label{sec:convergence}

This section is dedicated to the convergence analysis of our proposed OrthALM method (Algorithm~\ref{alg:alm} and Algorithm~\ref{alg:bcd}) for solving the regularized problem \eqref{prob.orth.three.sets.vari}.

\subsection{Well-definedness of Step 3 of Algorithm~\ref{alg:alm}}\label{subsec:welldefined}

This subsection establishes that step~3 of Algorithm~\ref{alg:alm} is well-defined when utilizing Algorithm~\ref{alg:bcd}. Specifically, we prove that Algorithm~\ref{alg:bcd} reaches its stopping criterion in finitely many iterations and returns variables $\Theta^k$ satisfying \eqref{eq:3.6}.
We begin by analyzing the sufficient descent properties of Algorithm~\ref{alg:bcd}, which is to quantify the objective progress achieved during a single iteration.

For a fixed outer iteration, write $\beta=\beta^{k-1}$ and $\Lambda=\Lambda^{k-1}$.
Let $\Theta_{A_l}^{(j)}$ denote the full tuple immediately after updating block $A_l$, where $A\in\{W,Z,b,V,U\}$. For example,
\begin{align*}
\Theta_{W_l}^{(j)}=\l( (W_1^{(j)},\ldots,W_l^{(j)},W_{l+1}^{(j-1)},\ldots,W_L^{(j-1)}), \caZ^{(j-1)},\caB^{(j-1)},\caV^{(j-1)},\caU^{(j-1)}\r).
\end{align*}
The boundary states are
\begin{align*}
\Theta_{W_0}^{(j)}&=\Theta^{(j-1)},&
\Theta_{Z_0}^{(j)}&=\Theta_{W_L}^{(j)},&
\Theta_{b_0}^{(j)}&=\Theta_{Z_L}^{(j)},\\
\Theta_{V_0}^{(j)}&=\Theta_{b_L}^{(j)},&
\Theta_{U_0}^{(j)}&=\Theta_{V_L}^{(j)},&
\Theta_{U_L}^{(j)}&=\Theta^{(j)}.
\end{align*}
The gradients with respect to $W_l$, $b_l$, and $V_l$ are unchanged by subsequent updates of other blocks in the same family. Thus their optimality conditions can also be evaluated at $\Theta_{W_L}^{(j)}$, $\Theta_{b_L}^{(j)}$, and $\Theta_{V_L}^{(j)}$, respectively.

\begin{itemize}
\item \textbf{The $W_l$-subproblem:} The optimality condition is:
\begin{align*}
\nabla_{W_l}\caL_\beta(\Theta_{W_l}^{(j)},\Lambda)
={}&\big[\Lambda_{1,l}+\beta E_{1,l}(\Theta_{W_l}^{(j)})\big]{V_{l-1}^{(j-1)}}^{\T}\\
&+\mathcal T_l^*\big[\Lambda_{3,l}+\beta E_{3,l}(\Theta_{W_l}^{(j)})\big]=0.
\end{align*}
Expanding the quadratic terms and using the first-order
optimality condition of the $W_l$-subproblem, we obtain
\begin{align*}
&\caL_\beta(\Theta_{W_l}^{(j)},\Lambda)-\caL_\beta(\Theta_{W_{l-1}}^{(j)},\Lambda)\\
={}&\frac\beta2\big[\|E_{1,l}(\Theta_{W_l}^{(j)})\|_{\F}^2
-\|E_{1,l}(\Theta_{W_{l-1}}^{(j)})\|_{\F}^2\big]
+\langle\Lambda_{1,l},\Delta_jW_lV_{l-1}^{(j-1)}\rangle\\
&+\frac\beta2\big[\|E_{3,l}(\Theta_{W_l}^{(j)})\|_{\F}^2
-\|E_{3,l}(\Theta_{W_{l-1}}^{(j)})\|_{\F}^2\big]
+\langle\Lambda_{3,l},\mathcal T_l(\Delta_jW_l)\rangle\\
={}&\big\langle\Lambda_{1,l}+\beta E_{1,l}(\Theta_{W_l}^{(j)}),\Delta_jW_lV_{l-1}^{(j-1)}\big\rangle
+\big\langle\Lambda_{3,l}+\beta E_{3,l}(\Theta_{W_l}^{(j)}),\mathcal T_l(\Delta_jW_l)\big\rangle\\
&-\frac\beta2\|\Delta_jW_lV_{l-1}^{(j-1)}\|_{\F}^2-\frac\beta2\|\Delta_jW_l\|_{\F}^2\\
={}&\langle\nabla_{W_l}\caL_\beta(\Theta_{W_l}^{(j)},\Lambda),\Delta_jW_l\rangle
-\frac\beta2\|\Delta_jW_lV_{l-1}^{(j-1)}\|_{\F}^2-\frac\beta2\|\Delta_jW_l\|_{\F}^2\\
\le{}&-\frac\beta2\|\Delta_jW_l\|_{\F}^2.
\end{align*}
Here $\mathcal T_l$ preserves the Frobenius norm, for both orientations of $W_l$.

\item \textbf{The $Z_l$-subproblem:} Global optimality in \eqref{eq:3.15}, compared with the previous block value, gives
\[
\caL_\beta(\Theta_{Z_l}^{(j)},\Lambda)-\caL_\beta(\Theta_{Z_{l-1}}^{(j)},\Lambda)
\le-\frac{\mu_1}{2}\|\Delta_jZ_l\|_{\F}^2.
\]

\item \textbf{The $b_l$-subproblem:} Its optimality condition is
\[
\nabla_{b_l}\caL_\beta(\Theta_{b_l}^{(j)},\Lambda)
=\lambda_1b_l^{(j)}+\big[\Lambda_{1,l}+\beta E_{1,l}(\Theta_{b_l}^{(j)})\big]1_m=0.
\]
The corresponding descent calculation is
\begin{align*}
&\caL_\beta(\Theta_{b_l}^{(j)},\Lambda)-\caL_\beta(\Theta_{b_{l-1}}^{(j)},\Lambda)\\
={}&\frac{\lambda_1}{2}\big(\|b_l^{(j)}\|^2-\|b_l^{(j-1)}\|^2\big)
+\langle\Lambda_{1,l},\Delta_jb_l1_m^{\T}\rangle + \frac\beta2\big[\|E_{1,l}(\Theta_{b_l}^{(j)})\|_{\F}^2
-\|E_{1,l}(\Theta_{b_{l-1}}^{(j)})\|_{\F}^2\big]\\
={}&\big\langle\Lambda_{1,l}+\beta E_{1,l}(\Theta_{b_l}^{(j)}),\Delta_jb_l1_m^{\T}\big\rangle
-\frac\beta2\|\Delta_jb_l1_m^{\T}\|_{\F}^2
+\frac{\lambda_1}{2}\big(\|b_l^{(j)}\|^2-\|b_l^{(j-1)}\|^2\big)\\
={}&\langle\nabla_{b_l}\caL_\beta(\Theta_{b_l}^{(j)},\Lambda),\Delta_jb_l\rangle
-\frac\beta2\|\Delta_jb_l1_m^{\T}\|_{\F}^2\\
&+\frac{\lambda_1}{2}\big(\|b_l^{(j)}\|^2-\|b_l^{(j-1)}\|^2\big)
-\langle\lambda_1b_l^{(j)},\Delta_jb_l\rangle\\
={}&-\frac\beta2\|\Delta_jb_l1_m^{\T}\|_{\F}^2-\frac{\lambda_1}{2}\|\Delta_jb_l\|^2
=-\frac{m\beta+\lambda_1}{2}\|\Delta_jb_l\|^2.
\end{align*}

\item \textbf{The $V_l$-subproblem:} For $l<L$, the optimality condition is
\begin{align*}
\nabla_{V_l}\caL_\beta(\Theta_{V_l}^{(j)},\Lambda)
={}&{W_{l+1}^{(j)}}^{\T}\big[\Lambda_{1,l+1}+\beta E_{1,l+1}(\Theta_{V_l}^{(j)})\big] + \Lambda_{2,l}+\beta E_{2,l}(\Theta_{V_l}^{(j)})=0.
\end{align*}
The corresponding descent calculation is
\begin{align*}
&\caL_\beta(\Theta_{V_l}^{(j)},\Lambda)-\caL_\beta(\Theta_{V_{l-1}}^{(j)},\Lambda)\\
={}&\frac\beta2\big[\|E_{1,l+1}(\Theta_{V_l}^{(j)})\|_{\F}^2
-\|E_{1,l+1}(\Theta_{V_{l-1}}^{(j)})\|_{\F}^2\big]\\
&+\frac\beta2\big[\|E_{2,l}(\Theta_{V_l}^{(j)})\|_{\F}^2
-\|E_{2,l}(\Theta_{V_{l-1}}^{(j)})\|_{\F}^2\big]
+\langle\Lambda_{1,l+1},W_{l+1}^{(j)}\Delta_jV_l\rangle+\langle\Lambda_{2,l},\Delta_jV_l\rangle\\
={}&\big\langle\Lambda_{1,l+1}+\beta E_{1,l+1}(\Theta_{V_l}^{(j)}),W_{l+1}^{(j)}\Delta_jV_l\big\rangle
+\big\langle\Lambda_{2,l}+\beta E_{2,l}(\Theta_{V_l}^{(j)}),\Delta_jV_l\big\rangle\\
&-\frac\beta2\|W_{l+1}^{(j)}\Delta_jV_l\|_{\F}^2-\frac\beta2\|\Delta_jV_l\|_{\F}^2\\
={}&-\frac\beta2\|W_{l+1}^{(j)}\Delta_jV_l\|_{\F}^2-\frac\beta2\|\Delta_jV_l\|_{\F}^2
\le-\frac\beta2\|\Delta_jV_l\|_{\F}^2.
\end{align*}
For $l=L$, we have
\[
\nabla_{V_L}\caL_\beta(\Theta_{V_L}^{(j)},\Lambda)
=\frac2m(V_L^{(j)}-Y)+\Lambda_{2,L}+\beta E_{2,L}(\Theta_{V_L}^{(j)})=0,
\]
and hence
\begin{align*}
&\caL_\beta(\Theta_{V_L}^{(j)},\Lambda)-\caL_\beta(\Theta_{V_{L-1}}^{(j)},\Lambda)\\
={}&\frac1m\big(\|V_L^{(j)}-Y\|_{\F}^2-\|V_L^{(j-1)}-Y\|_{\F}^2\big)
+\langle\Lambda_{2,L},\Delta_jV_L\rangle\\
&+\frac\beta2\big[\|E_{2,L}(\Theta_{V_L}^{(j)})\|_{\F}^2
-\|E_{2,L}(\Theta_{V_{L-1}}^{(j)})\|_{\F}^2\big]\\
={}&\big\langle\Lambda_{2,L}+\beta E_{2,L}(\Theta_{V_L}^{(j)}),\Delta_jV_L\big\rangle
+\frac2m\langle V_L^{(j)}-Y,\Delta_jV_L\rangle
-\left(\frac\beta2+\frac1m\right)\|\Delta_jV_L\|_{\F}^2\\
={}&-\left(\frac\beta2+\frac1m\right)\|\Delta_jV_L\|_{\F}^2.
\end{align*}

\item \textbf{The $U_l$-subproblem:} Global optimality in \eqref{eq:3.17} yields
\[
\caL_\beta(\Theta_{U_l}^{(j)},\Lambda)-\caL_\beta(\Theta_{U_{l-1}}^{(j)},\Lambda)
\le-\frac{\mu_2}{2}\|\Delta_jU_l\|_{\F}^2.
\]
\end{itemize}

We now derive upper bounds for the partial-gradient differences of the augmented Lagrangian. These bounds are essential for formulating the algorithmic stopping criterion and conducting the subsequent convergence analysis.

\begin{lemma}\label{le.partial.bound}

Fix an outer iteration $k$ and set $\beta=\beta^{k-1}$ and $\Lambda=\Lambda^{k-1}$. For $j\ge1$, suppose that the full tuples
\[
\Theta^{(j-1)},\quad \Theta_{W_L}^{(j)},\quad
\Theta_{b_L}^{(j)},\quad \Theta_{V_L}^{(j)},\quad
\Theta^{(j)}
\]
belong to $\Omega_{\caL_\beta}(\hat r^{k-1})$.
Then the following partial-gradient difference bounds hold. 
    \begin{enumerate}
    \item For $W_l$:
\begin{equation}\label{upper.bound.w.l}
\begin{aligned}
    & \left\| \nabla_{W_l} \caL_\beta \l( \Theta_{W_L}^{(j)}, \Lambda \r) - \nabla_{W_l} \caL_\beta \l( \Theta^{(j)}, \Lambda \r) \right\|_{\F}^2 \\
    \le \; & \begin{cases}
        \begin{aligned}
            & L_{1,1}^{k-1} \left( \|\Delta_j Z_{1}\|^2_{\F} + \|\Delta_j b_{1}\|^2 +  \|\Delta_j U_{1}\|^2_{\F} \right), 
        \end{aligned} & \text{for } l = 1, \\[4ex]
        \begin{aligned}
            & L_{1,l}^{k-1} \left( \|\Delta_j V_{l-1}\|^2_{\F} + \|\Delta_j Z_{l}\|^2_{\F} + \|\Delta_j b_{l}\|^2 + \|\Delta_j U_{l}\|^2_{\F} \right),
        \end{aligned} & \text{for } l = 2,\dots,L.
    \end{cases}
\end{aligned}
\end{equation}

    \item For $b_l$:
    \begin{equation}\label{upper.bound.b.l}
    \begin{aligned}
       & \left\| \nabla_{b_l}\caL_\beta \l( \Theta_{b_L}^{(j)}, \Lambda \r) - \nabla_{b_l} \caL_\beta \l( \Theta^{(j)}, \Lambda \r) \right\|^2 \\
         \le & \begin{cases}
            L_{2,1}^{k-1} \|\Delta_j U_{1}\|^2_{\F}, & \text{for } l = 1, \\
            L_{2,l}^{k-1} \left( \|\Delta_j V_{l-1}\|^2_{\F} + \|\Delta_j U_{l}\|^2_{\F} \right), & \text{for } l = 2,\dots,L.
        \end{cases}
    \end{aligned}
    \end{equation}
    
    \item For $V_l$:
    \begin{equation}\label{upper.bound.v.l}
    \begin{aligned}
       & \left\| \nabla_{V_l}\caL_\beta \l( \Theta_{V_L}^{(j)}, \Lambda \r) - \nabla_{V_l} \caL_\beta \l( \Theta^{(j)}, \Lambda \r) \right\|_{\F}^2 \\
         \le & \begin{cases}
            L_{3,l}^{k-1} \l( \|\Delta_j U_{l+1}\|^2_{\F} + \|\Delta_j U_{l}\|^2_{\F} \r), & \text{for } l = 1,...,L-1, \\
            L_{3,L}^{k-1} \|\Delta_j U_{L}\|^2_{\F}, & \text{for } l = L.
        \end{cases}                                
    \end{aligned}
    \end{equation}

    \end{enumerate}
    \noindent Here, the Lipschitz-like constants are defined as:
    \begin{equation}\label{upper.bound.l.constants}
    \begin{aligned}
        L_{1,l}^{k-1} & = L_{1,l} \l( \lam, \beta, \hat{r}^{k-1} \r) = \begin{cases}
            3 \beta^2 \max \l\{ m\|X\|^2_{\F}, \, 1 \r\}, & \text{for } l = 1, \\
            8 \max \l\{ \delta^{k-1}_{2,l}, \, \beta^2, \, m\beta^2\delta_{1,l-1}^{k-1} \r\}, & \text{for } l = 2,\dots,L,
        \end{cases} \\
        L_{2,l}^{k-1} & = L_{2,l} \l( \lam, \beta, \hat{r}^{k-1} \r) = \begin{cases}
             m \beta^2, & \text{for } l = 1, \\
             2 m \beta^2 \max \l\{ \delta_{0,l}^{k-1}, \, 1 \r\}, & \text{for } l = 2,\dots,L,
         \end{cases} \\
         L_{3,l}^{k-1} & = L_{3,l} \l( \lam, \beta, \hat{r}^{k-1} \r) = \begin{cases}
             2 \beta^2 \max \l\{ \delta_{0,l+1}^{k-1}, \, 1 \r\}, & \text{for } l = 1,...,L-1, \\
             \beta^2, & \text{for } l = L,
         \end{cases}
    \end{aligned}
    \end{equation}
    where
    \begin{align*}
        \delta^{k-1} &= \hat{r}^{k-1} + \frac{1}{2\beta}\sum_{l=1}^L \left( \|\Lambda_{1,l}\|_{\F}^2 + \|\Lambda_{2,l}\|_{\F}^2 + \|\Lambda_{3,l}\|_{\F}^2 \right), \quad n_{\max} = \max \{ n_0,\dots, n_L \}, \\
        \delta_{0,l}^{k-1} &= \left( \sqrt{\frac{2\delta^{k-1}}{\beta}} + \sqrt{n_{\max}} + \frac{\|\Lambda_{3,l}\|_{\F}}{\beta} \right)^2,  \quad {\rm for} \; l = 1,\dots,L, \\
        \delta_{1,l}^{k-1} &= \left( \sqrt{\frac{2\delta^{k-1}}{\beta}} + \sqrt{\frac{2\delta^{k-1}}{\lambda_2}} + \frac{\|\Lambda_{2,l}\|_{\F}}{\beta} \right)^2, \quad {\rm for} \; l = 1,\dots,L,\\
        \delta_{2,l}^{k-1} &= \|\Lambda_{1,l}\|^2_{\F} + \beta^2 \left( 2\delta_{0,l}^{k-1} \delta_{1,l-1}^{k-1} + \frac{2m\delta^{k-1}}{\lambda_1} + \frac{2\delta^{k-1}}{\lambda_2} \right) \quad {\rm for} \; l = 2,\dots,L. 
    \end{align*}
\end{lemma}

\begin{proof}
    The $W_l$ gradient is
\[
\nabla_{W_l}\caL_\beta(\Theta,\Lambda)
=\big[\Lambda_{1,l}+\beta E_{1,l}(\Theta)\big]V_{l-1}^{\T}
+\mathcal T_l^*\big[\Lambda_{3,l}+\beta E_{3,l}(\Theta)\big].
\]
For $l>1$, adding and subtracting the mixed product
$E_{1,l}(\Theta_{W_L}^{(j)}){V_{l-1}^{(j)}}^{\T}$ gives
\begin{align*}
&\nabla_{W_l}\caL_\beta(\Theta_{W_L}^{(j)},\Lambda)
-\nabla_{W_l}\caL_\beta(\Theta^{(j)},\Lambda)\\
={}&\Lambda_{1,l}{V_{l-1}^{(j-1)}}^{\T}
+\beta E_{1,l}(\Theta_{W_L}^{(j)}){V_{l-1}^{(j-1)}}^{\T}
-\Lambda_{1,l}{V_{l-1}^{(j)}}^{\T}
-\beta E_{1,l}(\Theta^{(j)}){V_{l-1}^{(j)}}^{\T}\\
&+\mathcal T_l^*\big[\Lambda_{3,l}+\beta E_{3,l}(\Theta_{W_L}^{(j)})
-\Lambda_{3,l}-\beta E_{3,l}(\Theta^{(j)})\big]\\
={}&-\Lambda_{1,l}(\Delta_jV_{l-1})^{\T}+\beta\mathcal T_l^*(\Delta_jZ_l) + \beta\big[E_{1,l}(\Theta_{W_L}^{(j)}){V_{l-1}^{(j-1)}}^{\T} \\
& - E_{1,l}(\Theta_{W_L}^{(j)}){V_{l-1}^{(j)}}^{\T} + E_{1,l}(\Theta_{W_L}^{(j)}){V_{l-1}^{(j)}}^{\T} - E_{1,l}(\Theta^{(j)}){V_{l-1}^{(j)}}^{\T}\big]\\
={}&-\Lambda_{1,l}(\Delta_jV_{l-1})^{\T}+\beta\mathcal T_l^*(\Delta_jZ_l)\\
&+\beta\big[-W_l^{(j)}V_{l-1}^{(j-1)}(\Delta_jV_{l-1})^{\T}
-b_l^{(j-1)}1_m^{\T}(\Delta_jV_{l-1})^{\T}
+U_l^{(j-1)}(\Delta_jV_{l-1})^{\T}\\
&\hspace{2.8em}-W_l^{(j)}\Delta_jV_{l-1}{V_{l-1}^{(j)}}^{\T}
-\Delta_jb_l1_m^{\T}{V_{l-1}^{(j)}}^{\T}
+\Delta_jU_l{V_{l-1}^{(j)}}^{\T}\big].
\end{align*}
For $l=1$, $V_0^{(j)}=V_0^{(j-1)}=X$, and therefore
\begin{align*}
\nabla_{W_1}\caL_\beta(\Theta_{W_L}^{(j)},\Lambda)
-\nabla_{W_1}\caL_\beta(\Theta^{(j)},\Lambda)
= \beta\big[(-\Delta_jb_1 1_m^{\T}+\Delta_jU_1)X^{\T}
+\mathcal T_1^*(\Delta_jZ_1)\big].
\end{align*}
Every full tuple specified in the lemma belongs to
$\Omega_{\caL_\beta}(\hat r^{k-1})$. Completing the squares in the augmented Lagrangian consequently gives, at each such tuple,
\[
\|b_l\|^2\le\frac{2\delta^{k-1}}{\lambda_1},\qquad
\|U_l\|_{\F}^2\le\frac{2\delta^{k-1}}{\lambda_2},\qquad
\|Z_l\|_{\F}^2=\min\{n_l,n_{l-1}\}\le n_{\max}.
\]
The reverse triangle inequality also yields
\begin{align*}
\|W_l\|_{\F} - \l\|Z_l - \frac{\lam_{3,l}}{\beta} \r\|_{\F}
&\le \l\| \mathcal T_l(W_l)-Z_l + \frac{\lam_{3,l}}{\beta} \r\|_{\F}
\le\sqrt{\frac{2\delta^{k-1}}{\beta}},\\
\|V_l\|_{\F}- \l\|\sigma(U_l) - \frac{\lam_{2,l}}{\beta} \r\|_{\F}
&\le \l\| V_l-\sigma(U_l) + \frac{\lam_{2,l}}{\beta} \r\|_{\F}
\le\sqrt{\frac{2\delta^{k-1}}{\beta}},\\
\|\sigma(U_l)\|_{\F}&\le\|U_l\|_{\F}.
\end{align*}
Hence
\begin{align*}
\|W_l\|_{\F}^2
&\le\left(\sqrt{\frac{2\delta^{k-1}}\beta}+\sqrt{n_{\max}}
+\frac{\|\Lambda_{3,l}\|_{\F}}\beta\right)^2=\delta_{0,l}^{k-1},\\
\|V_l\|_{\F}^2
&\le\left(\sqrt{\frac{2\delta^{k-1}}\beta}+\sqrt{\frac{2\delta^{k-1}}{\lambda_2}}
+\frac{\|\Lambda_{2,l}\|_{\F}}\beta\right)^2=\delta_{1,l}^{k-1}.
\end{align*}
In particular, $\delta_{0,l}^{k-1}\ge1$. Applying
$\|\sum_{i=1}^{8}A_i\|_{\F}^2\le8\sum_{i=1}^{8}\|A_i\|_{\F}^2$
to the eight terms of the gradient difference, we obtain for $l>1$
\begin{align*}
&\|\nabla_{W_l}\caL_\beta(\Theta_{W_L}^{(j)},\Lambda)
-\nabla_{W_l}\caL_\beta(\Theta^{(j)},\Lambda)\|_{\F}^2\\
\le{}&8\Big\{\Big[\|\Lambda_{1,l}\|_{\F}^2
+\beta^2\big(\|W_l^{(j)}\|_{\F}^2\|V_{l-1}^{(j-1)}\|_{\F}^2
+m\|b_l^{(j-1)}\|^2+\|U_l^{(j-1)}\|_{\F}^2\\
&\hspace{7em}+\|W_l^{(j)}\|_{\F}^2\|V_{l-1}^{(j)}\|_{\F}^2\big)\Big]
\|\Delta_jV_{l-1}\|_{\F}^2+\beta^2\|\Delta_jZ_l\|_{\F}^2\\
&\hspace{2em}+m\beta^2\|V_{l-1}^{(j)}\|_{\F}^2\|\Delta_jb_l\|^2
+\beta^2\|V_{l-1}^{(j)}\|_{\F}^2\|\Delta_jU_l\|_{\F}^2\Big\}\\
\le{}&8\Big[\delta_{2,l}^{k-1}\|\Delta_jV_{l-1}\|_{\F}^2
+\beta^2\|\Delta_jZ_l\|_{\F}^2
+m\beta^2\delta_{1,l-1}^{k-1}\|\Delta_jb_l\|^2
+\beta^2\delta_{1,l-1}^{k-1}\|\Delta_jU_l\|_{\F}^2\Big]\\
\le{}&8\max\{\delta_{2,l}^{k-1},\beta^2,m\beta^2\delta_{1,l-1}^{k-1}\}
\big(\|\Delta_jV_{l-1}\|_{\F}^2+\|\Delta_jZ_l\|_{\F}^2
+\|\Delta_jb_l\|^2+\|\Delta_jU_l\|_{\F}^2\big).
\end{align*}
Here the coefficient of $\|\Delta_jV_{l-1}\|_{\F}^2$ is bounded by
\[
\delta_{2,l}^{k-1}=\|\Lambda_{1,l}\|_{\F}^2
+\beta^2\left(2\delta_{0,l}^{k-1}\delta_{1,l-1}^{k-1}
+\frac{2m\delta^{k-1}}{\lambda_1}+\frac{2\delta^{k-1}}{\lambda_2}\right),
\]
as defined above. For $l=1$, applying the corresponding three-term bound gives
\begin{align*}
&\|\nabla_{W_1}\caL_\beta(\Theta_{W_L}^{(j)},\Lambda)
-\nabla_{W_1}\caL_\beta(\Theta^{(j)},\Lambda)\|_{\F}^2\\
\le{}&3\beta^2\big(\|\Delta_jZ_1\|_{\F}^2
+m\|X\|_{\F}^2\|\Delta_jb_1\|^2
+\|X\|_{\F}^2\|\Delta_jU_1\|_{\F}^2\big)\\
\le{}&3\beta^2\max\{m\|X\|_{\F}^2,1\}
\big(\|\Delta_jZ_1\|_{\F}^2+\|\Delta_jb_1\|^2+\|\Delta_jU_1\|_{\F}^2\big).
\end{align*}
These inequalities prove \eqref{upper.bound.w.l} with the constants in \eqref{upper.bound.l.constants}.

    For $b_l$ ($l = 1, \dots, L$), we know that the gradient of $\caL_\beta$ with respect to $b_l$ is:
    \begin{align*}
        & \nabla_{b_l}\caL_\beta(\Theta, \Lambda) = \lambda_1 b_l + \Lambda_{1,l} 1_m + \beta\left( b_l 1_m^{\T} + W_l V_{l-1} - U_l \right) 1_m. 
    \end{align*}
    Therefore, for $l > 1$, the difference is:
    \begin{align*}
        & \nabla_{b_l} \caL_\beta \l( \Theta_{b_L}^{(j)}, \Lambda \r) - \nabla_{b_l} \caL_\beta \l( \Theta^{(j)}, \Lambda \r) \\
        = \; & \lambda_1 b_l^{(j)} + \Lambda_{1,l} 1_m + \beta \left( b_l^{(j)} 1_m^{\T} + W_l^{(j)} V^{(j-1)}_{l-1} - U^{(j-1)}_l \right) 1_m - \lambda_1 b_l^{(j)} - \Lambda_{1,l} 1_m \\
        & - \beta \left( b_l^{(j)} 1_m^{\T} + W_l^{(j)} V^{(j)}_{l-1} - U^{(j)}_l \right) 1_m \\
        = \; & \beta \left[ W_l^{(j)} \l( -\Delta_j V_{l-1} \r) + \l( \Delta_j U_{l} \r)\right]1_m.
    \end{align*}
    For $l = 1$, this becomes:
    \begin{align*}
        & \nabla_{b_1} \caL_\beta \l( \Theta_{b_L}^{(j)}, \Lambda \r) - \nabla_{b_1} \caL_\beta \l( \Theta^{(j)}, \Lambda \r) = \beta \l( \Delta_j U_{1} \r) 1_m. 
    \end{align*}
    Thus, for $l > 1$, the bound is:
    \begin{align*}
        & \left\| \nabla_{b_l} \caL_\beta \l( \Theta_{b_L}^{(j)}, \Lambda \r) - \nabla_{b_l} \caL_\beta \l( \Theta^{(j)}, \Lambda \r) \right\|^2 \\
        \le \; & 2 \beta^2 m \left( \delta^{k-1}_{0,l} \|\Delta_j V_{l-1}\|^2_{\F} + \|\Delta_j U_{l}\|^2_{\F} \right) \\
        \le \; & 2 \beta^2 m \max \l\{ \delta^{k-1}_{0,l}, \, 1 \r\} \left( \|\Delta_j V_{l-1}\|^2_{\F} + \|\Delta_j U_{l}\|^2_{\F} \right),
    \end{align*}     
    and for $l = 1$:
    \begin{align*}
         \left\| \nabla_{b_1} \caL_\beta \l( \Theta_{b_L}^{(j)}, \Lambda \r) - \nabla_{b_1} \caL_\beta \l( \Theta^{(j)}, \Lambda \r) \right\|^2 \le  m \beta^2 \|\Delta_j U_{1}\|^2_{\F}.
    \end{align*}     
    Therefore, we can get \eqref{upper.bound.b.l} with $L_{2,l}\l( \lam, \beta, \hat{r}^{k-1} \r)$ defined in \eqref{upper.bound.l.constants}.

    For $V_l$ ($l = 1, \dots, L-1$), the gradient of $\caL_\beta$ with respect to $V_l$ is:
    \begin{align*}
        \nabla_{V_l} \caL_\beta \left(\Theta, \Lambda\right) = W_{l+1}^{\T} \Lambda_{1,l+1} + \beta W_{l+1}^{\T} \left( W_{l+1} V_l + b_{l+1} 1_m^{\T} - U_{l+1} \right) + \Lambda_{2,l} + \beta \left(V_l - \sigma(U_l)\right).
    \end{align*}
    Therefore, for $l = 1, \dots, L-1$, the difference is: 
    \begin{align*}
    \nabla_{V_l} \caL_\beta \l( \Theta_{V_L}^{(j)}, \lam \r) - \nabla_{V_l} \caL_{\beta} \l( \Theta^{(j)}, \lam \r) = \beta \l[ {W_{l+1}^{(j)}}^{\T} \l( \Delta_j U_{l+1} \r) + \sigma \l(U_l^{(j)}\r) - \sigma \l(U_l^{(j-1)}\r) \r].
    \end{align*}
    The bound is:
    \begin{align*}
    & \l\| \nabla_{V_l} \caL_\beta \l( \Theta_{V_L}^{(j)}, \lam \r) - \nabla_{V_l} \caL_{\beta} \l( \Theta^{(j)}, \lam \r) \r\|_{\F}^2\\
    \leq \; & 2 \beta^2  \l( \delta_{0,l+1}^{k-1} \|\Delta_j U_{l+1}\|_{\F}^2 + \l\| \sigma \l(U_l^{(j)}\r) - \sigma \l(U_l^{(j-1)}\r) \r\|_{\F}^2 \r) \\
    \leq \; & 2 \beta^2 \max\l\{ \delta_{0,l+1}^{k-1}, 1 \r\} \l( \|\Delta_j U_{l+1}\|_{\F}^2 + \|\Delta_j U_{l}\|_{\F}^2 \r).
    \end{align*} 
    For $V_L$, we have the following:
    \begin{align*}
    \nabla_{V_L} \caL_\beta \l( \Theta_{V_L}^{(j)}, \lam \r) - \nabla_{V_L} \caL_{\beta} \l( \Theta^{(j)}, \lam \r) = \beta \l( \sigma \l(U_L^{(j)}\r) - \sigma \l(U_L^{(j-1)}\r) \r).
    \end{align*}
    The bound is:
    \begin{align*}
    \l\| \nabla_{V_L} \caL_\beta \l( \Theta_{V_L}^{(j)}, \lam \r) - \nabla_{V_L} \caL_{\beta} \l( \Theta^{(j)}, \lam \r) \r\|_{\F}^2 \leq \beta^2 \|\Delta_j U_{L}\|_{\F}^2.
    \end{align*}    
    Therefore, we can get \eqref{upper.bound.v.l} with $L_{3,l}\l( \lam, \beta, \hat{r}^{k-1} \r)$ defined in \eqref{upper.bound.l.constants}.
    This completes the proof.
\end{proof}

The following lemma verifies the level-set memberships required in Lemma~\ref{le.partial.bound}, including those of the intermediate states.

\begin{lemma}\label{le.s.level.set}

Fix an outer iteration $k$, with $\beta = \beta^{k-1}$ and $\Lambda=\Lambda^{k-1}$. Let $\{\Theta^{(j)}\}$ be generated by Algorithm~\ref{alg:bcd}, initialized by \eqref{eq:3.13}, and let $\hat r^{k-1}=\caL_\beta(\Theta^{(0)},\Lambda)$. Then $\Theta^{(j)}\in\Omega_{\caL_\beta}(\hat r^{k-1})$ for every
$j\ge0$. Every intermediate state $\Theta_{W_l}^{(j)},\Theta_{Z_l}^{(j)},\Theta_{b_l}^{(j)},
\Theta_{V_l}^{(j)},\Theta_{U_l}^{(j)}$, for $j\ge1$ and $l=1,\ldots,L$, also belongs to this compact level set. Moreover,
\[
\Omega_{\caL_\beta}(\hat r^{k-1})
\subseteq\Omega_{\caL_\beta}(\hat c).
\]
\end{lemma}

\begin{proof}
Each $W$, $b$, and $V$ update minimizes $\caL_\beta$ over its block, so it cannot increase $\caL_\beta$. For a $Z$ or $U$ update, global optimality of the proximal subproblem, evaluated against the preceding
block value, gives
\[
\caL_\beta(\Theta_{\rm new},\Lambda)
+\frac{\mu}{2}\|x_{\rm new}-x_{\rm old}\|_{\F}^2
\le\caL_\beta(\Theta_{\rm old},\Lambda),
\]
where $(x,\mu)$ is $(Z_l,\mu_1)$ or $(U_l,\mu_2)$. Thus $\caL_\beta$ is nonincreasing throughout every sweep, including all intermediate states. In particular,
\begin{equation}\label{eq.suff.des.bcd.aug}
\caL_\beta(\Theta^{(j)},\Lambda)
\le\caL_\beta(\Theta^{(j-1)},\Lambda),\qquad j\ge1.
\end{equation}
The initialization \eqref{eq:3.13} guarantees $\hat r^{k-1}\le\hat c$. Induction over the individual block updates therefore yields
\begin{equation}\label{eq.bcd.aug.bounded.gamma}
\caL_\beta(\Theta^{(j)},\Lambda)
\le\hat r^{k-1}\le\hat c,\qquad j\ge0,
\end{equation}
and the same bound holds at every intermediate state. All these states belong to $\mathcal D$, since the $Z$ updates remain on their prescribed Stiefel manifolds. Lemma~\ref{le.bcd.basic} gives compactness of
$\Omega_{\caL_\beta}(\hat r^{k-1})$, completing the proof.
\end{proof}

Finally, the following theorem establishes that step $3$ of Algorithm~\ref{alg:alm} is well-defined when executed via Algorithm~\ref{alg:bcd}. Furthermore, we prove that condition \eqref{eq:3.6} is satisfied whenever the stopping criterion \eqref{eq:3.19} is met.

\begin{theorem}\label{th.8.3.1}
    At the $k$-th iteration of the augmented Lagrangian method (Algorithm~\ref{alg:alm}), the alternating minimization method (Algorithm~\ref{alg:bcd}) for solving the subproblem \eqref{eq:3.5} satisfies the stopping criterion \eqref{eq:3.19} in a finite number of steps, requiring $\mathcal{O}(1/\epsilon^2)$ iterations. Furthermore, the resulting output $\Theta^k$ generated by Algorithm~\ref{alg:bcd} guarantees the satisfaction of condition \eqref{eq:3.6} required by Algorithm~\ref{alg:alm}.
\end{theorem}

\begin{proof}
    Combining the sufficient descent properties established before, we obtain the following overall descent inequality after one complete iteration of Algorithm~\ref{alg:bcd}:
\begin{align*}
    & \caL_\beta\left(\Theta^{(j - 1)}, \Lambda \right) - \caL_\beta\left(\Theta^{(j)}, \Lambda \right) \\
    \ge \; & \frac{\beta}{2}\sum_{l = 1}^L \|\Delta_j W_{l}\|_{\F}^2 + \frac{m\beta + \lambda_1}{2} \sum_{l = 1}^L \|\Delta_j b_{l}\|^2 + \frac{\mu_1}{2}\sum_{l = 1}^L\|\Delta_j Z_{l}\|_{\F}^2 \\
    & + \frac{\beta}{2} \sum_{l = 1}^{L-1} \|\Delta_j V_{l}\|_{\F}^2 + \left( \frac{\beta}{2} + \frac{1}{m} \right) \|\Delta_j V_{L}\|_{\F}^2 + \frac{\mu_2}{2}\sum_{l = 1}^L \|\Delta_j U_{l}\|_{\F}^2 \\
    \ge \; & L_{\mathrm{des}} \sum_{l = 1}^L \bigg( \|\Delta_j b_{l}\|^2 + \|\Delta_j V_{l}\|_{\F}^2 + \|\Delta_j Z_{l}\|_{\F}^2 + \|\Delta_j U_{l}\|_{\F}^2 \bigg),
\end{align*}
    where $L_{\mathrm{des}} := \min \left\{\beta, \, \mu_1, \, \mu_2 \right\} / 2$. 
    Summing $\caL_\beta\left(\Theta^{(j - 1)}, \Lambda \right) - \caL_\beta\left(\Theta^{(j)}, \Lambda \right)$ from $j = 1$ to $J$: 
\begin{align*}
    & \caL_\beta\left(\Theta^{(0)}, \Lambda \right) - \caL_\beta\left(\Theta^{(J)}, \Lambda \right) \\
    \ge \; & L_{\mathrm{des}} \sum_{j = 1}^J \Bigg[ \sum_{l = 1}^L \bigg( \|\Delta_j b_{l}\|^2 + \|\Delta_j V_{l}\|_{\F}^2 + \|\Delta_j Z_{l}\|_{\F}^2 + \|\Delta_j U_{l}\|_{\F}^2 \bigg) \Bigg] \\
    \ge \; & J L_{\mathrm{des}} \min_{1 \le j \le J} \Bigg\{ \sum_{l = 1}^L \bigg( \|\Delta_j b_{l}\|^2 + \|\Delta_j V_{l}\|_{\F}^2 + \|\Delta_j Z_{l}\|_{\F}^2 + \|\Delta_j U_{l}\|_{\F}^2 \bigg) \Bigg\}.
\end{align*}  
    This, together with Lemma~\ref{le.bcd.basic}, yields that: 
\begin{align*}
    & \min_{1 \le j \le J} \Bigg\{ \sum_{l = 1}^L \bigg( \|\Delta_j b_{l}\|^2 + \|\Delta_j V_{l}\|_{\F}^2 + \|\Delta_j Z_{l}\|_{\F}^2 + \|\Delta_j U_{l}\|_{\F}^2 \bigg) \Bigg\} \\
    \le \; & \frac{\caL_\beta\left(\Theta^{(0)}, \Lambda \right) + \frac{1}{2\beta}\sum_{l = 1}^L \left( \left\| \Lambda_{1,l} \right\|_{\F}^2 + \left\| \Lambda_{2,l} \right\|_{\F}^2 + \left\| \Lambda_{3,l} \right\|_{\F}^2 \right)}{J L_{\mathrm{des}}}. 
\end{align*}
    It follows that there exists $j\in\{1,\ldots,J\}$ satisfying the stopping criterion \eqref{eq:3.19} whenever: 
    \[
    J \ge \hat J:=\max\left\{1,\left\lceil\frac{L_{\max}^{k-1}\delta^{k-1}}{L_{\mathrm{des}}\epsilon^2}\right\rceil\right\}.
    \]
    Here $\delta^{k-1}=\caL_\beta(\Theta^{(0)},\Lambda)+\sum_{l=1}^L\sum_{a=1}^3\|\Lambda_{a,l}\|_{\F}^2/(2\beta)$ is the nonnegative gap bound defined in Lemma~\ref{le.partial.bound}. If this gap is zero, the first sweep already has zero increments, which explains the maximum with $1$. Therefore, at the $k$-th iteration of Algorithm~\ref{alg:alm}, the alternating minimization method (Algorithm~\ref{alg:bcd}) can be stopped in at most $\hat{J}$ iterations and output $\Theta^k$, which corresponds to a complexity of order $\mathcal{O}(1/\epsilon^2)$.

We now prove that satisfying the stopping criterion \eqref{eq:3.19} guarantees the satisfaction of condition \eqref{eq:3.6}. Following the update rules of Algorithm~\ref{alg:bcd}, the first-order optimality conditions for the respective subproblems for $l = 1, \dots, L$ are given by: 
\begin{align*}
    & \nabla_{W_l}\caL_\beta \l( \Theta_{W_l}^{(j)}, \Lambda \r) = 0, \qquad {\mit\Pi}_{Z_l^{(j)}}\Big(-\Lambda_{3,l}
      +\beta\big[Z_l^{(j)}-\mathcal T_l(W_l^{(j)})\big]
      +\mu_1\Delta_jZ_l\Big)=0, \\
    & \nabla_{b_l} \caL_\beta \l( \Theta_{b_l}^{(j)}, \Lambda \r) = 0, \qquad \nabla_{V_l} \caL_\beta \l( \Theta_{V_l}^{(j)}, \Lambda \r) = 0, \qquad 0 \in \mu_2 \l( \Delta_j U_{l} \r) + \partial_{U_l} \caL_\beta \l( \Theta^{(j)}, \Lambda \r).
\end{align*}  
The $U_l$ inclusion holds at the full iterate $\Theta^{(j)}$: all $V$ blocks are updated before the $U$ blocks, and the $U$ subproblems are separable across layers. Thus $-\mu_2\Delta_jU_l$ is an actual member of the limiting subdifferential at that iterate. No comparison of subdifferentials at two different $U$ values is required.

By substituting these first-order optimality conditions, the squared Frobenius norm of the gradients of $\caL_{\beta}$ with respect to $W_l$, $b_l$, and $V_l$, evaluated at the fully updated iterate $\Theta^{(j)}$, can be expressed for all $l = 1, \dots, L$ as:
\begin{align*}
    \l\| \nabla_{W_l}\caL_\beta ( \Theta^{(j)}, \Lambda ) \r\|_{\F}^2 &= \l\| \nabla_{W_l}\caL_\beta ( \Theta^{(j)}, \Lambda ) - \nabla_{W_l}\caL_\beta ( \Theta_{W_l}^{(j)}, \Lambda ) \r\|_{\F}^2, \\
    \l\| \nabla_{b_l}\caL_\beta ( \Theta^{(j)}, \Lambda ) \r\|^2 &= \l\| \nabla_{b_l}\caL_\beta ( \Theta^{(j)}, \Lambda ) - \nabla_{b_l}\caL_\beta ( \Theta_{b_l}^{(j)}, \Lambda ) \r\|^2, \\
    \l\| \nabla_{V_l}\caL_\beta ( \Theta^{(j)}, \Lambda ) \r\|_{\F}^2 & = \l\| \nabla_{V_l}\caL_\beta ( \Theta^{(j)}, \Lambda ) - \nabla_{V_l} \caL_\beta \l( \Theta_{V_l}^{(j)}, \Lambda \r)\r\|_{\F}^2.
\end{align*}
Using the upper bounds \eqref{upper.bound.w.l} derived in Lemma~\ref{le.partial.bound}, we have the following inequality for $W_l$:
\begin{align*}
    & \sum_{l = 1}^L \left\| \nabla_{W_l}\caL_\beta \left( \Theta^{(j)}, \Lambda \right) \right\|_{\F}^2 \\
    = & \; \sum_{l = 1}^L \left\|\nabla_{W_l}\caL_\beta \left( \Theta^{(j)}, \Lambda \right) - \nabla_{W_l}\caL_\beta \left( \Theta_{W_l}^{(j)}, \Lambda \right) \right\|_{\F}^2 \\
    \le & \; \sum_{l=1}^L L_{1,l}^{k-1} \left( \|\Delta_j Z_{l}\|^2_{\F} + \|\Delta_j b_{l}\|^2 + \|\Delta_j U_{l}\|_{\F}^2 \right) + \sum_{l = 2}^L L_{1,l}^{k-1} \|\Delta_j V_{l-1}\|_{\F}^2.
\end{align*}
For $Z_l$ ($l=1,\ldots,L$), the gradient depends only on $W_l$ and $Z_l$, which have already been updated. Linearity and non-expansiveness of the tangent projection give
\begin{align*}
&\sum_{l = 1}^L \l\|{\mit\Pi}_{Z_l^{(j)}}\big(\nabla_{Z_l}\caL_\beta(\Theta^{(j)},\Lambda)\big)\r\|_{\F}^2\\
={}&\sum_{l = 1}^L \l\|{\mit\Pi}_{Z_l^{(j)}}\big(\nabla_{Z_l}\caL_\beta(\Theta^{(j)},\Lambda)+\mu_1\Delta_jZ_l\big)
-{\mit\Pi}_{Z_l^{(j)}}\big(\nabla_{Z_l}\caL_\beta(\Theta^{(j)},\Lambda)\big)\r\|_{\F}^2\\
={}& \mu_1^2 \sum_{l=1}^L \l\|{\mit\Pi}_{Z_l^{(j)}}(\Delta_jZ_l)\r\|_{\F}^2
\le \mu_1^2 \sum_{l=1}^L \|\Delta_jZ_l\|_{\F}^2.
\end{align*}
    Using the upper bounds \eqref{upper.bound.b.l} derived in Lemma~\ref{le.partial.bound}, we have the following inequality for $b_l$:
    \begin{align*}
        \sum_{l = 1}^L \left\| \nabla_{b_l} \caL_\beta \left( \Theta^{(j)}, \Lambda \right) \right\|^2 & = \sum_{l = 1}^L \left\| \nabla_{b_l}\caL_\beta \left( \Theta^{(j)}, \Lambda \right) - \nabla_{b_l}\caL_\beta \left( \Theta_{b_l}^{(j)}, \Lambda \right) \right\|^2 \\
        & \le \sum_{l = 1}^L L_{2,l}^{k-1} \|\Delta_j U_{l}\|_{\F}^2 + \sum_{l = 2}^L L_{2,l}^{k-1} \|\Delta_j V_{l-1}\|_{\F}^2.
    \end{align*}
    Using \eqref{upper.bound.v.l} and the convention $L_{3,0}^{k-1}=0$, we obtain:
    \begin{align*}
    \sum_{l = 1}^L \l\| \nabla_{V_l} \caL_\beta \left( \Theta^{(j)}, \Lambda \right) \r\|_{\F}^2 & = \sum_{l = 1}^L \left\| \nabla_{V_l}\caL_\beta \left( \Theta^{(j)}, \Lambda \right) - \nabla_{V_l}\caL_\beta \left( \Theta_{V_l}^{(j)}, \Lambda \right) \right\|_{\F}^2 \\
    & \le \sum_{l = 1}^L \l( L_{3,l}^{k-1} + L_{3,l-1}^{k-1} \r) \|\Delta_j U_{l}\|_{\F}^2.
    \end{align*}
    For $U_l$, we can derive the following: 
\begin{align*}
    \sum_{l = 1}^L \mathrm{dist} \left( 0, \partial_{U_l}\caL_\beta \left( \Theta^{(j)}, \Lambda \right)\right)^2 \le \mu_2^2 \sum_{l = 1}^L \|\Delta_j U_{l}\|_{\F}^2.
\end{align*}
    Combining all these upper bounds, we evaluate the following squared distance at $\Theta^{(j)}$: 
\begin{align*}
    & \mathrm{dist} \left( 0, \partial_{\Theta} \caL_\beta \left( \Theta^{(j)}, \Lambda \right) \right)^2 \\
    = \; & \sum_{l = 1}^L \bigg( \left\| \nabla_{W_l} \caL_{\beta}\left( \Theta^{(j)}, \Lambda \right) \right\|_{\F}^2 + \left\| {\mit{\Pi}}_{Z_l^{(j)}} \l( \nabla_{Z_l} \caL_{\beta} \l( \Theta^{(j)}, \lam \r) \r) \right\|_{\F}^2 + \left\| \nabla_{b_l} \caL_{\beta}\left( \Theta^{(j)}, \Lambda \right) \right\|^2 \\
        & \qquad \quad + \min_{\tilde{U}_l^{(j)} \in \partial_{U_l} \caL_{\beta}\left( \Theta^{(j)}, \Lambda \right)}\left\| \tilde{U}^{(j)}_l \right\|_{\F}^2 + \left\| \nabla_{V_l} \caL_{\beta}\left( \Theta^{(j)}, \Lambda \right) \right\|_{\F}^2 \bigg) \\
    \le \; & \sum_{l = 1}^L \l( L^{k-1}_{1,l} + \mu_1^2 \r) \|\Delta_j Z_{l}\|_{\F}^2 + \sum_{l = 1}^L L^{k-1}_{1,l} \|\Delta_j b_{l}\|^2 + \sum_{l = 1}^{L-1} \l( L_{1,l+1}^{k-1} + L_{2,l+1}^{k-1} \r) \|\Delta_j V_{l}\|_{\F}^2 \\
    & + \sum_{l = 1}^{L} \l( L^{k-1}_{1,l} + L^{k-1}_{2,l} + L_{3,l}^{k-1} + L_{3,l-1}^{k-1} + \mu_2^2\r) \|\Delta_j U_{l}\|_{\F}^2 \\
    \le \; & L_{\max}^{k-1} \sum_{l = 1}^L \bigg( \|\Delta_j Z_{l}\|_{\F}^2 + \|\Delta_j b_{l}\|^2 + \|\Delta_j U_{l}\|_{\F}^2 + \|\Delta_j V_{l}\|_{\F}^2 \bigg),
\end{align*}
    where  
    \begin{align*}
        L_{\max}^{k-1} = \max \bigg\{ & \max_{1 \le l \le L} \l\{ L_{1,l}^{k-1} \r\} + \mu_1^2, \; \max_{1 \le l \le L-1} \l\{ L^{k-1}_{1,l+1} + L^{k-1}_{2,l+1} \r\}, \\ & \max_{1 \le l \le L} \l\{ L_{1,l}^{k-1} + L_{2,l}^{k-1} + L_{3,l}^{k-1} + L_{3,l-1}^{k-1} \r\} + \mu_2^2 \bigg\}.
    \end{align*}
    Therefore, the satisfaction of condition \eqref{eq:3.19}, which states:
\begin{align*}
    \sum_{l = 1}^L \bigg( \|\Delta_j Z_{l}\|_{\F}^2 + \|\Delta_j b_{l}\|^2 + \|\Delta_j U_{l}\|_{\F}^2 + \|\Delta_j V_{l}\|_{\F}^2 \bigg) \le \frac{\epsilon^2}{L_{\max}^{k-1}},
\end{align*}
    together with $\Theta^k = \Theta^{(j)}$, implies that: 
    \[
    \mathrm{dist}\left( 0, \, \partial_{\Theta} \caL_\beta \left( \Theta^k, \Lambda \right) \right) \le \epsilon.
    \]
    Note that the update rule \eqref{z.l.sub.solu} ensures that the variables $Z_l^{(j)}$ ($l = 1,\dots,L$) inherently remain on the Stiefel manifold. 
    This confirms that the condition \eqref{eq:3.6} is fully satisfied. 
\end{proof}

\subsection{Convergence of Algorithm~\ref{alg:alm}}\label{subsec:globalconv}

This subsection focuses on establishing the convergence properties of the augmented Lagrangian method (Algorithm~\ref{alg:alm}), where step~$3$ of Algorithm~\ref{alg:alm} is performed by Algorithm~\ref{alg:bcd} with initialization \eqref{eq:3.13}, fixed positive proximal parameters, and the stated stopping criterion. We begin by proving that any accumulation point generated by Algorithm~\ref{alg:alm} is a feasible point.

\begin{theorem}\label{th.alm.feas}
    Let $\{ \Theta^k \}$ be the sequence generated by Algorithm~\ref{alg:alm}.
    Then the following statements hold:
    \begin{enumerate}
        \item For $h\in\{1,2,3\}$ and $l\in\{1,\ldots,L\}$, we have
            \[
                \lim_{k\to\infty}\|E_{h,l}^k\|_{\F}=0.
            \]

        \item There exists at least one accumulation point of $\{ \Theta^k \}$, and any accumulation point is a feasible point of problem \eqref{prob.orth.three.sets.vari}.
    \end{enumerate}
\end{theorem}

\begin{proof}
\begin{enumerate}
\item Define the set of penalty-increase indices
\begin{equation}\label{set.beta.k}
 \caK:=\left\{k\in\mathbb N_+:
 \beta^k=\max\left\{\frac{\beta^{k-1}}{\eta_2},
       \max_{h,l}\|\Lambda_{h,l}^k\|_{\F}^{1+\eta_3}\right\}\right\}.
\end{equation}
If $\caK$ is finite, choose $K$ so that $k\notin\caK$ for every $k>K$.
The penalty update rule gives
\begin{equation}\label{eq.feas.finite}
 \max_{h,l}\|E_{h,l}^k\|_{\F}
 \le\eta_1\max_{h,l}\|E_{h,l}^{k-1}\|_{\F}
 \le\eta_1^{k-K}\max_{h,l}\|E_{h,l}^{K}\|_{\F}.
\end{equation}
Since $0<\eta_1<1$, it follows that $\caR^k\to0$, proving the first
statement in this case.

Suppose next that $\caK$ is infinite. Each increase multiplies the
penalty by at least $1/\eta_2>1$, so the nondecreasing sequence
$\beta^{k-1}$ tends to infinity. For $k-1\in\caK$, the update gives
\[
 \begin{aligned}
 \max_{h,l}\frac{\|\Lambda_{h,l}^{k-1}\|_{\F}}{\beta^{k-1}}
 \le(\beta^{k-1})^{\frac{-\eta_3}{1+\eta_3}}, \qquad
 \max_{h,l}\frac{\|\Lambda_{h,l}^{k-1}\|_{\F}^2}{\beta^{k-1}}
 \le(\beta^{k-1})^{\frac{1-\eta_3}{1+\eta_3}}.
 \end{aligned}
\]
Thus $\eta_3>1$ implies
\begin{equation}\label{eq.lag.mul.norm.k-1}
 \lim_{\substack{k\to\infty\\ k-1\in\caK}}
 \max_{h,l}\left\{
 \frac{\|\Lambda_{h,l}^{k-1}\|_{\F}}{\beta^{k-1}},
 \frac{\|\Lambda_{h,l}^{k-1}\|_{\F}^2}{\beta^{k-1}}\right\}=0.
\end{equation}
Completing the squares in \eqref{eq.aug.lag} and using
\eqref{eq.bcd.aug.bounded.gamma} and the initialization of
Algorithm~\ref{alg:bcd}, we obtain, for every $k$,
\begin{equation}\label{eq.aug.norm.bound}
 \begin{aligned}
 0\le\sum_{l=1}^L\sum_{h=1}^3
  \left\|E_{h,l}^k+\frac{\Lambda_{h,l}^{k-1}}{\beta^{k-1}}\right\|_{\F}^2 \le \frac{2\bigl(\hat c-\caF(\Theta^k)\bigr)}{\beta^{k-1}}
  +\sum_{l=1}^L\sum_{h=1}^3
      \frac{\|\Lambda_{h,l}^{k-1}\|_{\F}^2}{(\beta^{k-1})^2}.
 \end{aligned}
\end{equation}
Because $\caF(\Theta^k)\ge0$, equations
\eqref{eq.lag.mul.norm.k-1}--\eqref{eq.aug.norm.bound} and
$\beta^{k-1}\to\infty$ imply
\begin{equation}\label{eq.cons.norm.bound}
 \lim_{\substack{k\to\infty\\ k-1\in\caK}}
       \|E_{h,l}^k\|_{\F}=0,
 \qquad h=1,2,3,\quad l=1,\ldots,L.
\end{equation}

To extend this limit to all indices, let $p_k$ be the largest member
of $\caK$ strictly less than $k$, for all sufficiently large $k$.
Then $p_k\to\infty$. If $p_k=k-1$, the required estimates follow
directly from \eqref{eq.lag.mul.norm.k-1} and
\eqref{eq.cons.norm.bound}. If $p_k<k-1$, put
$\mathcal J_k:=\{i:p_k<i<k\}$. The penalty remains constant on this
interval: $\beta^i=\beta^{i-1}=\beta^{p_k}$ for $i\in\mathcal J_k$.
For every such $i$, the multiplier update yields
\begin{equation}\label{eq.lam.beta.inqua}
 \begin{aligned}
 \frac{\|\Lambda_{h,l}^i\|_{\F}}{\beta^i}
 =\frac{\|\Lambda_{h,l}^{i-1}+\beta^{i-1}E_{h,l}^i\|_{\F}}
              {\beta^{i-1}} \le \frac{\|\Lambda_{h,l}^{i-1}\|_{\F}}{\beta^{i-1}}
      +\|E_{h,l}^i\|_{\F}.
 \end{aligned}
\end{equation}
Summing over $i\in\mathcal J_k$ gives
\begin{equation}\label{eq.lam.beta.inqua.sum}
 \frac{\|\Lambda_{h,l}^{k-1}\|_{\F}}{\beta^{k-1}}
 \le\frac{\|\Lambda_{h,l}^{p_k}\|_{\F}}{\beta^{p_k}}
      +\sum_{i=p_k+1}^{k-1}\|E_{h,l}^i\|_{\F}.
\end{equation}
For $p_k+1\le i\le k-1$, the successful residual-reduction tests give
$\|E_{h,l}^i\|_{\F}\le\caR^i
 \le\eta_1^{i-p_k-1}\caR^{p_k+1}$. Hence
\[
 \begin{aligned}
 \frac{\|\Lambda_{h,l}^{k-1}\|_{\F}}{\beta^{k-1}}
 \le\frac{\|\Lambda_{h,l}^{p_k}\|_{\F}}{\beta^{p_k}}
      +\caR^{p_k+1}\sum_{j=0}^{k-p_k-2}\eta_1^j \le \frac{\|\Lambda_{h,l}^{p_k}\|_{\F}}{\beta^{p_k}}
      +\frac{\caR^{p_k+1}}{1-\eta_1}.
 \end{aligned}
\]
Since $p_k\in\caK$ and $p_k\to\infty$, equations
\eqref{eq.lag.mul.norm.k-1} and \eqref{eq.cons.norm.bound} show that
both terms on the last line tend to zero. Combining the cases
$p_k=k-1$ and $p_k<k-1$, we conclude that
\begin{equation}\label{eq.lim.lambda.bounded.inf}
 \lim_{k\to\infty}
 \frac{\|\Lambda_{h,l}^{k-1}\|_{\F}}{\beta^{k-1}}=0,
 \qquad h=1,2,3,\quad l=1,\ldots,L.
\end{equation}
Applying \eqref{eq.aug.norm.bound} again, now along the full sequence,
and using $\caF(\Theta^k)\ge0$ and $\beta^{k-1}\to\infty$, gives
\begin{equation}\label{eq.cons.norm.bound.all.k}
 \lim_{k\to\infty}\|E_{h,l}^k\|_{\F}=0,
 \qquad h=1,2,3,\quad l=1,\ldots,L.
\end{equation}
This proves the first statement also when $\caK$ is infinite.

\item First suppose that $\caK$ is finite, and choose $K$ as above.
Then $\beta^{k-1}=\beta^K$ for $k>K$. Summing the multiplier updates
and using \eqref{eq.feas.finite}, we obtain
\[
 \begin{aligned}
 \max_{h,l}\|\Lambda_{h,l}^{k-1}\|_{\F}
 \le\max_{h,l}\|\Lambda_{h,l}^{K}\|_{\F}
       +\beta^K\sum_{i=K+1}^{k-1}\caR^i \le \max_{h,l}\|\Lambda_{h,l}^{K}\|_{\F}
       +\frac{\beta^K\eta_1}{1-\eta_1}\caR^K.
 \end{aligned}
\]
Thus the multipliers are bounded, and so are
$\|\Lambda_{h,l}^{k-1}\|_{\F}^2/\beta^{k-1}$ for all $h,l$.
The augmented-Lagrangian bound gives
\begin{equation}\label{eq.bounded.for.s}
 \begin{aligned}
 \caF(\Theta^k)+\frac{\beta^{k-1}}2
   \sum_{l=1}^L\sum_{h=1}^3
   \left\|E_{h,l}^k+\frac{\Lambda_{h,l}^{k-1}}{\beta^{k-1}}\right\|_{\F}^2 \le \hat c+\frac{1}{2\beta^{k-1}}
       \sum_{l=1}^L\sum_{h=1}^3\|\Lambda_{h,l}^{k-1}\|_{\F}^2.
 \end{aligned}
\end{equation}
As in Lemma~\ref{le.bcd.basic}(2), the positive regularization terms
bound $\caB^k$ and $\caU^k$. The shifted $E_{2,l}^k$ terms then bound
$\caV^k$, while the shifted $E_{3,l}^k$ terms and the compactness of
the Stiefel manifolds bound $\caW^k$ and $\caZ^k$.
Here the multiplier shifts are bounded and $\beta^{k-1}\ge\beta^0>0$.
Consequently, $\{\Theta^k\}$ is bounded and has an accumulation point.

If $\caK$ is infinite, \eqref{eq.lag.mul.norm.k-1} gives boundedness
of $\|\Lambda_{h,l}^{k-1}\|_{\F}^2/\beta^{k-1}$ along $k-1\in\caK$.
In particular, the multiplier shifts are bounded because
\[
 \left\|\frac{\Lambda_{h,l}^{k-1}}{\beta^{k-1}}\right\|_{\F}^2
 \le\frac{1}{\beta^0}
       \frac{\|\Lambda_{h,l}^{k-1}\|_{\F}^2}{\beta^{k-1}}.
\]
Equation~\eqref{eq.bounded.for.s} and the same blockwise boundedness
argument therefore show that $\{\Theta^k\}_{k-1\in\caK}$ is bounded.
This infinite subsequence has an accumulation point, which is also
an accumulation point of $\{\Theta^k\}$.

Finally, all $Z_l^k$ belong to their respective Stiefel manifolds,
which are closed. Continuity of the residuals, together with the
first statement, shows that every accumulation point satisfies all
equality constraints. Hence every accumulation point is feasible.
\end{enumerate}
\end{proof}

Finally, we establish the overall convergence of the augmented Lagrangian method (Algorithm~\ref{alg:alm}).

\begin{theorem}\label{theo.conver.orthalm.kkt}
    Any accumulation point of the sequence $\{ \Theta^k \}$ generated by Algorithm~\ref{alg:alm} is a KKT point of the regularized problem~\eqref{prob.orth.three.sets.vari}.
\end{theorem}

\begin{proof}

    Let $\l\{\Theta^{k_i}\r\}$ be a subsequence of $\l\{\Theta^k\r\}$ converging to $\Theta^\ast$. Then we know that $\Theta^\ast$ is a feasible point of \eqref{prob.orth.three.sets.vari} by Theorem~\ref{th.alm.feas}. We claim that
    \begin{equation}\label{eqvi.augl.lag}
        \partial_{\Theta} \caL_{\beta^{k_i - 1}} \l( \Theta^{k_i}, \lam^{k_i - 1} \r) = \partial_{\Theta} \caL \l( \Theta^{k_i}, \lam^{k_i} \r).
    \end{equation}
    First, by direct computation and the multiplier updates, for $W_l$ ($l=1,\ldots,L$), we have
\begin{align*}
 &\nabla_{W_l}\caL_{\beta^{k_i-1}}
       (\Theta^{k_i},\Lambda^{k_i-1})\\
 ={}&\Lambda_{1,l}^{k_i-1}{V_{l-1}^{k_i}}^{\T}
      +\beta^{k_i-1}E_{1,l}^{k_i}{V_{l-1}^{k_i}}^{\T}
      +\mathcal T_l^*\big(\Lambda_{3,l}^{k_i-1}
                         +\beta^{k_i-1}E_{3,l}^{k_i}\big)\\
 ={}&\Lambda_{1,l}^{k_i}{V_{l-1}^{k_i}}^{\T}
      +\mathcal T_l^*(\Lambda_{3,l}^{k_i}) = \nabla_{W_l}\caL(\Theta^{k_i},\Lambda^{k_i}).
\end{align*}
Here $\mathcal T_l^*$ represents the two orientations in
\eqref{orthdnn.kkt.w.i}--\eqref{orthdnn.kkt.w.i.c}.
For $Z_l$ ($l=1,\ldots,L$), we have
\begin{align*}
 & {\mit\Pi}_{Z_l^{k_i}}\big(
       \nabla_{Z_l}\caL_{\beta^{k_i-1}}
          (\Theta^{k_i},\Lambda^{k_i-1})\big) \\
          = {}& -{\mit\Pi}_{Z_l^{k_i}}\big(
          \Lambda_{3,l}^{k_i-1}+\beta^{k_i-1}E_{3,l}^{k_i}\big) = -{\mit\Pi}_{Z_l^{k_i}}(\Lambda_{3,l}^{k_i}) = {\mit\Pi}_{Z_l^{k_i}}\big(
         \nabla_{Z_l}\caL(\Theta^{k_i},\Lambda^{k_i})\big).
\end{align*}
    For $b_l$ ($l = 1,\dots,L$), we have
    \[
    \begin{split}
        & \nabla_{b_l} \caL_{\beta^{k_i - 1}} \l( \Theta^{k_i}, \lam^{k_i - 1} \r) \\
        = \; & \lambda_1 b_l^{k_i} + \lam_{1,l}^{k_i - 1} 1_m + \beta^{k_i - 1} \l( W^{k_i}_l V_{l-1}^{k_i} + b_l^{k_i} 1_m^{\T} - U_l^{k_i} \r) 1_m\\
        = \; & \lambda_1 b_l^{k_i} + \lam_{1,l}^{k_i} 1_m = \nabla_{b_l} \caL \l( \Theta^{k_i}, \lam^{k_i} \r).
    \end{split}
    \]    
    For $V_l$ ($l = 1,\dots,L-1$), we have
\[
    \begin{aligned}
        & \nabla_{V_l} \caL_{\beta^{k_i - 1}} \left( \Theta^{k_i}, \lam^{k_i - 1} \right) \\
        = \; & {W_{l+1}^{k_i}}^{\T} \lam_{1,l+1}^{k_i-1} + \beta^{k_i - 1} {W_{l+1}^{k_i}}^{\T} \left( W_{l+1}^{k_i} V_l^{k_i} + b_{l+1}^{k_i} 1_m^{\T} - U_{l+1}^{k_i} \right) \\
        & + \lam_{2,l}^{k_i - 1} + \beta^{k_i - 1} \left( V_l^{k_i} - \sigma(U_l^{k_i}) \right) \\
        = \; & {W_{l+1}^{k_i}}^{\T} \lam_{1,l+1}^{k_i} + \lam_{2,l}^{k_i} = \nabla_{V_l} \caL \left( \Theta^{k_i}, \lam^{k_i} \right).
    \end{aligned}
\]
    For $V_L$, we have
    \[
    \begin{split}
        & \nabla_{V_L} \caL_{\beta^{k_i - 1}} \l( \Theta^{k_i}, \lam^{k_i - 1} \r) \\
        = \; & \frac{2}{m} \l( V_L^{k_i} - Y \r) + \lam_{2,L}^{k_i - 1} + \beta^{k_i - 1} \left( V_L^{k_i} - \sigma(U_L^{k_i}) \right) = \nabla_{V_L} \caL \l( \Theta^{k_i}, \lam^{k_i} \r).
    \end{split}
    \]    
    To prove that \eqref{eqvi.augl.lag} holds for $U_l$ ($l = 1,\dots,L$), we conduct an element-wise analysis of its associated nonsmooth components. To streamline this analysis, we introduce the following scalar notation for the $(r,s)$-th entries: 
\[
 \bar u = {U_l^{k_i}}_{(r,s)}, \quad \bar v = {V_l^{k_i}}_{(r,s)}, \quad
 \lambda^- = {\Lambda_{2,l}^{k_i-1}}_{(r,s)},\quad
 \lambda^+ = {\Lambda_{2,l}^{k_i}}_{(r,s)},\quad \beta^- = \beta^{k_i-1}. 
\]
    Then the multiplier update can be written as $\lambda^+ = \lambda^- + \beta^- ( \bar v - \sigma(\bar u))$. For fixed $\bar v$, $\lambda^-$, and $\beta^-$, define
\[
 q(u) = \lambda^- ( \bar v - \sigma(u) ) + \frac{\beta^-}{2} ( \bar v - \sigma(u) )^2.
\]
We need to prove 
\begin{equation}\label{element.u.sub.prove}
\partial q(\bar u) = \partial ( - \lambda^+ \sigma ) ( \bar u ).
\end{equation}
For ReLU, if $\bar u > 0$, then both sides equal $\{-\lambda^+\}$. If $\bar u < 0$, we have both sides equal $\{0\}$.
If $ \bar u = 0$, then the multiplier update becomes $ \lambda^+ = \lambda^- + \beta^- \bar v$ and we have the following:
\[
 q(u) = \lambda^- \bar v + \frac{\beta^-}{2} {\bar v}^2
       - \lambda^ + \max\{0,u\} + \frac{\beta^-}{2}(\max\{0,u\})^2 .
\]
The last term is continuously differentiable with derivative zero at
zero. Consequently,
\begin{align}
 \partial q(0) & =
 \begin{cases}
   \{0,-\lambda^+\}, & \lambda^+ \ge 0, \\
   [0,-\lambda^+], & \lambda^+ < 0,
 \end{cases} \\
 & = \partial ( - \lambda^+ \sigma_{\rm relu})(0).
\end{align}
Thus, \eqref{element.u.sub.prove} is proved. 
For leaky ReLU $\sigma_{\rm lrelu}(u)=\max\{\alpha u, u\}$ with $\alpha\in(0,1)$, by using the following:
\[
 \partial(-\eta \sigma_{\rm lrelu})(u)=
 \begin{cases}
   \{- \alpha \eta\}, & u < 0,\\
   \{- \eta\}, & u > 0,\\
   \{- \alpha \eta, - \eta\}, & u = 0,\ \eta\ge0,\\
   [- \alpha \eta, -\eta], & u = 0,\ \eta<0,
 \end{cases}
\]
we can prove \eqref{element.u.sub.prove} similarly: at $\bar u=0$,
the remainder $(\beta^-/2)\sigma_{\rm lrelu}(u)^2$ is continuously
differentiable with derivative zero, and away from zero both functions
are smooth.

The scalar formulas above also establish the joint closedness needed
when the multipliers vary. Indeed, if $(u_j,\eta_j,g_j)\to(u,\eta,g)$
and $g_j\in\partial(-\eta_j\sigma)(u_j)$, then
$g\in\partial(-\eta\sigma)(u)$. For $u\ne0$ this follows from the
fixed-sign formula. At $u=0$, when $\eta>0$ any limit is one of the
two endpoint slopes, whereas when $\eta<0$ every possible limit lies
in the stated interval; when $\eta=0$, the bound $|g_j|\le|\eta_j|$
forces $g=0$. Entrywise separability therefore shows that
\[
 (U,A)\longmapsto\partial_U\langle-A,\sigma(\cdot)\rangle(U)
\]
has a closed graph for both activations. This is the form of outer
semicontinuity used below when both the primal variables and the
multipliers converge. 
Entrywise separability and the smooth sum rule therefore give
\[
\begin{aligned}
 &\partial_{U_l}\caL_{\beta^{k_i-1}}
          \l( \Theta^{k_i},\Lambda^{k_i-1} \r)\\
 = & \; \lambda_2 U_l^{k_i}-\Lambda_{1,l}^{k_i-1} -\beta^{k_i-1}
       \l( W_l^{k_i}V_{l-1}^{k_i}+b_l^{k_i}1_m^{\T}-U_l^{k_i} \r) + \partial_{U_l} \l\langle -\Lambda_{2,l}^{k_i}, \sigma(\cdot) \r\rangle
                \l( U_l^{k_i} \r)\\
 = & \; \lambda_2 U_l^{k_i} - \Lambda_{1,l}^{k_i} + \partial_{U_l} \l\langle -\Lambda_{2,l}^{k_i}, \sigma(\cdot) \r\rangle \l( U_l^{k_i} \r) = \partial_{U_l} \caL \l( \Theta^{k_i}, \Lambda^{k_i} \r).
\end{aligned}
\]
Together, the block identities prove \eqref{eqvi.augl.lag}.
Consequently, by \eqref{eq:3.6} and the equivalence in \eqref{eqvi.augl.lag}, there exists a sequence $\{ \tau^{k_i} \}$ with $\|\tau^{k_i}\|_{\F} \le \epsilon^{k_i-1}\to0$ such that
\begin{equation}\label{err.lag}
 \tau^{k_i}\in\partial_\Theta\caL(\Theta^{k_i},\Lambda^{k_i}).
\end{equation}

    Next, we need to show the boundedness of Lagrange multipliers $ \{ \lam^{k_i} \}$ in \eqref{err.lag}. Define 
    \[
    \xi^i = \max_{1 \leq l \leq L} \l\{ \l\| \lam_{1,l}^{k_i} \r\|_{\max} , \l\| \lam_{2,l}^{k_i} \r\|_{\max} ,  \l\| \lam_{3,l}^{k_i} \r\|_{\max} \r\},
    \]
    and assume that $\{ \xi^i \}$ is unbounded. Passing to a further subsequence, we may assume that $\xi^i\to\infty$. Note that the normalized sequences $\{ \lam_{1,l}^{k_i} / \xi^i \}$, $\{ \lam_{2,l}^{k_i} / \xi^i \}$, $\{ \lam_{3,l}^{k_i} / \xi^i \}$ for $l = 1,\dots,L$ are bounded. Without loss of generality, we assume for $l = 1,\dots,L$
    \[
    \l\{ \lam_{1,l}^{k_i} / \xi^i \r\} \rightarrow \bar \lam_{1,l}, \quad \l\{\lam_{2,l}^{k_i} / \xi^i \r\} \rightarrow \bar \lam_{2,l}, \quad \l\{ \lam_{3,l}^{k_i} / \xi^i \r\} \rightarrow \bar \lam_{3,l}
    \]
    as $i \rightarrow \infty$ and have
    \begin{equation}\label{eq.lam.bar.0}
    \max_{1 \leq l \leq L} \l\{ \l\| \bar \lam_{1,l} \r\|_{\max} , \l\| \bar \lam_{2,l} \r\|_{\max} , \l\| \bar \lam_{3,l} \r\|_{\max} \r\} = 1.
    \end{equation}
    Dividing \eqref{err.lag} by $\xi^i$ and using
$\partial(tf)=t\partial f$ for $t>0$, we may absorb the positive
scaling into the multipliers in each nonsmooth term. The objective
gradients divided by $\xi^i$ vanish because $\Theta^{k_i}\to\Theta^\ast$
and $\xi^i\to\infty$. Taking $i\to\infty$, using the joint closedness
proved above, continuity of the smooth terms, and
$\tau^{k_i}/\xi^i\to0$, we obtain
    \begin{equation}\label{eq.partial.phat.bar}
        \begin{aligned}
            0 & = \bar \Lambda_{1,l} {V_{l-1}^\ast}^{\T} + \bar \Lambda_{3,l}, \quad && l \in \mathcal{I}, \\
            0 & = \bar \Lambda_{1,l} {V_{l-1}^\ast}^{\T} + {\bar \Lambda_{3,l}}^{\T}, \quad && l \in \mathcal{I}^c, \\
            0 & = \bar \Lambda_{1,l} 1_m, \quad && l = 1,\dots,L, \\
            0 & \in - \bar \Lambda_{1,l} + \partial_{U_l} \l\langle - \bar \lam_{2,l}, \sigma(\cdot) \r\rangle \l( U_l^\ast \r), \quad && l = 1,\dots,L, \\
            0 & = {W_{l+1}^\ast}^{\T} \bar \Lambda_{1,l+1} + \bar \Lambda_{2,l}, \quad && l = 1,\dots,L-1, \\
            0 & = \bar \Lambda_{2,L}.
        \end{aligned}
    \end{equation}
    Combining \eqref{eq.partial.phat.bar} and Lemma~\ref{le.partial.phat.0} yields that 
    \[
    \bar \lam_{1,l} = 0, \quad \bar \lam_{2,l} = 0, \quad \bar \lam_{3,l} = 0, 
    \]
    for $l = 1,\dots,L$, which contradicts \eqref{eq.lam.bar.0}. Therefore, the sequence $\{ \lam^{k_i}\}$ is bounded. Without loss of generality, we assume 
    \[
    \l\{ \lam_{1,l}^{k_i} \r\} \rightarrow \lam^\ast_{1,l}, \quad \l\{ \lam_{2,l}^{k_i} \r\}  \rightarrow \lam^\ast_{2,l}, \quad \l\{ \lam_{3,l}^{k_i} \r\} \rightarrow \lam^\ast_{3,l}, 
    \]
    for $l = 1,\dots,L$ as $i \rightarrow \infty$, after passing to
a further subsequence if necessary. Letting $i\to\infty$ in
\eqref{err.lag}, the same joint closedness, continuity of the smooth
gradients and Stiefel tangent projections, and $\tau^{k_i}\to0$ give
    \[
    0 \in \partial_\Theta \caL \l( \Theta^\ast, \lam^\ast \r).
    \]
    Together with the feasibility of $\Theta^\ast$, this proves that $\Theta^\ast$ is a KKT point of the problem~\eqref{prob.orth.three.sets.vari}.

\end{proof}

\section{Numerical Evaluation of Algorithms for OrthDFNs}\label{sec:experiments}

In this section, we evaluate the numerical performance of our proposed OrthALM and OrthADMM algorithms for training OrthDFNs. 
To demonstrate their practical effectiveness, we conduct comprehensive numerical experiments across three real-world benchmark datasets: the Extended Yale B face recognition database \cite{lee2005acquiring}, the MNIST dataset \cite{lecun1998mnist}, and the Fashion-MNIST dataset \cite{xiao2017fashion}.
All computations were conducted on a MacBook Pro equipped with an Apple M3 Pro processor (11-core CPU, 14-core GPU) and 18 GB of unified memory.

\subsection{Experimental Settings}\label{se9.1}

In our experiments, the weight matrices $W_l$ ($l = 1,\dots,L$) are initialized using random orthogonal initialization \cite{saxe2013exact}, while the bias vectors $b_l$ ($l = 1,\dots,L$) are initialized to zero.
In all reported experiments, both OrthALM and OrthADMM use leaky ReLU with $\alpha = 0.01$. The regularization parameters for OrthALM are $\lambda_1 = \lambda_2 = 10^{-7}$; OrthADMM uses the unregularized objective. For the proposed OrthALM, the parameters are chosen as $\mu_1 = \mu_2 = 10^{-4}$, $\eta_1 = 0.99$, $\eta_2 = 0.8$, $\eta_3 = 1.5$, $\eta_4 = 0.9$, $\epsilon^0 = 10^{-1}$, and $\beta^0 = 5 \cdot 10^{-3}$. For OrthADMM, the penalty parameter is fixed at $\hat{\beta} = 5 \cdot 10^{-3}$. In both algorithms, all initial multiplier matrices are initialized to zero.

We evaluate our proposed method against five state-of-the-art baselines:
\begin{enumerate}
    \item Approximate OrthDFNs via SVB \cite{li2019orthogonal};
    \item Riemannian SGD with QR retraction on the layer-wise Stiefel manifolds;
    \item OLM utilizing group based orthogonalization \cite{huang2018orthogonal}; 
    \item IALAM \cite{liu2023inexact} for solving leaky ReLU DFNs, which incorporates an $\ell_{2,1}$-norm weight regularizer into the squared loss objective;
    \item ADMM (sigmoid) \cite{zeng2021admm} for deep sigmoid networks, which incorporates a Frobenius norm weight regularizer into the squared loss objective. 
\end{enumerate}
The source code for all five baseline algorithms was obtained from the respective authors' GitHub repositories. The three SGD-based approaches were carried out using Python 3.10.18, whereas the remaining two baselines, alongside our proposed OrthADMM and OrthALM, were implemented in MATLAB R2024b. We specifically adopt the squared loss function for the three SGD-based methods. All these baseline algorithms are executed using their default settings, and the training process is terminated after a predefined maximum number of epochs. 
For the non-SGD-based methods, we define a single epoch as one iteration of OrthADMM or ADMM (sigmoid), and as one outer iteration for OrthALM or IALAM. Thus, equal epoch counts do not imply equal computational work; the reported times refer to these implementations and epoch budgets.

Across all evaluated datasets, we record the classification accuracy on both the training and test sets (denoted as ``TrainAcc'' and ``TestAcc'', respectively) and report the total computational CPU time in seconds (``Time''). Furthermore, we compute the empirical training loss (``TrainErr'') alongside three normalized constraint violations (ConsVio1, ConsVio2, and ConsVio3), defined for nonzero denominators as follows:
\begin{align*}
    \text{TrainErr} & = \frac{1}{m} \sum_{i = 1}^{m} \Big\| \sigma \Big( W_L \sigma \big( \dots \sigma ( W_1 x_i + b_1 ) \dots \big) + b_L \Big) - y_i \Big\|^2, \\[1ex]
    \text{ConsVio1}      & = \sum_{l=1}^{L}
      \frac{\|E_{1,l}(\Theta)\|_{\F}}{\|U_l\|_{\F}}, \qquad
    \text{ConsVio2} =\sum_{l=1}^{L}
      \frac{\|E_{2,l}(\Theta)\|_{\F}}{\|V_l\|_{\F}}, \\[1ex]
    \text{ConsVio3}      & =\sum_{l=1}^{L}
      \frac{\|E_{3,l}(\Theta)\|_{\F}}{\|Z_l\|_{\F}}.
\end{align*}

\subsection{Extended Yale B Face Recognition Database}\label{se9.2}

The Extended Yale B face recognition database \cite{georghiades2001from, lee2005acquiring} comprises a total of $2432$ images of $38$ subjects, with $64$ images per subject. During preprocessing, we resize the images from their original resolution of $192 \times 168$ pixels to $48 \times 42$ pixels. Each downsampled image is then vectorized into a $2016$-dimensional data vector. 
The $64$ images for each subject are randomly partitioned into equal-sized training and test sets, yielding exactly $1216$ samples for each. Following this split, we normalize the dataset by scaling all raw pixel intensities by a factor of $1/255$.

First, we evaluate the performance of our proposed OrthADMM and OrthALM frameworks across varying network architectures. Tables~\ref{ta_yaleb_admm} and \ref{ta_yaleb_alm} provide a detailed comparison of both algorithms after $200$ epochs, using networks of varying depths (up to four layers) and widths (number of neurons per layer).
Note that $L$ counts all weight-bearing layers, including the output layer. The hidden-width tuples list only hidden-layer widths; a dash denotes an unused hidden-layer position. The output-layer width is not included in these tuples.
For a fixed network depth, increasing the number of neurons predictably elevates the CPU time required by both methods, whereas the remaining performance metrics often exhibit only marginal improvements.

\begin{table}[htbp]
\caption{Experimental results of OrthADMM under different network settings.}\label{ta_yaleb_admm}
  \setlength{\tabcolsep}{4pt} 
  \begin{tabular}{ccccccccc}
    \toprule
       $L$ & Hidden widths & ConsVio1 & ConsVio2 & ConsVio3 & TrainErr & TrainAcc & TestAcc & Time (s) \\
       \midrule
       2 & {(50,-,-)} & {7.78e-3} & {4.31e-3} & {5.58e-3} & {5.51e-2} & {0.994} & {0.950} & {14.4} \\
       2 & {(100,-,-)} & {7.04e-3} & {3.11e-3} & {4.01e-3} & {4.07e-2} & {0.995} & {0.952} & {17.6} \\
       2 & {(200,-,-)} & {5.32e-3} & {2.66e-3} & {4.37e-3} & {3.06e-2} & {0.996} & {0.957} & {22.7} \\
       3 & {(50,25,-)} & {8.79e-3} & {5.94e-3} & {8.43e-2} & {1.70e-1} & {0.993} & {0.933} & {14.5} \\
       3 & {(100,50,-)} & {7.99e-3} & {5.67e-3} & {7.72e-3} & {2.06e-2} & {0.994} & {0.953} & {18.2} \\
       3 & {(200,100,-)} & {6.29e-3} & {4.27e-3} & {7.53e-3} & {2.86e-2} & {0.995} & {0.956} & {24.3} \\
       4 & {(50,25,25)} & {4.25e-2} & {5.15e-2} & {5.44e-1} & {2.09e-1} & {0.991} & {0.931} & {14.8} \\
       4 & {(100,50,50)} & {1.12e-2} & {8.33e-3} & {3.69e-2} & {4.43e-2} & {0.994} & {0.950} & {18.7} \\
       4 & {(200,100,100)} & {1.01e-2} & {7.65e-3} & {2.35e-2} & {3.15e-2} & {0.994} & {0.954} & {26.1} \\
       4 & {(500,250,250)} & {6.36e-3} & {4.92e-3} & {1.77e-2} & {2.04e-2} & {0.996} & {0.953} & {52.3} \\
     \botrule
  \end{tabular}
\end{table}

\begin{table}[htbp]
\caption{Experimental results of OrthALM under different network settings.}\label{ta_yaleb_alm}
  \setlength{\tabcolsep}{4pt} 
  \begin{tabular}{ccccccccc}
    \toprule
       $L$ & Hidden widths & ConsVio1 & ConsVio2 & ConsVio3 & TrainErr & TrainAcc & TestAcc & Time (s) \\
       \midrule
       2 & {(50,-,-)} & {6.04e-4} & {5.85e-4} & {1.07e-3} & {2.24e-2} & {0.997} & {0.942} & {97.0} \\
       2 & {(100,-,-)} & {5.74e-4} & {7.01e-5} & {2.81e-4} & {1.27e-2} & {0.998} & {0.951} & {117.8} \\
       2 & {(200,-,-)} & {4.03e-4} & {4.77e-5} & {1.94e-4} & {1.05e-2} & {0.997} & {0.947} & {151.8} \\
       3 & {(50,25,-)} & {2.49e-3} & {1.90e-3} & {4.88e-3} & {6.92e-2} & {0.993} & {0.924} & {96.5} \\
       3 & {(100,50,-)} & {5.04e-4} & {4.35e-4} & {1.59e-3} & {9.19e-3} & {0.998} & {0.951} & {121.9} \\
       3 & {(200,100,-)} & {9.18e-4} & {4.02e-4} & {7.44e-4} & {9.19e-3} & {0.998} & {0.947} & {165.2} \\
       4 & {(50,25,25)} & {2.45e-3} & {1.89e-3} & {8.11e-3} & {8.45e-2} & {0.993} & {0.923} & {97.4} \\
       4 & {(100,50,50)} & {1.29e-3} & {8.67e-4} & {1.34e-2} & {1.45e-2} & {0.998} & {0.947} & {127.8} \\
       4 & {(200,100,100)} & {1.22e-3} & {6.79e-4} & {3.13e-3} & {2.92e-2} & {0.978} & {0.928} & {176.4} \\
       4 & {(500,250,250)} & {8.91e-4} & {5.69e-4} & {2.01e-3} & {2.22e-2} & {0.983} & {0.930} & {335.7} \\
    \botrule
  \end{tabular}
\end{table}

\begin{figure}[htbp]
    \centering
    \begin{subfigure}[b]{0.325\textwidth}
        \centering
        \includegraphics[width=\linewidth]{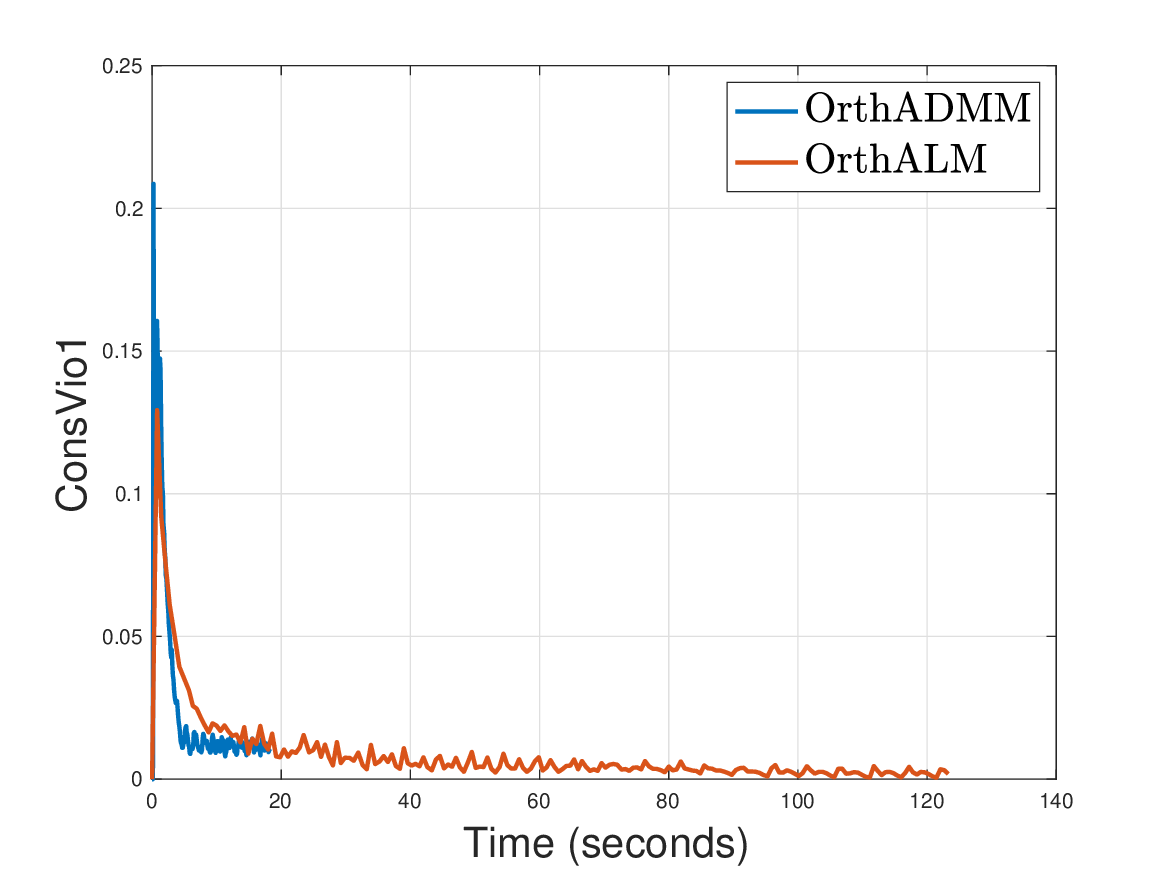}
        \caption{ConsVio1}
        \label{fig:yaleb_consvio_1}
    \end{subfigure}
    \hfill
    \begin{subfigure}[b]{0.325\textwidth}
        \centering
        \includegraphics[width=\linewidth]{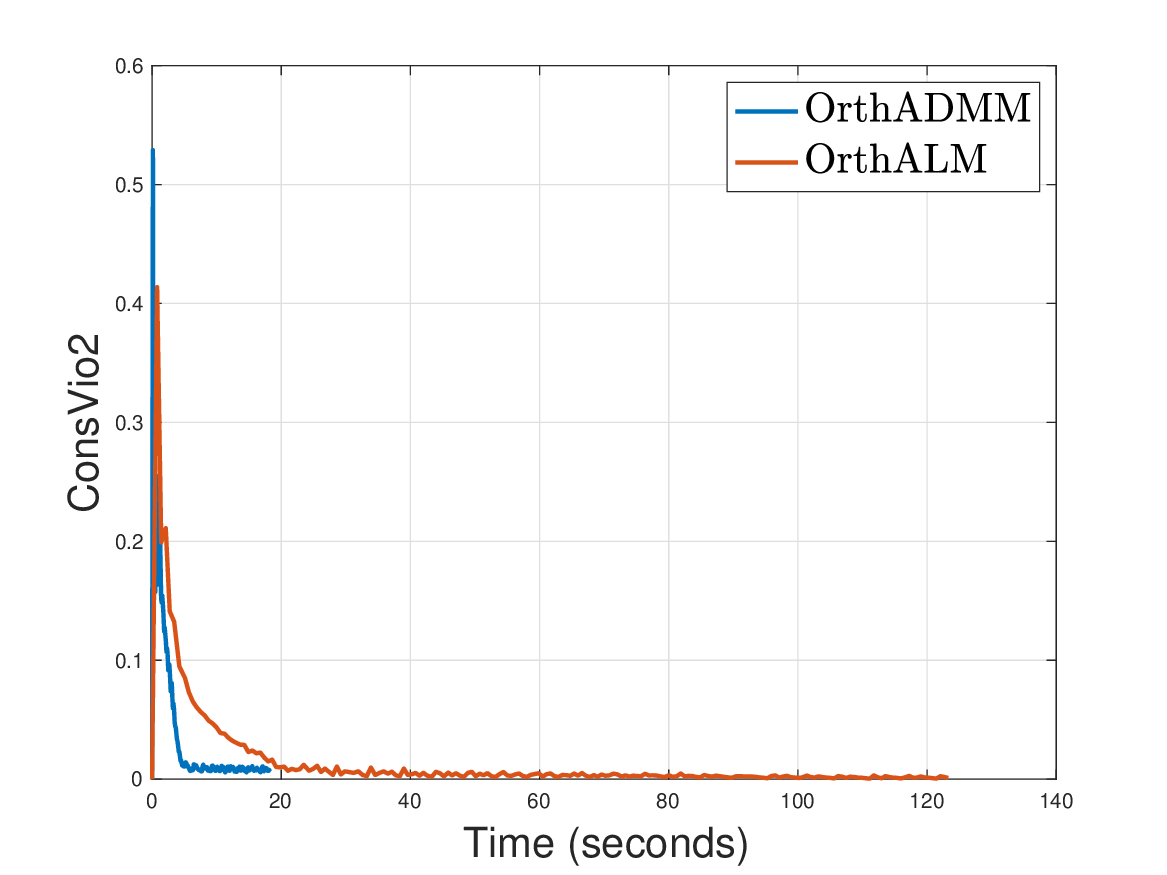}
        \caption{ConsVio2}
        \label{fig:yaleb_consvio_2}
    \end{subfigure}
    \hfill
    \begin{subfigure}[b]{0.325\textwidth}
        \centering
        \includegraphics[width=\linewidth]{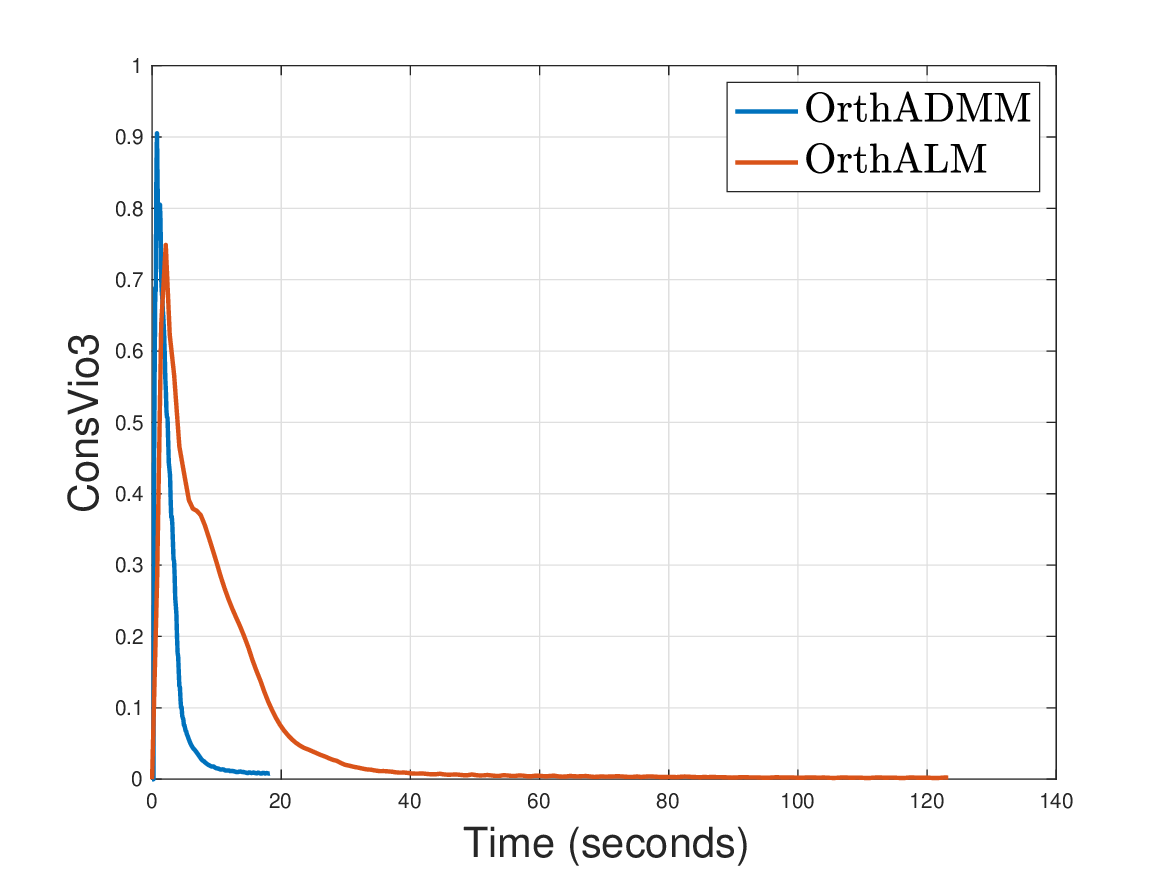}
        \caption{ConsVio3}
        \label{fig:yaleb_consvio_3}
    \end{subfigure}
    \caption{Evolution of the three normalized constraint violations against CPU time for OrthADMM and OrthALM}
    \label{fig:yaleb_consvio}
\end{figure}

For the architecture $L=3$, $n_1=100$, and $n_2=50$, Figure~\ref{fig:yaleb_consvio} plots the three normalized constraint violations against CPU time. OrthADMM shows a rapid initial decrease, while OrthALM attains smaller final values for all three residuals in the reported runs, as also shown in Tables~\ref{ta_yaleb_admm} and~\ref{ta_yaleb_alm}.

\begin{figure}
     \centering
     \begin{subfigure}[b]{0.325\textwidth}
         \centering
         \includegraphics[width=\textwidth]{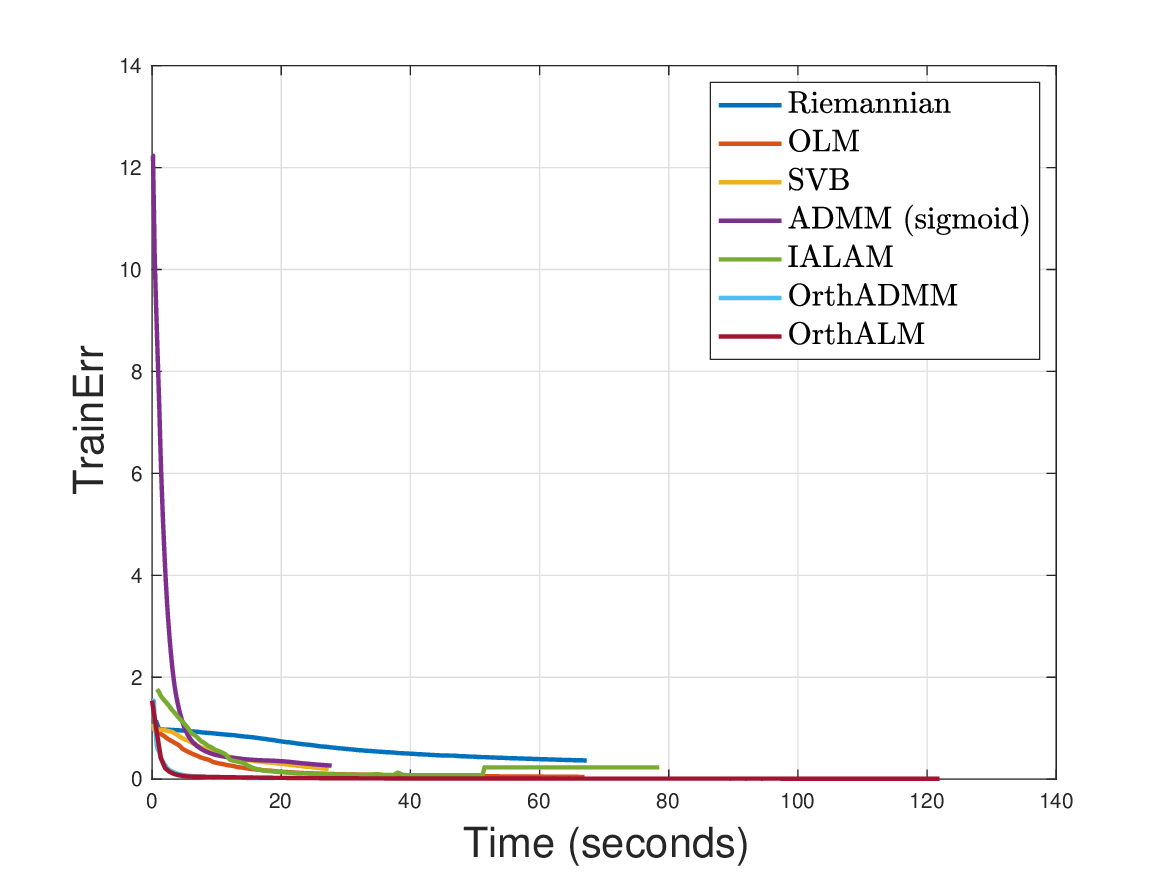}
         \caption{TrainErr vs Time}
         \label{9.1a}
     \end{subfigure}
     \hfill
     \begin{subfigure}[b]{0.325\textwidth}
         \centering
         \includegraphics[width=\textwidth]{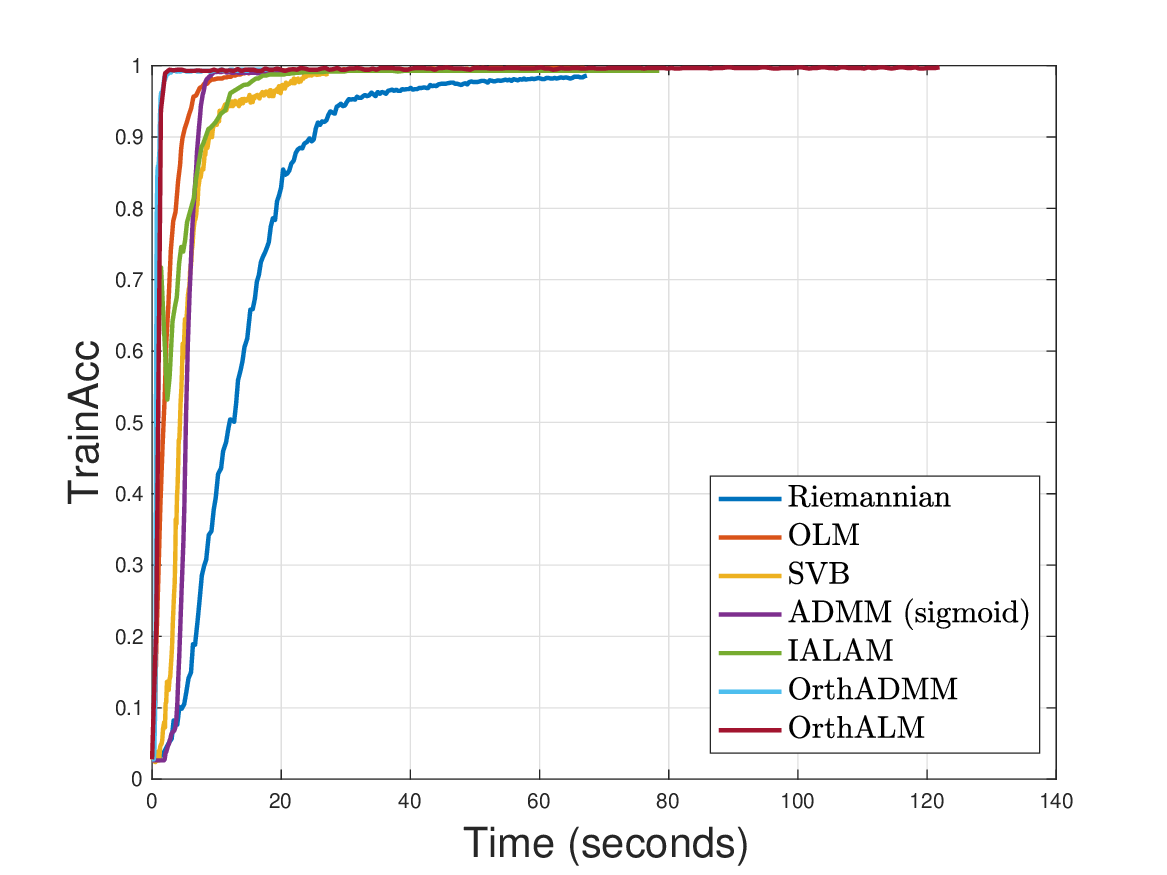}
         \caption{TrainAcc vs Time}
         \label{9.1c}
     \end{subfigure}
     \hfill
     \begin{subfigure}[b]{0.325\textwidth}
         \centering
         \includegraphics[width=\textwidth]{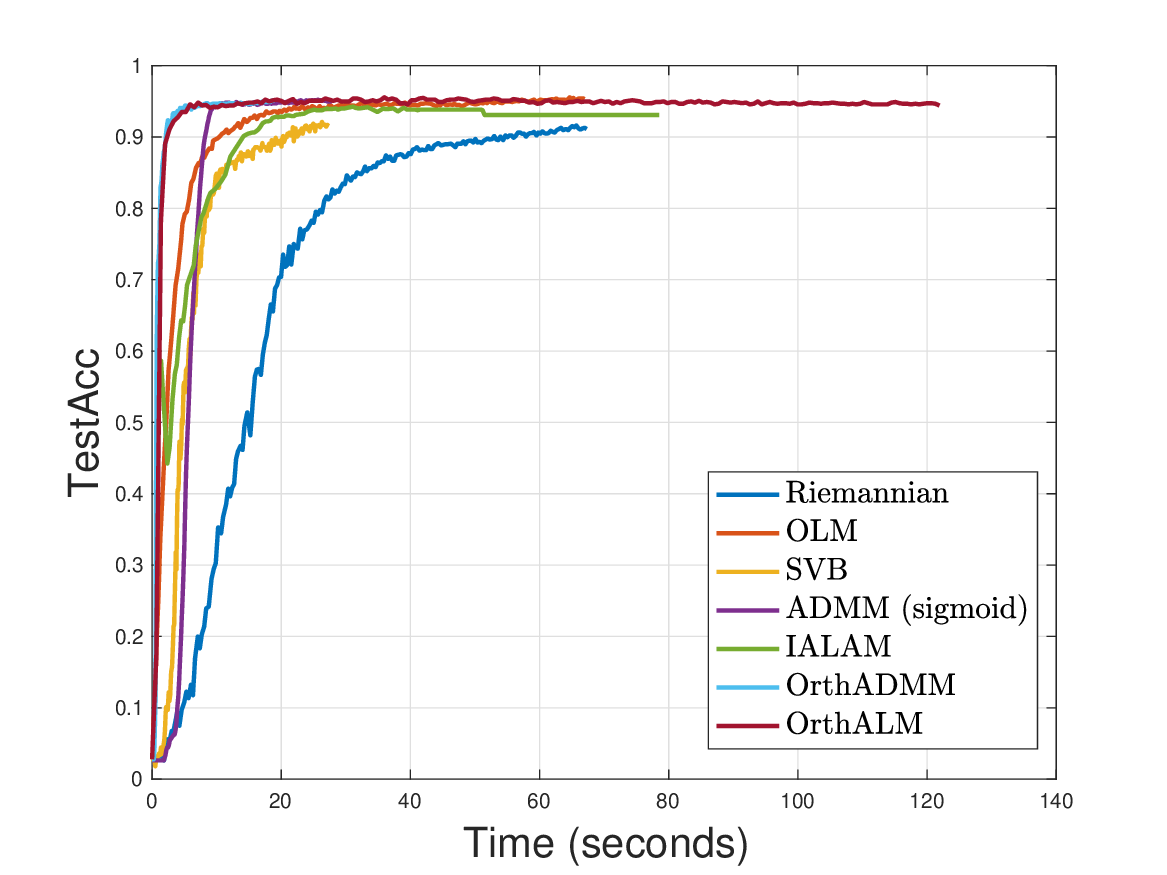}
         \caption{TestAcc vs Time}
         \label{9.1e}
     \end{subfigure}
     \hfill
     \begin{subfigure}[b]{0.325\textwidth}
         \centering
         \includegraphics[width=\textwidth]{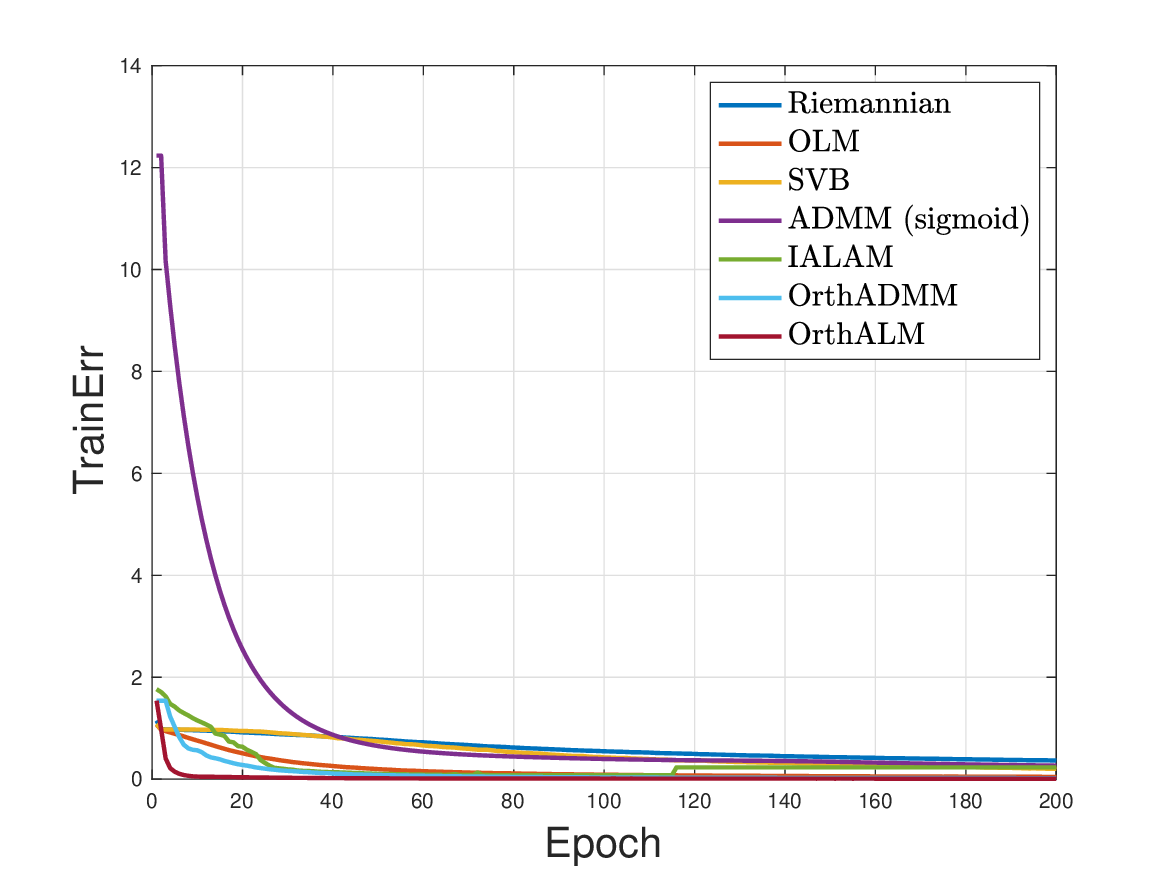}
         \caption{TrainErr vs Epoch}
         \label{9.1b}
     \end{subfigure}
     \hfill
     \begin{subfigure}[b]{0.325\textwidth}
         \centering
         \includegraphics[width=\textwidth]{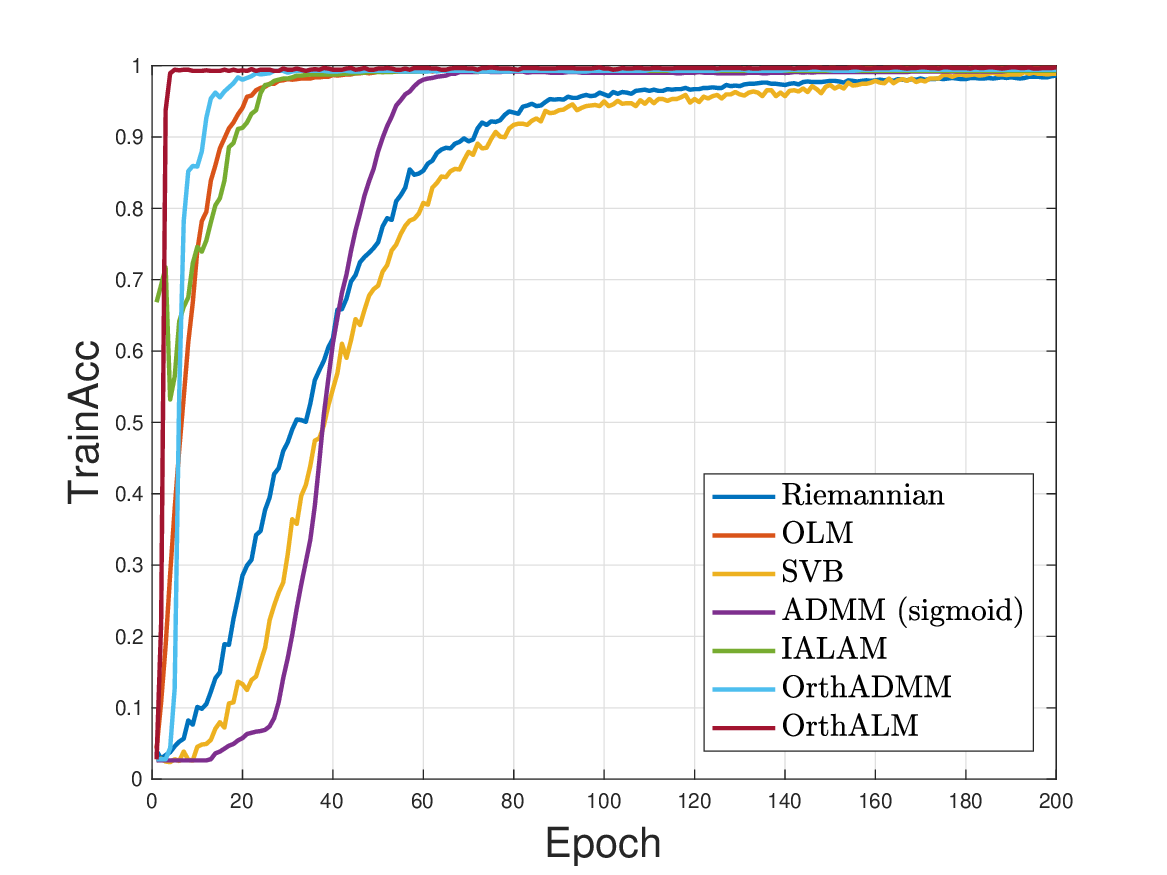}
         \caption{TrainAcc vs Epoch}
         \label{9.1d}
     \end{subfigure}
     \hfill
     \begin{subfigure}[b]{0.325\textwidth}
         \centering
         \includegraphics[width=\textwidth]{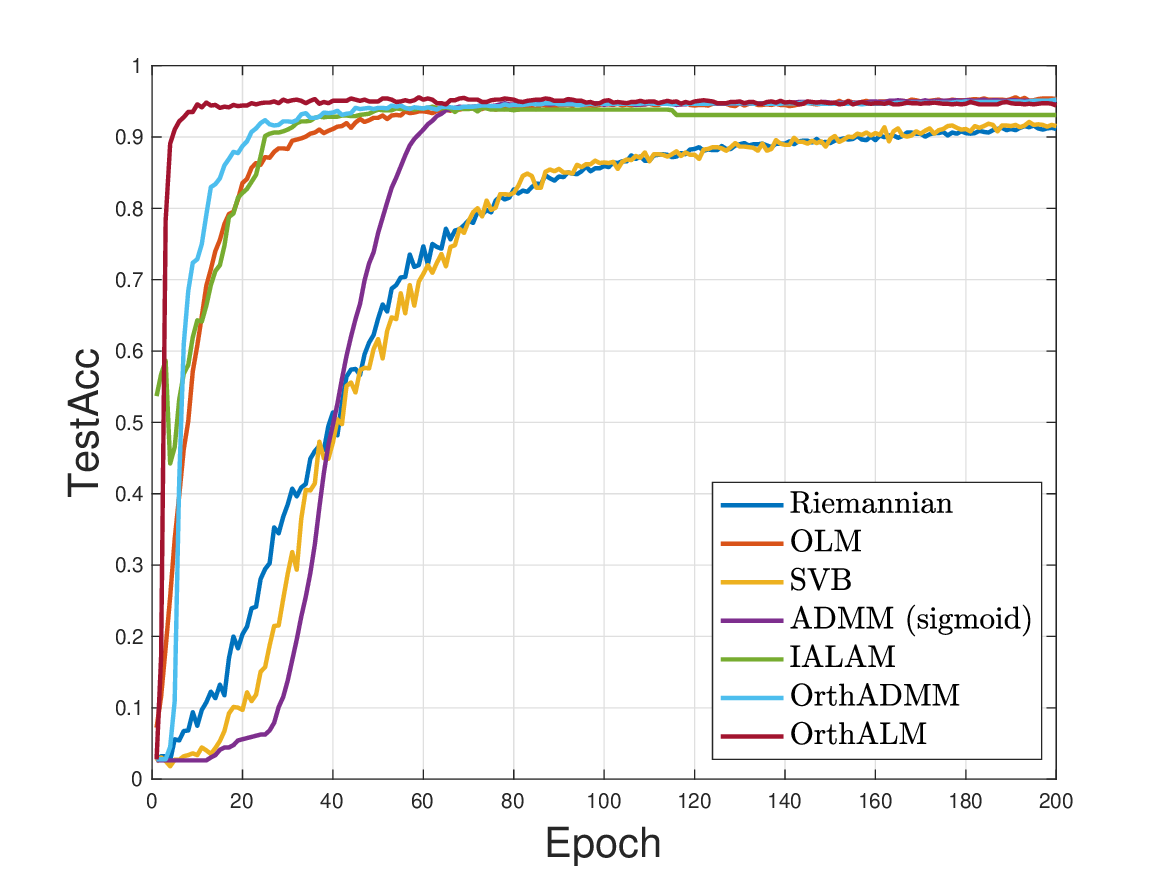}
         \caption{TestAcc vs Epoch}
         \label{9.1f}
     \end{subfigure}
        \caption{Comparison among OrthADMM, OrthALM, and five other existing approaches on the Extended Yale B dataset.}
        \label{fig 9.1}
\end{figure}

Under this same architectural configuration, we evaluate the performance of OrthADMM and OrthALM against the five baseline approaches from previous Section~\ref{se9.1} after $200$ training epochs.
We record the corresponding numerical results in Table~\ref{ta_yaleb_compare}, and the performance trajectories are illustrated in Figure~\ref{fig 9.1}. 
Based on these results, we make the following observations:

\begin{enumerate}
    \item Both OrthADMM and OrthALM achieve ``TrainAcc'' and ``TestAcc'' scores comparable to those of the state-of-the-art baselines.
    \item OrthALM yields the lowest overall ``TrainErr'', followed by OrthADMM. 
    \item Although OrthALM requires the most CPU time to complete the $200$ training epochs, OrthADMM stands out as the most computationally efficient algorithm among all evaluated algorithms.
    \item Figures~\ref{9.1a}, \ref{9.1c}, and~\ref{9.1e} show rapid initial improvement in the reported metrics for OrthADMM and OrthALM. The total times for the prescribed $200$ epochs are reported in Table~\ref{ta_yaleb_compare}; OrthADMM has the lowest total time, whereas OrthALM has the highest total time in this comparison.
\end{enumerate}

\begin{table}[htbp]
\caption{Performance comparison on the Extended Yale B dataset.}\label{ta_yaleb_compare}
  \begin{tabular}{ccccc}
    \toprule
       Method & TrainErr & TrainAcc & TestAcc & Time (s) \\
       \midrule
       {OrthADMM} & {2.06e-2} & {0.994} & {0.953} & {18.2} \\
       {OrthALM} & {9.19e-3} & {0.998} & {0.951} & {121.9} \\
       {OLM} & {4.50e-2} & {0.997} & {0.956} & {66.9} \\
       {SVB} & {2.03e-1} & {0.991} & {0.921} & {27.3} \\
       {Riemannian} & {3.66e-1} & {0.986} & {0.916} & {67.3} \\
       {IALAM} & {3.85e-2} & {0.993} & {0.938} & {78.5} \\
       {ADMM(sigmoid)} & {2.66e-1} & {0.993} & {0.952} & {27.8} \\
    \botrule
  \end{tabular}
\end{table}

\subsection{MNIST Dataset}\label{se9.3}

The MNIST dataset \cite{lecun1998gradient} is a widely used benchmark for handwritten digit recognition. It comprises $70{,}000$ grayscale images across $10$ classes ($n_L = 10$), pre-partitioned into a training set of $60{,}000$ images and a test set of $10{,}000$ images. Prior to training, each $28 \times 28$ pixel image is vectorized into a $784$-dimensional feature vector ($n_0 = 784$), and its raw pixel intensities are then also normalized by a factor of $1/255$.

To evaluate the scalability and sample efficiency of our proposed OrthADMM and OrthALM methods across different training set sizes, we randomly sample $m \in \{500, 1000, 2000, 5000\}$ training instances from the full training set. The size of the corresponding test set, denoted as $m_{\text{test}}$, is fixed at $0.4m$. Here in practice, we randomly pick up data entries from each class of MNIST under uniform distribution to ensure a balanced representation.

\begin{figure}
     \centering
     \begin{subfigure}[b]{0.325\textwidth}
         \centering
         \includegraphics[width=\textwidth]{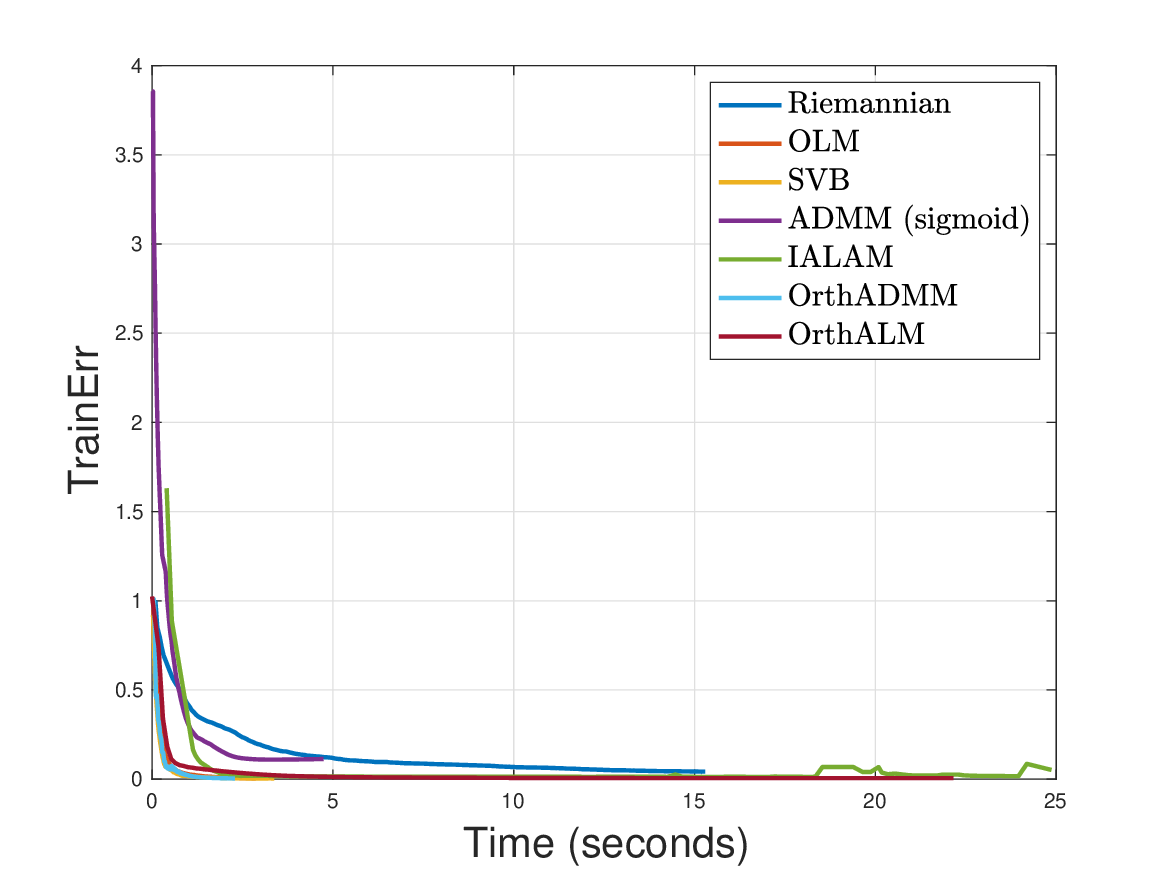}
         \caption{TrainErr}
         \label{9.2a}
     \end{subfigure}
     \hfill
     \begin{subfigure}[b]{0.325\textwidth}
         \centering
         \includegraphics[width=\textwidth]{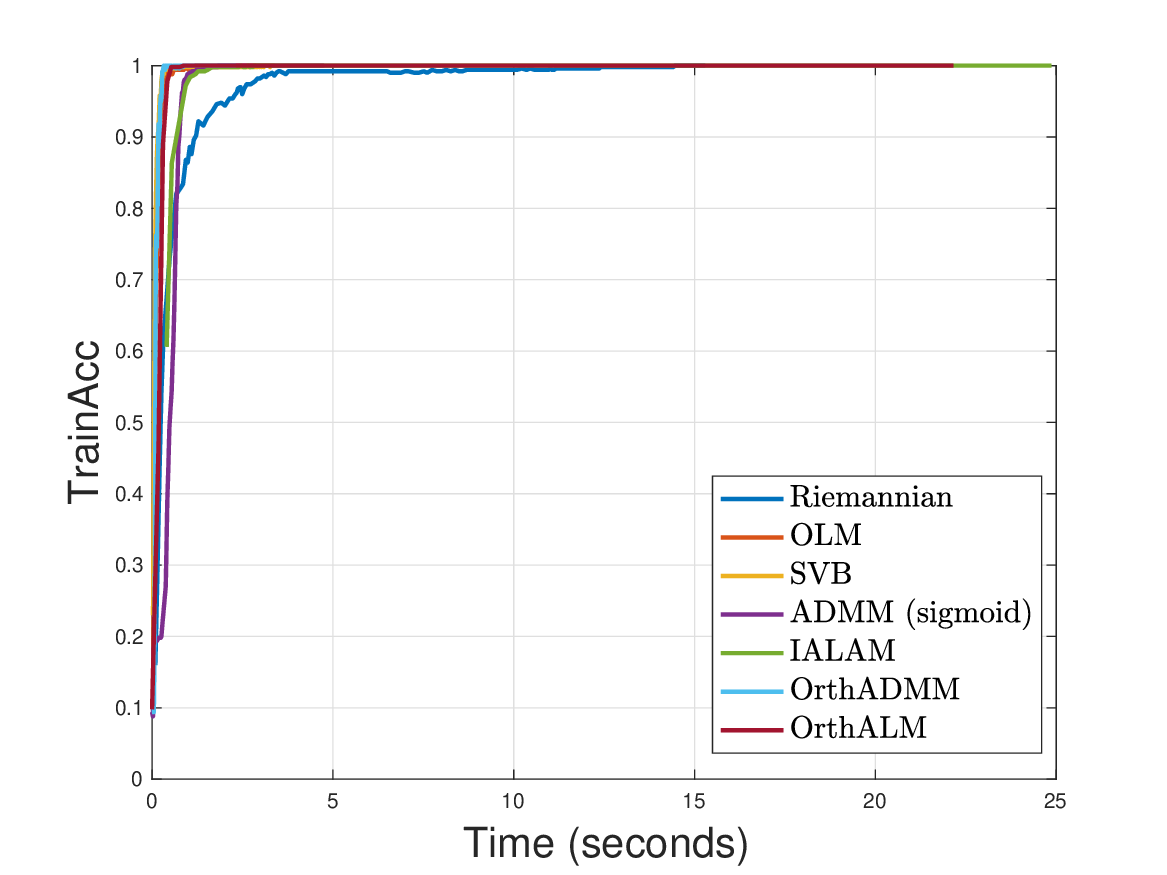}
         \caption{TrainAcc}
         \label{9.2b}
     \end{subfigure}
     \hfill
     \begin{subfigure}[b]{0.325\textwidth}
         \centering
         \includegraphics[width=\textwidth]{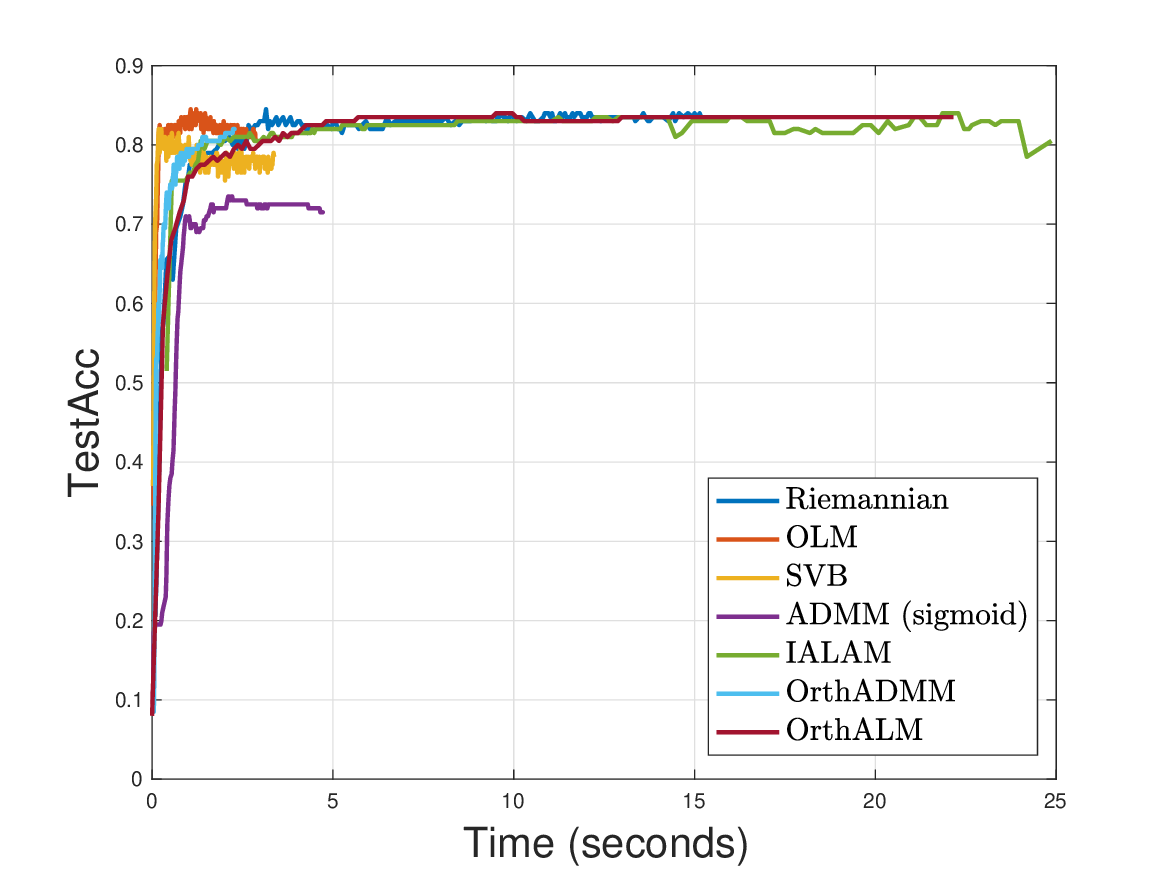}
         \caption{TestAcc}
         \label{9.2c}
     \end{subfigure}
     \hfill
     \begin{subfigure}[b]{0.325\textwidth}
         \centering
         \includegraphics[width=\textwidth]{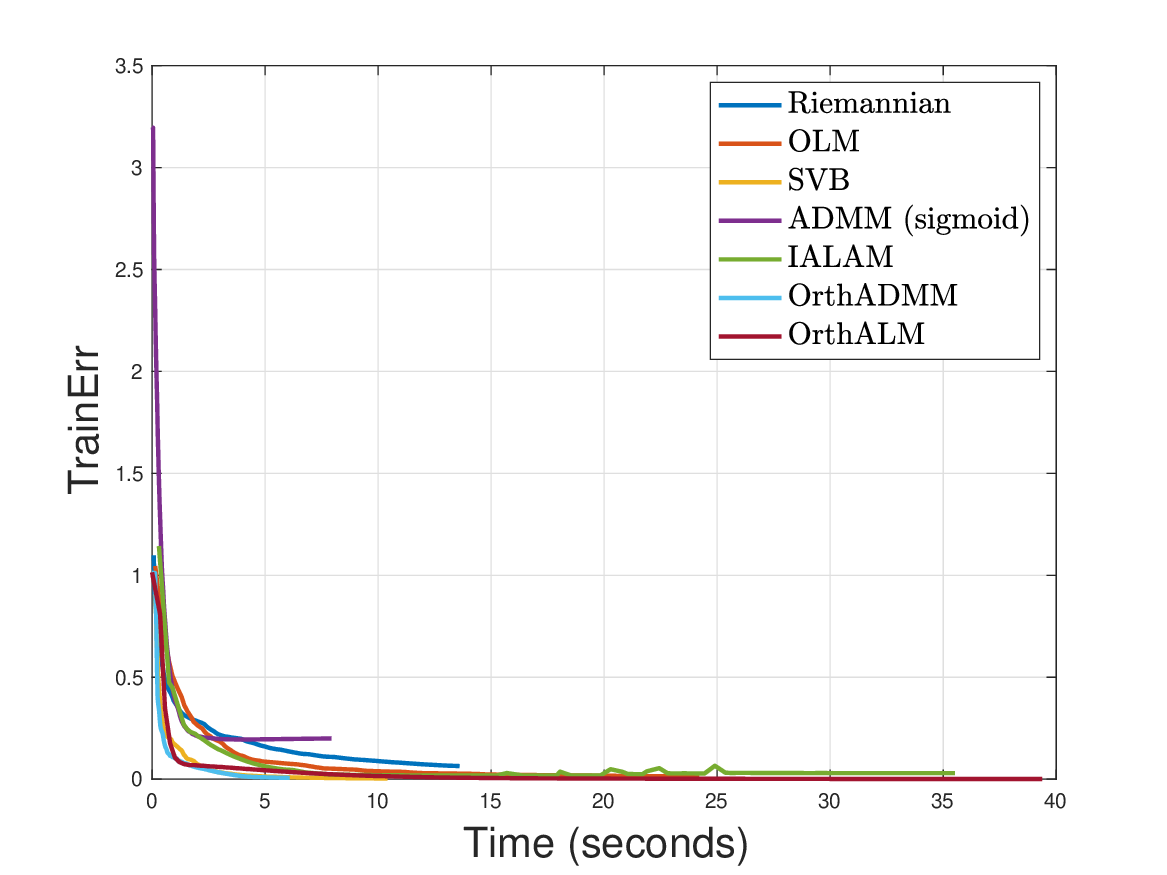}
         \caption{TrainErr}
         \label{9.2d}
     \end{subfigure}
     \hfill
     \begin{subfigure}[b]{0.325\textwidth}
         \centering
         \includegraphics[width=\textwidth]{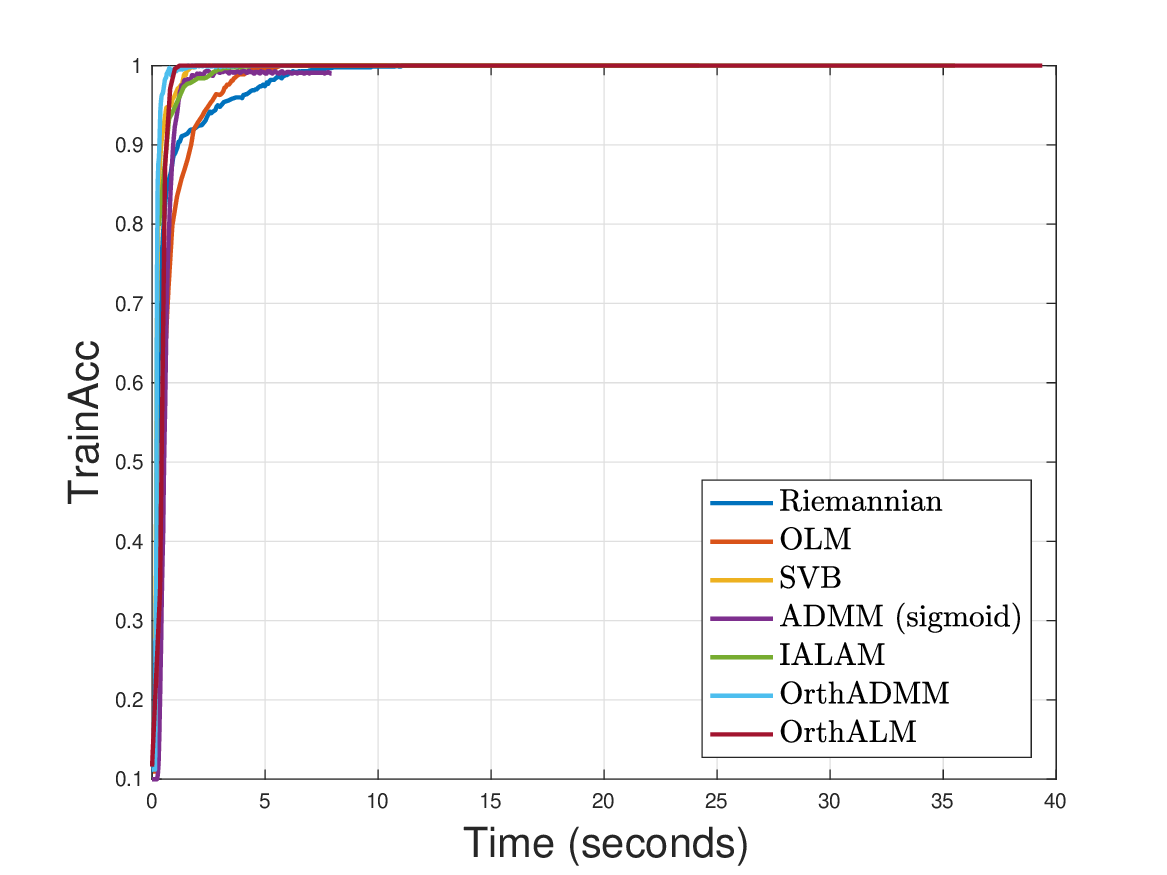}
         \caption{TrainAcc}
         \label{9.2e}
     \end{subfigure}
     \hfill
     \begin{subfigure}[b]{0.325\textwidth}
         \centering
         \includegraphics[width=\textwidth]{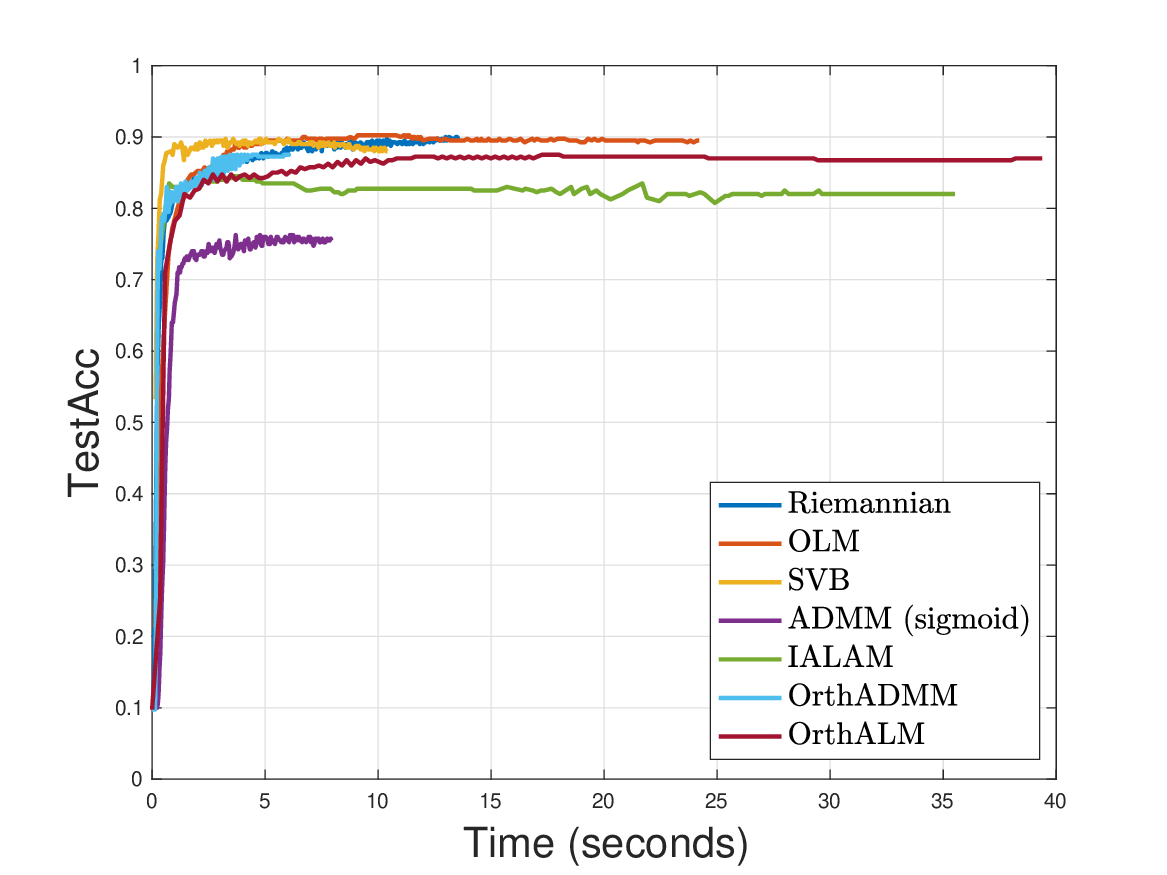}
         \caption{TestAcc}
         \label{9.2f}
     \end{subfigure}
     \hfill
     \begin{subfigure}[b]{0.325\textwidth}
         \centering
         \includegraphics[width=\textwidth]{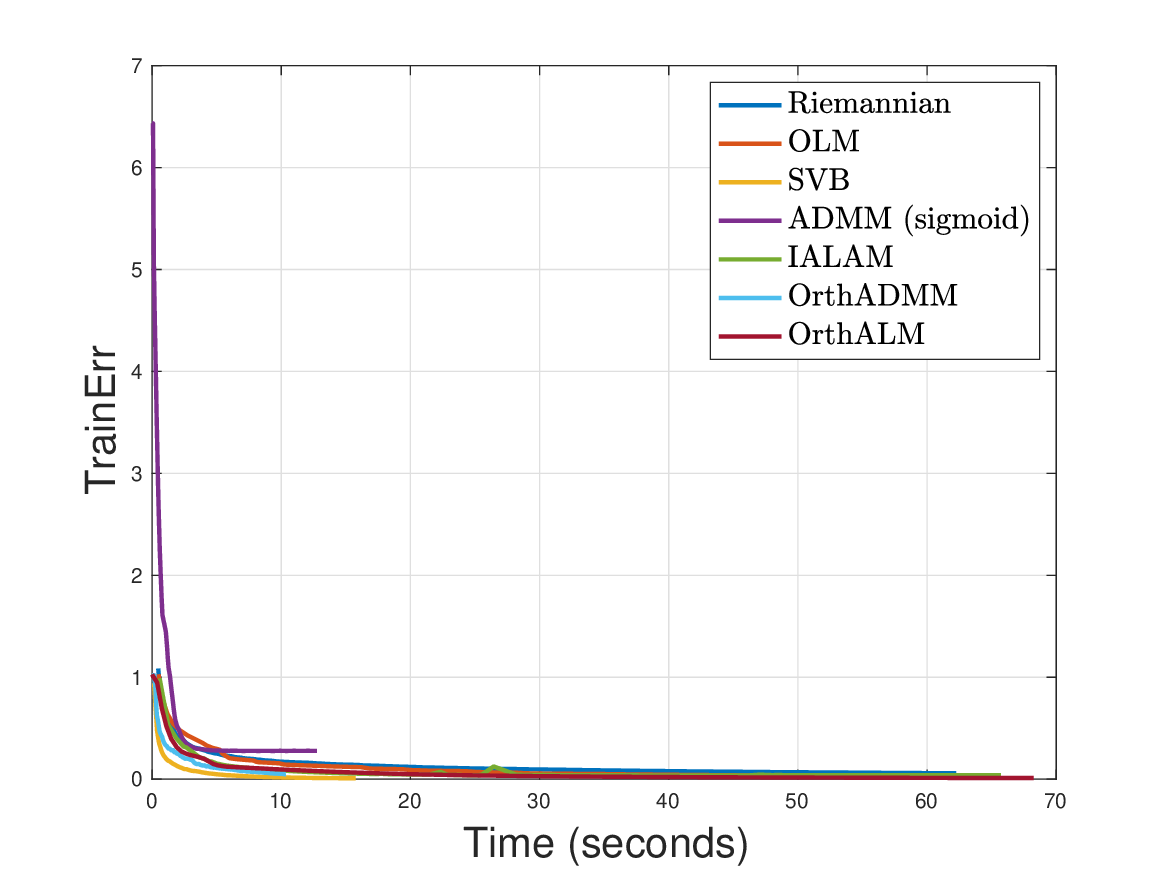}
         \caption{TrainErr}
         \label{9.2g}
     \end{subfigure}
     \hfill
     \begin{subfigure}[b]{0.325\textwidth}
         \centering
         \includegraphics[width=\textwidth]{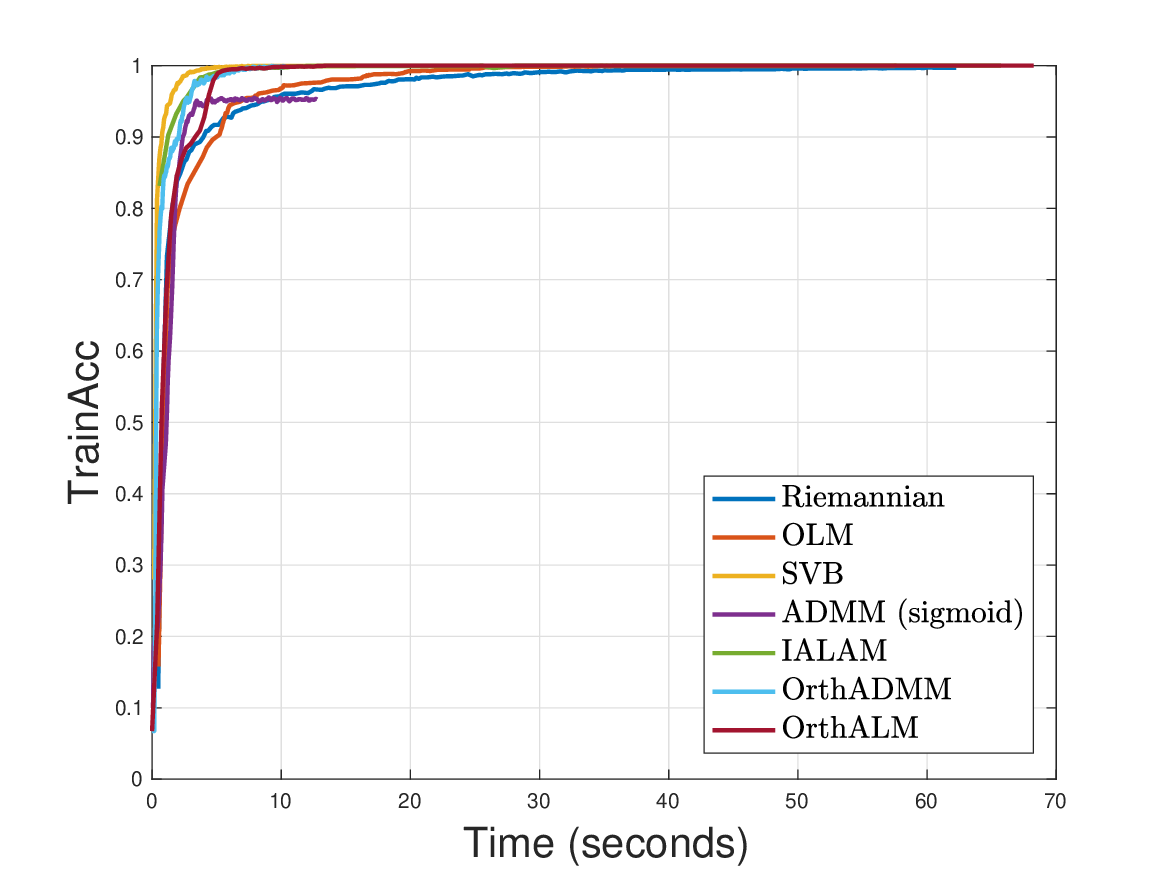}
         \caption{TrainAcc}
         \label{9.2h}
     \end{subfigure}
     \hfill
     \begin{subfigure}[b]{0.325\textwidth}
         \centering
         \includegraphics[width=\textwidth]{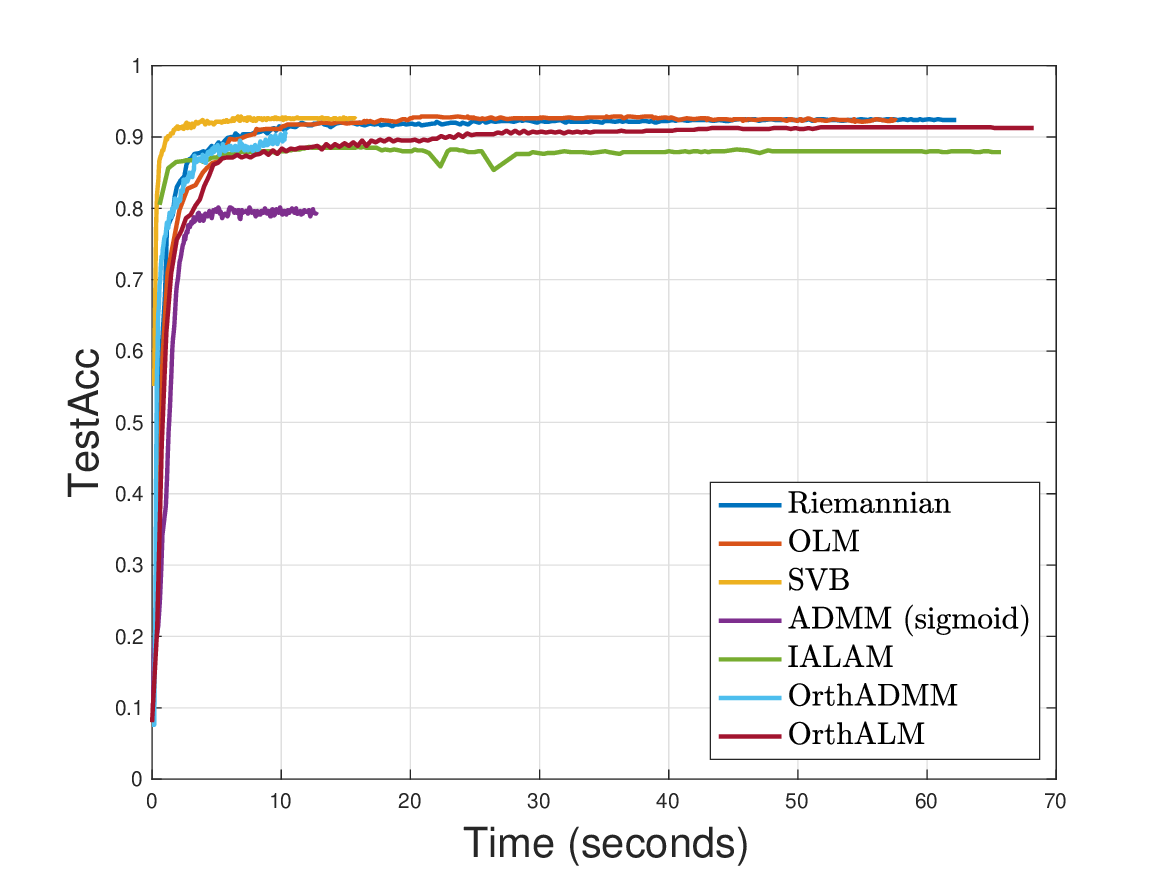}
         \caption{TestAcc}
         \label{9.2i}
     \end{subfigure}
     \hfill
     \begin{subfigure}[b]{0.325\textwidth}
         \centering
         \includegraphics[width=\textwidth]{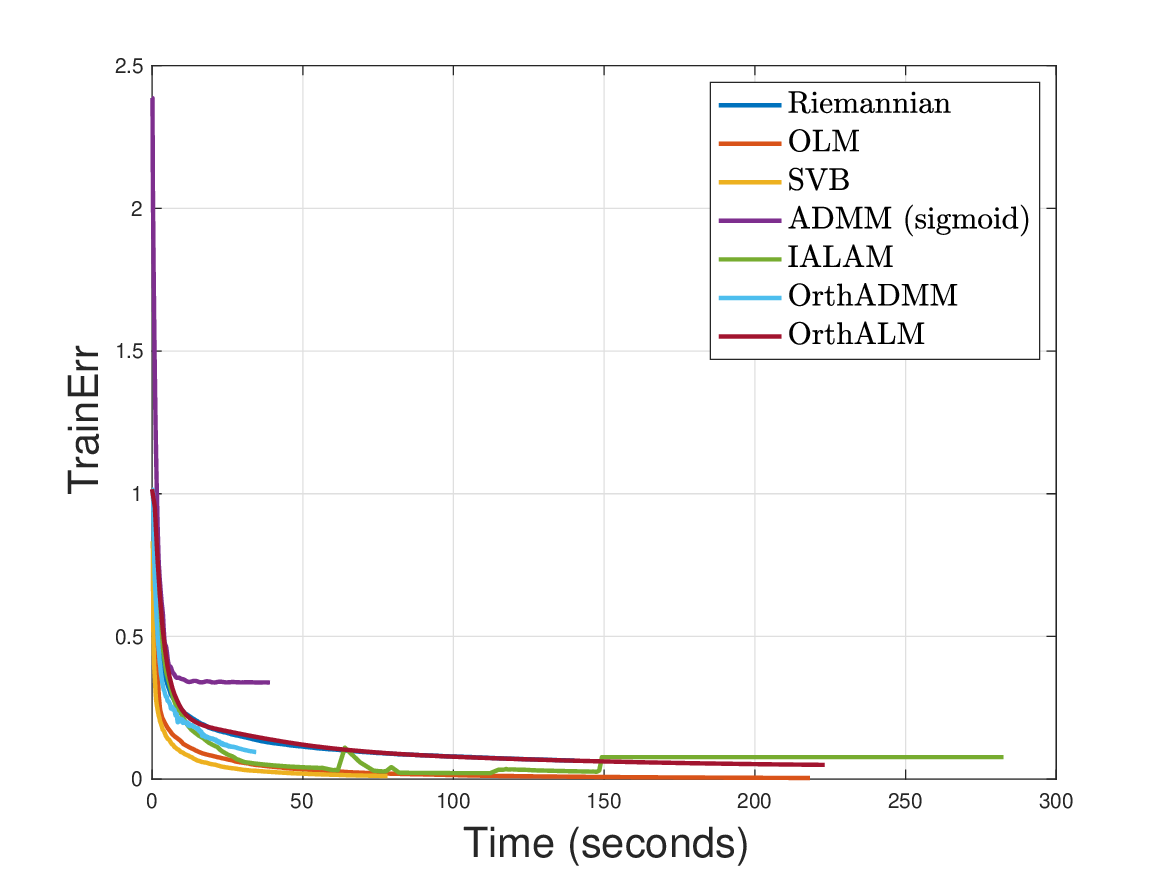}
         \caption{TrainErr}
         \label{9.2j}
     \end{subfigure}
     \hfill
     \begin{subfigure}[b]{0.325\textwidth}
         \centering
         \includegraphics[width=\textwidth]{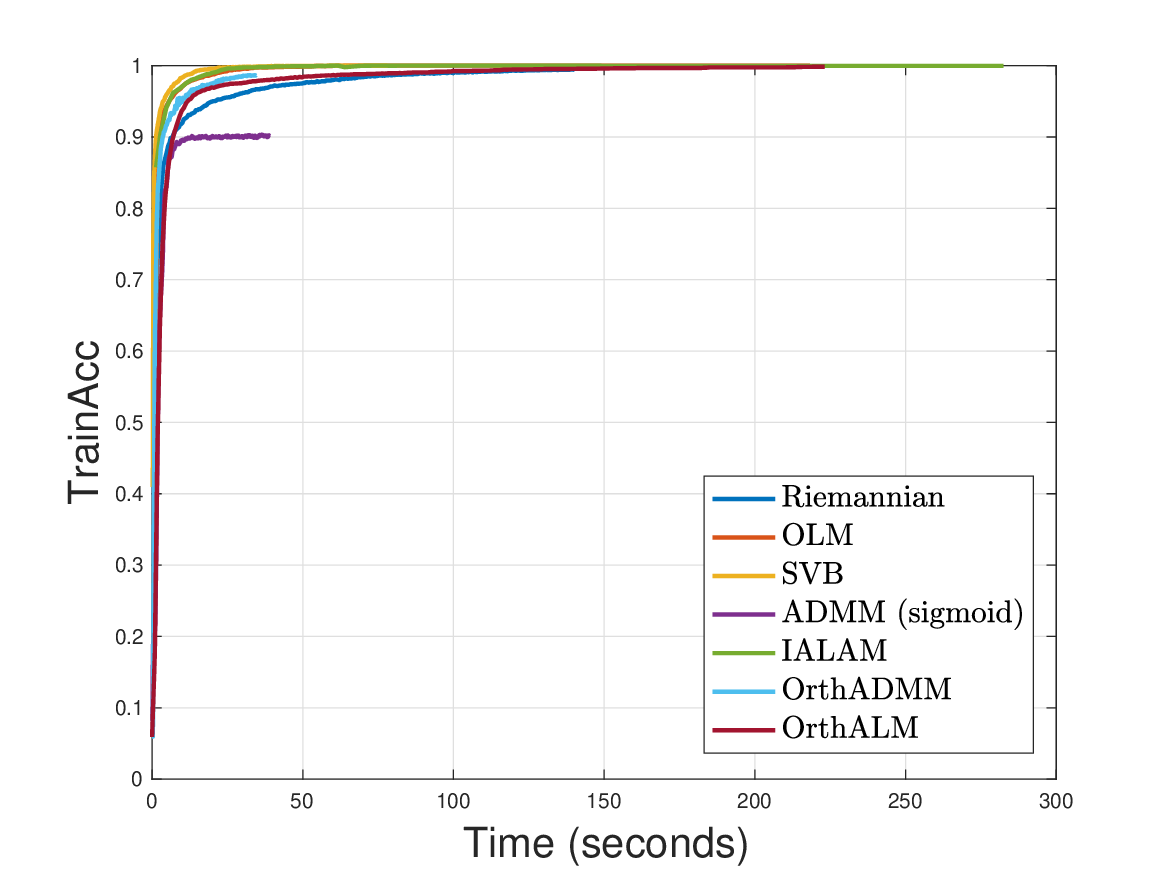}
         \caption{TrainAcc}
         \label{9.2k}
     \end{subfigure}
     \hfill
     \begin{subfigure}[b]{0.325\textwidth}
         \centering
         \includegraphics[width=\textwidth]{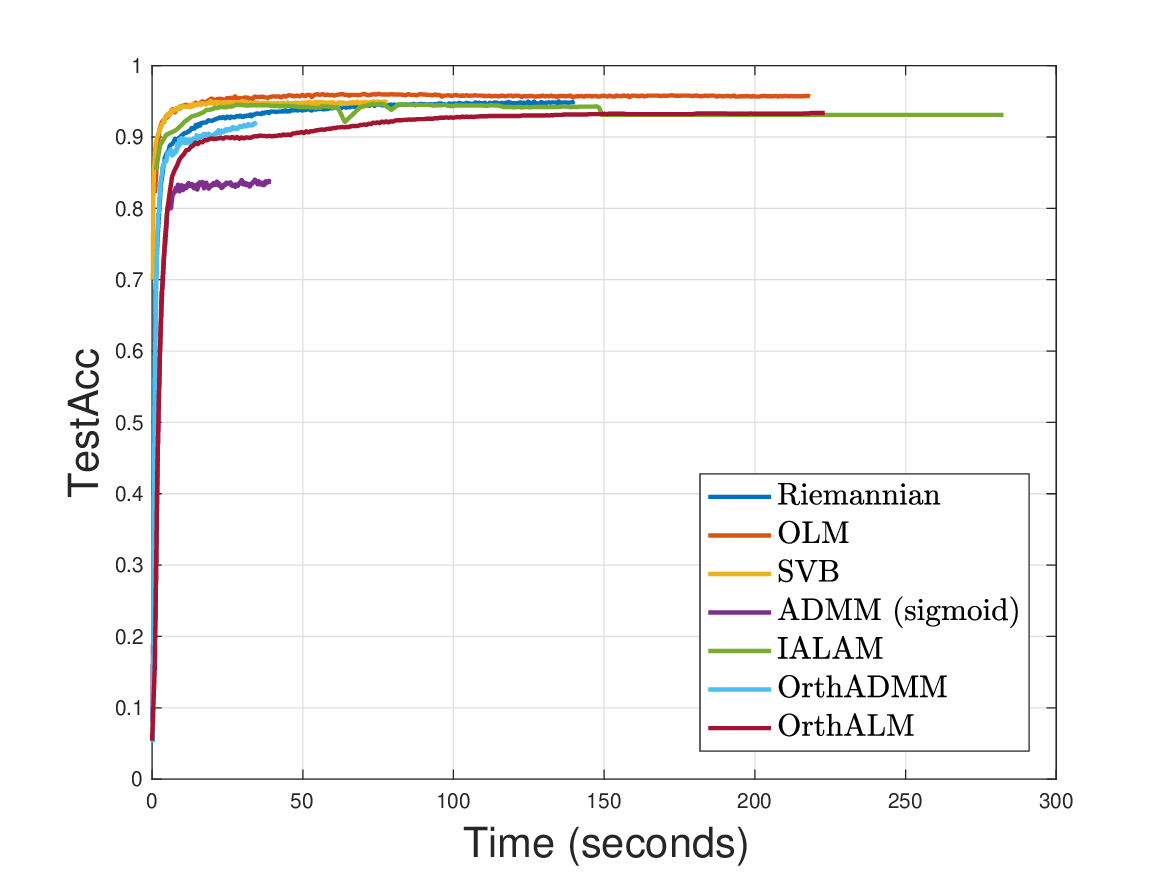}
         \caption{TestAcc}
         \label{9.2l}
     \end{subfigure}
        \caption{Comparison among OrthADMM, OrthALM, and five other existing approaches on the MNIST dataset across different training sample sizes: $m = 500$ (first row, (a)--(c)), $m = 1000$ (second row, (d)--(f)), $m = 2000$ (third row, (g)--(i)), and $m = 5000$ (fourth row, (j)--(l)).}
        \label{fig:mnist}
\end{figure}

We evaluate OrthADMM and OrthALM against the five previously described existing methods in Section~\ref{se9.1}. All these seven algorithms are executed until a predefined maximum number of epochs, denoted as $\text{epoch}_{\max}$, is reached. The experiments are conducted across four distinct network and dataset configurations, specified as follows:
\begin{equation}\label{mnist.network.para}
\begin{aligned}
    & m = 500,  && m_{\text{test}} = 200,  && L = 3, \quad n_1 = 25,  \quad n_2 = 20,  && \text{epoch}_{\max} = 200; \\
    & m = 1000, && m_{\text{test}} = 400,  && L = 3, \quad n_1 = 100, \quad n_2 = 50,  && \text{epoch}_{\max} = 200; \\
    & m = 2000, && m_{\text{test}} = 800,  && L = 3, \quad n_1 = 150, \quad n_2 = 75,  && \text{epoch}_{\max} = 200; \\
    & m = 5000, && m_{\text{test}} = 2000, && L = 3, \quad n_1 = 200, \quad n_2 = 100, && \text{epoch}_{\max} = 300.
\end{aligned}
\end{equation}
The numerical results are summarized in Table~\ref{tab:mnist}, and the corresponding convergence plots of the training loss, training accuracy, and test accuracy relative to CPU time are illustrated in Figure~\ref{fig:mnist}.
Based on these numerical results, we draw the following observations:
\begin{enumerate}

    \item As the training sample size $m$ scales from $500$ to $5000$, the test accuracy for all seven methods is ultimately higher at $m = 5000$, although the overall progression includes intermediate declines.
    However, the training accuracy for these methods either remains stable near $1.00$ or experiences a slight decline.
    
    \item Both OrthADMM and OrthALM achieve a $100\%$ training accuracy for $m \in \{500, 1000, 2000\}$ and maintain near-perfect accuracy ($0.987$ and $0.998$, respectively) at $m = 5000$. Their test accuracies are competitive with IALAM and the SGD-based methods for OrthDFNs, while consistently outperforming the ADMM (sigmoid) baseline across all evaluated sample sizes.
    
    \item Regarding the training error, OrthALM achieves a lower value than OrthADMM for $m \ge 1000$, while the ADMM (sigmoid) method exhibits the poorest overall performance.
    
    \item OrthADMM is the most computationally efficient method among all seven algorithms across every evaluated sample size $m \in \{500, 1000, 2000, 5000\}$. Meanwhile, the CPU time required by OrthALM is comparable to that of IALAM.
    
    \item As illustrated across all four rows of Figure~\ref{fig:mnist}, OrthADMM and OrthALM exhibit a rapid initial reduction in training loss alongside improvements in training and test accuracy. Table~\ref{tab:mnist} reports the total times for the prescribed epoch budgets.

    \item A direct comparison of the two proposed methods reveals that OrthALM generally achieves higher training and test accuracies, alongside a lower training error, than OrthADMM. However, this performance advantage comes at the expense of a higher computational cost and significantly greater CPU time.
    
\end{enumerate}

\begin{table}[htbp]
\caption{Performance comparison of our proposed OrthADMM and OrthALM against five state-of-the-art baselines on the MNIST dataset.}\label{tab:mnist}
\begin{tabular}{llcccc}
\toprule
Method & Sample Size & TrainErr & TrainAcc & TestAcc & Time (s) \\
\midrule
OrthADMM & $\mathbf{m = 500}$  & 3.97e-3 & 1.00  & 0.820 & 2.3   \\
         & $\mathbf{m = 1000}$ & 6.37e-3 & 1.00  & 0.875 & 6.1   \\
         & $\mathbf{m = 2000}$ & 4.77e-2 & 1.00  & 0.906 & 10.4  \\
         & $\mathbf{m = 5000}$ & 9.50e-2 & 0.987 & 0.921 & 34.5  \\
\midrule
OrthALM  & $\mathbf{m = 500}$  & 5.88e-3 & 1.00  & 0.835 & 22.1  \\
         & $\mathbf{m = 1000}$ & 7.74e-4 & 1.00  & 0.873 & 39.4  \\
         & $\mathbf{m = 2000}$ & 1.08e-2 & 1.00  & 0.914 & 68.2  \\
         & $\mathbf{m = 5000}$ & 5.00e-2 & 0.998 & 0.934 & 223.1 \\
\midrule
OLM      & $\mathbf{m = 500}$  & 5.33e-3 & 1.00  & 0.835 & 2.9   \\
         & $\mathbf{m = 1000}$ & 1.31e-2 & 1.00  & 0.895 & 24.2  \\
         & $\mathbf{m = 2000}$ & 2.64e-2 & 1.00  & 0.923 & 57.6  \\
         & $\mathbf{m = 5000}$ & 4.15e-3 & 1.00  & 0.957 & 218.3 \\
\midrule
SVB      & $\mathbf{m = 500}$  & 3.05e-3 & 1.00  & 0.795 & 3.4   \\
         & $\mathbf{m = 1000}$ & 3.28e-3 & 1.00  & 0.880 & 10.4  \\
         & $\mathbf{m = 2000}$ & 8.50e-3 & 1.00  & 0.925 & 15.8  \\
         & $\mathbf{m = 5000}$ & 1.10e-2 & 1.00  & 0.949 & 78.1  \\
\midrule
Riemannian & $\mathbf{m = 500}$  & 4.20e-2 & 1.00  & 0.840 & 15.3  \\
           & $\mathbf{m = 1000}$ & 6.41e-2 & 1.00  & 0.898 & 13.6  \\
           & $\mathbf{m = 2000}$ & 5.80e-2 & 0.997 & 0.924 & 62.2  \\
           & $\mathbf{m = 5000}$ & 6.53e-2 & 0.994 & 0.950 & 140.1 \\
\midrule
IALAM    & $\mathbf{m = 500}$  & 6.24e-3 & 1.00  & 0.840 & 24.9  \\
         & $\mathbf{m = 1000}$ & 8.76e-3 & 1.00  & 0.835 & 35.5  \\
         & $\mathbf{m = 2000}$ & 1.79e-2 & 1.00  & 0.883 & 65.7  \\
         & $\mathbf{m = 5000}$ & 1.06e-2 & 1.00  & 0.945 & 282.5 \\
\midrule
ADMM (sigmoid) & $\mathbf{m = 500}$  & 1.09e-1 & 1.00  & 0.730 & 4.7   \\
               & $\mathbf{m = 1000}$ & 1.95e-1 & 0.994 & 0.763 & 7.9   \\
               & $\mathbf{m = 2000}$ & 2.77e-1 & 0.956 & 0.801 & 12.8  \\
               & $\mathbf{m = 5000}$ & 3.38e-1 & 0.903 & 0.840 & 39.1  \\
\botrule
\end{tabular}
\end{table}

\subsection{Fashion-MNIST Dataset}\label{se9.4}

Our final evaluation addresses image classification using the Fashion-MNIST dataset \cite{xiao2017fashion}. Similar to MNIST, it comprises $70{,}000$ grayscale images of size $28 \times 28$---pre-partitioned into $60{,}000$ training images and $10{,}000$ test images---distributed across $10$ classes ($n_L = 10$). These categories include T-shirts/tops, trousers, pullovers, dresses, coats, sandals, shirts, sneakers, bags, and ankle boots, with each assigned an integer label from $0$ to $9$.
Following the preprocessing pipeline detailed in Section~\ref{se9.3}, each image is vectorized into a $784$-dimensional feature vector ($n_0 = 784$) and normalized by a factor of $1/255$.

\begin{figure}
     \centering
     \begin{subfigure}[b]{0.325\textwidth}
         \centering
         \includegraphics[width=\textwidth]{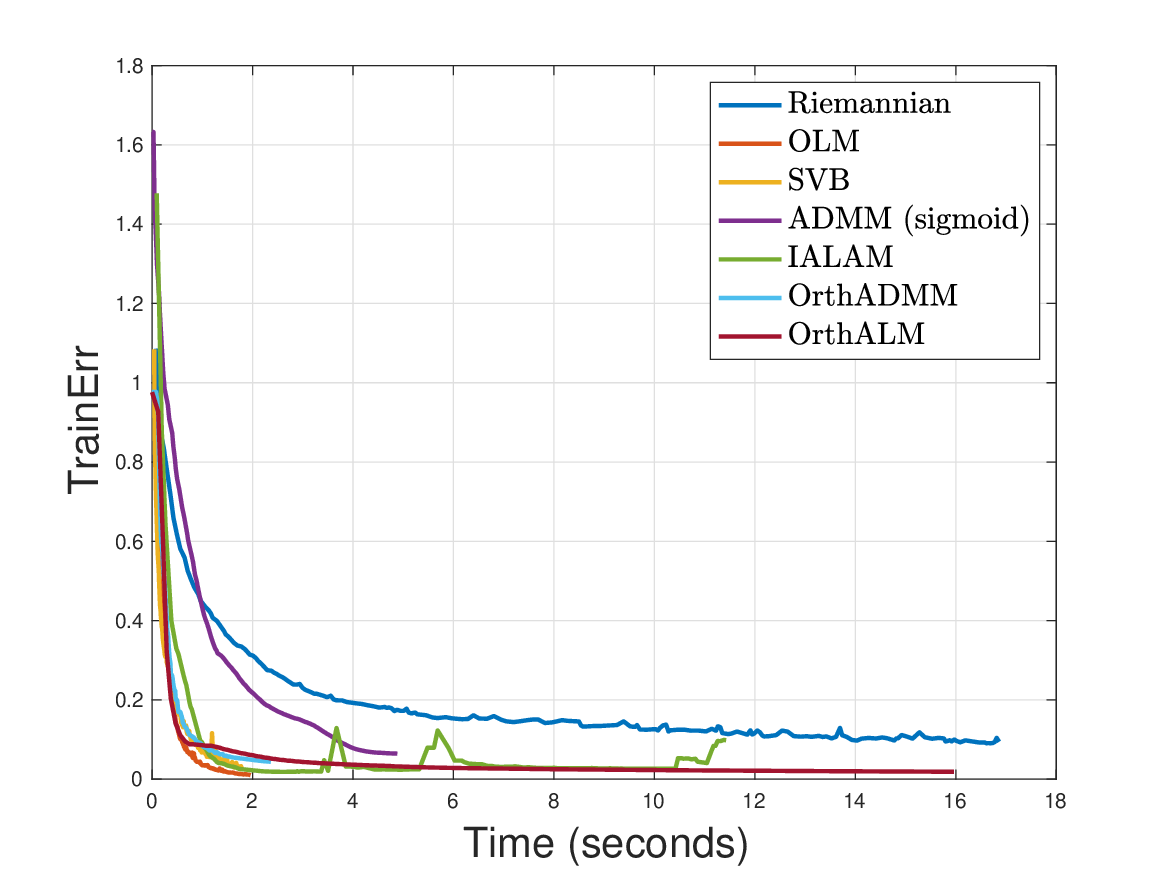}
         \caption{TrainErr}
         \label{9.3a}
     \end{subfigure}
     \hfill
     \begin{subfigure}[b]{0.325\textwidth}
         \centering
         \includegraphics[width=\textwidth]{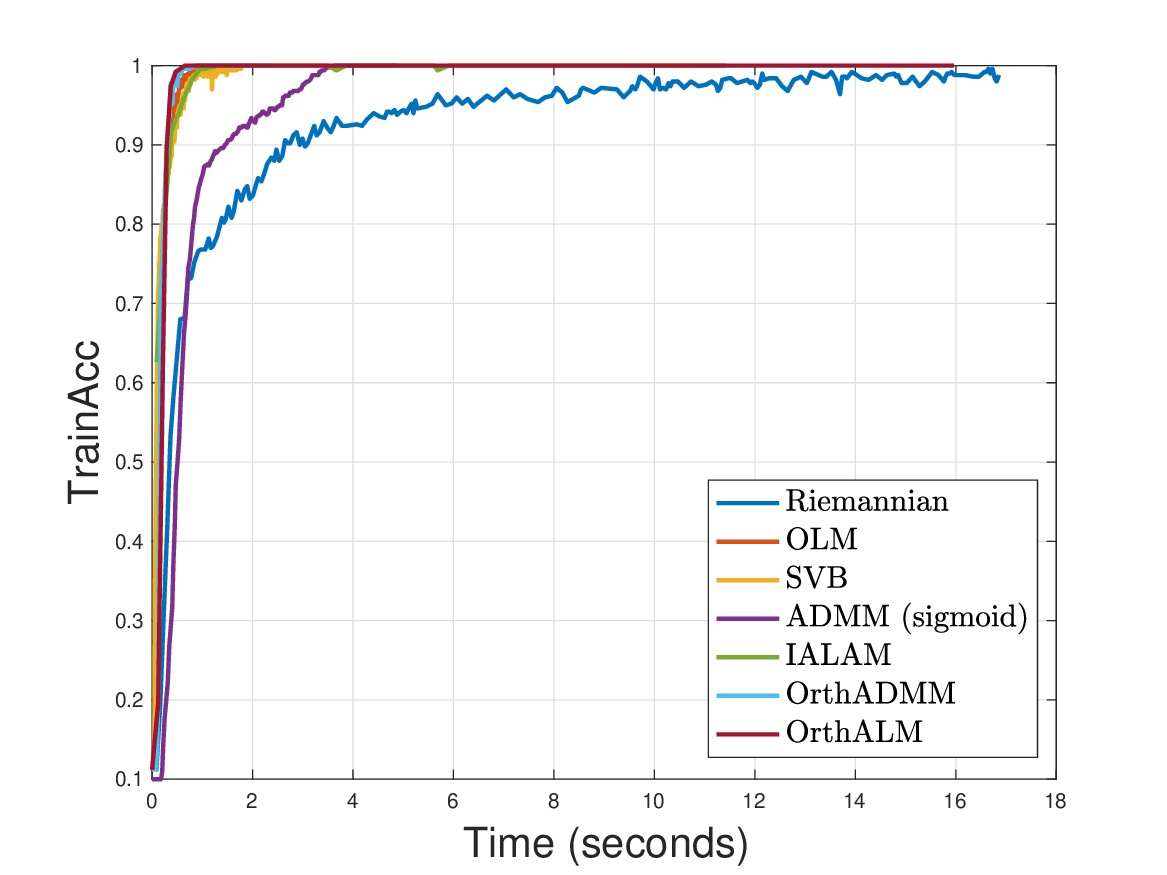}
         \caption{TrainAcc}
         \label{9.3b}
     \end{subfigure}
     \hfill
     \begin{subfigure}[b]{0.325\textwidth}
         \centering
         \includegraphics[width=\textwidth]{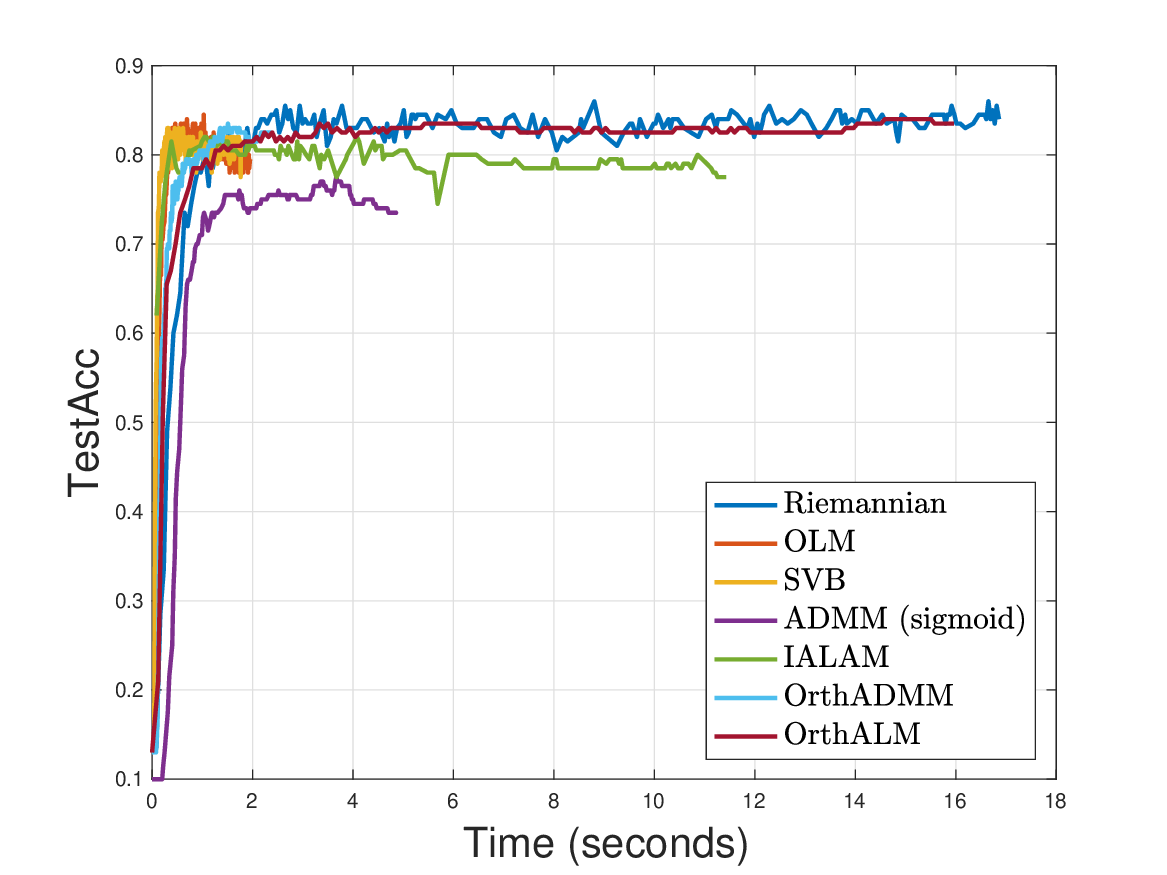}
         \caption{TestAcc}
         \label{9.3c}
     \end{subfigure}
     \hfill
     \begin{subfigure}[b]{0.325\textwidth}
         \centering
         \includegraphics[width=\textwidth]{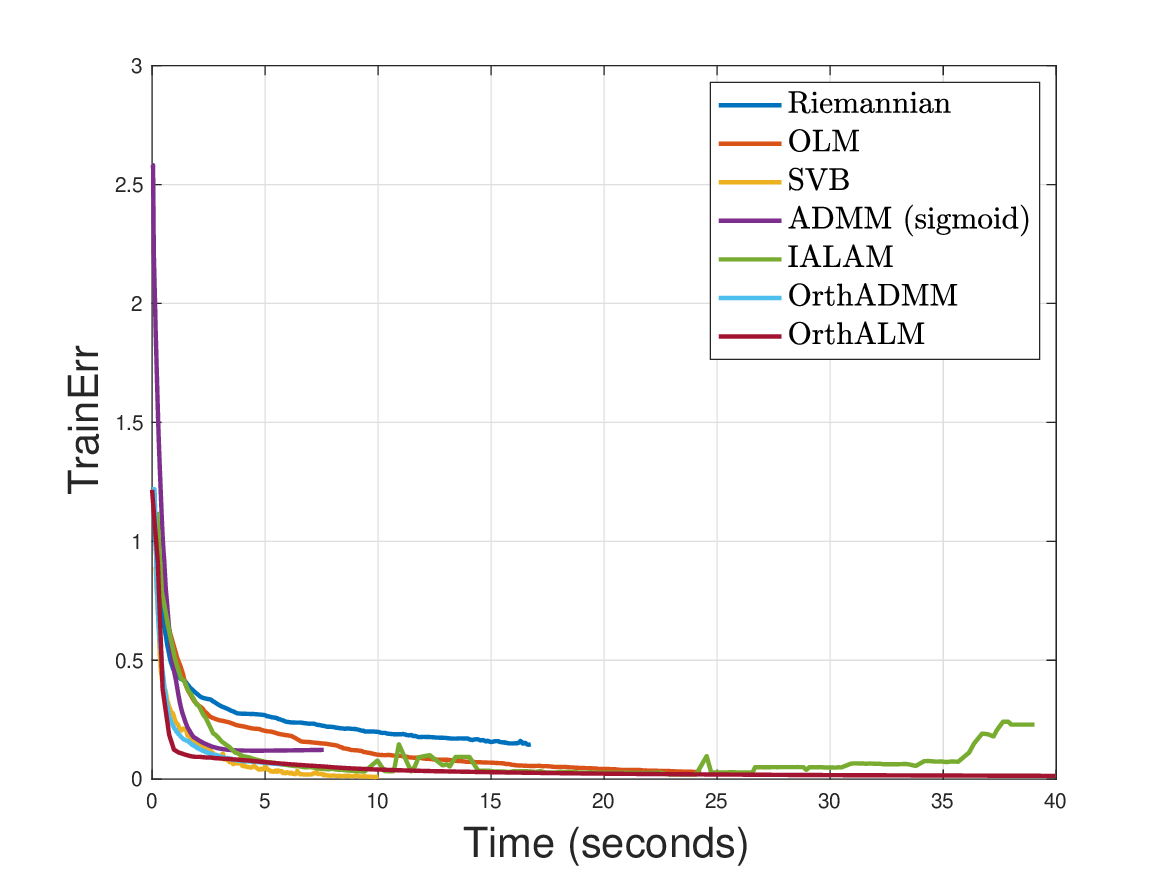}
         \caption{TrainErr}
         \label{9.3d}
     \end{subfigure}
     \hfill
     \begin{subfigure}[b]{0.325\textwidth}
         \centering
         \includegraphics[width=\textwidth]{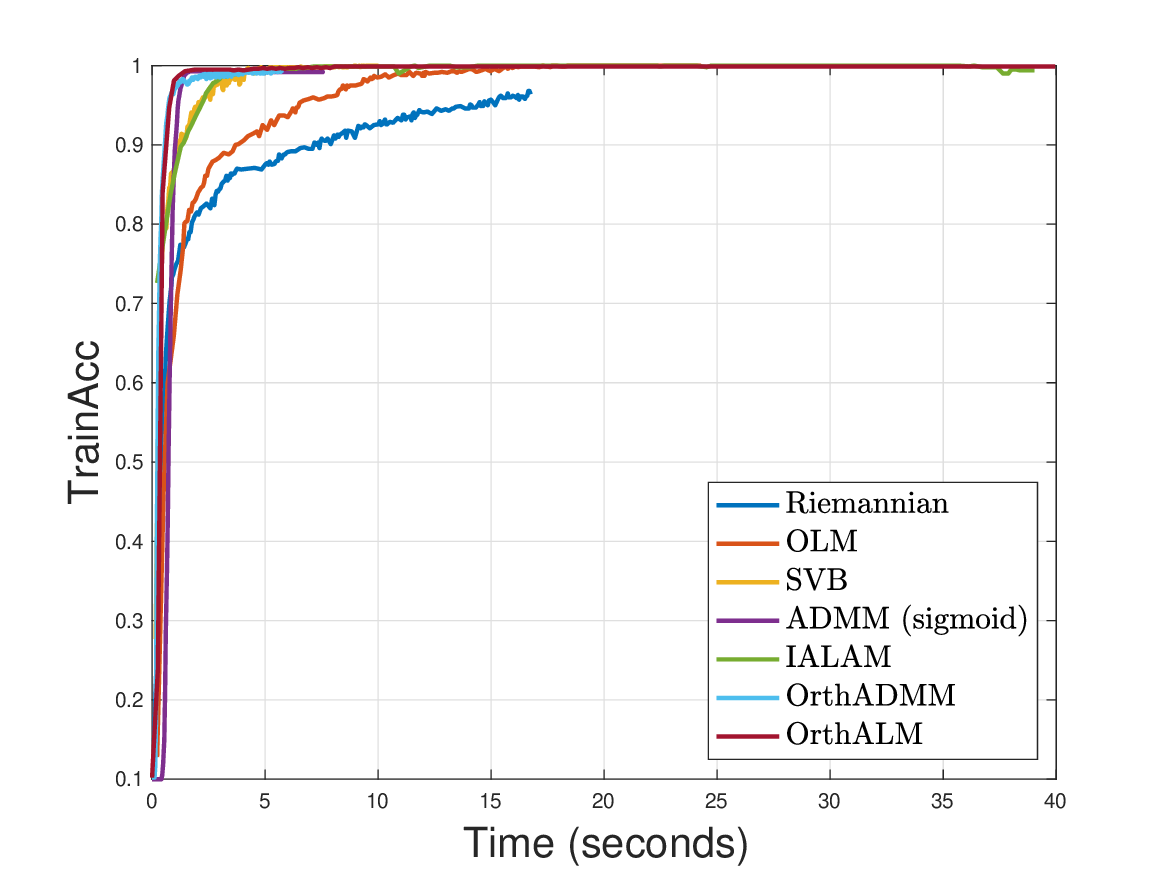}
         \caption{TrainAcc}
         \label{9.3e}
     \end{subfigure}
     \hfill
     \begin{subfigure}[b]{0.325\textwidth}
         \centering
         \includegraphics[width=\textwidth]{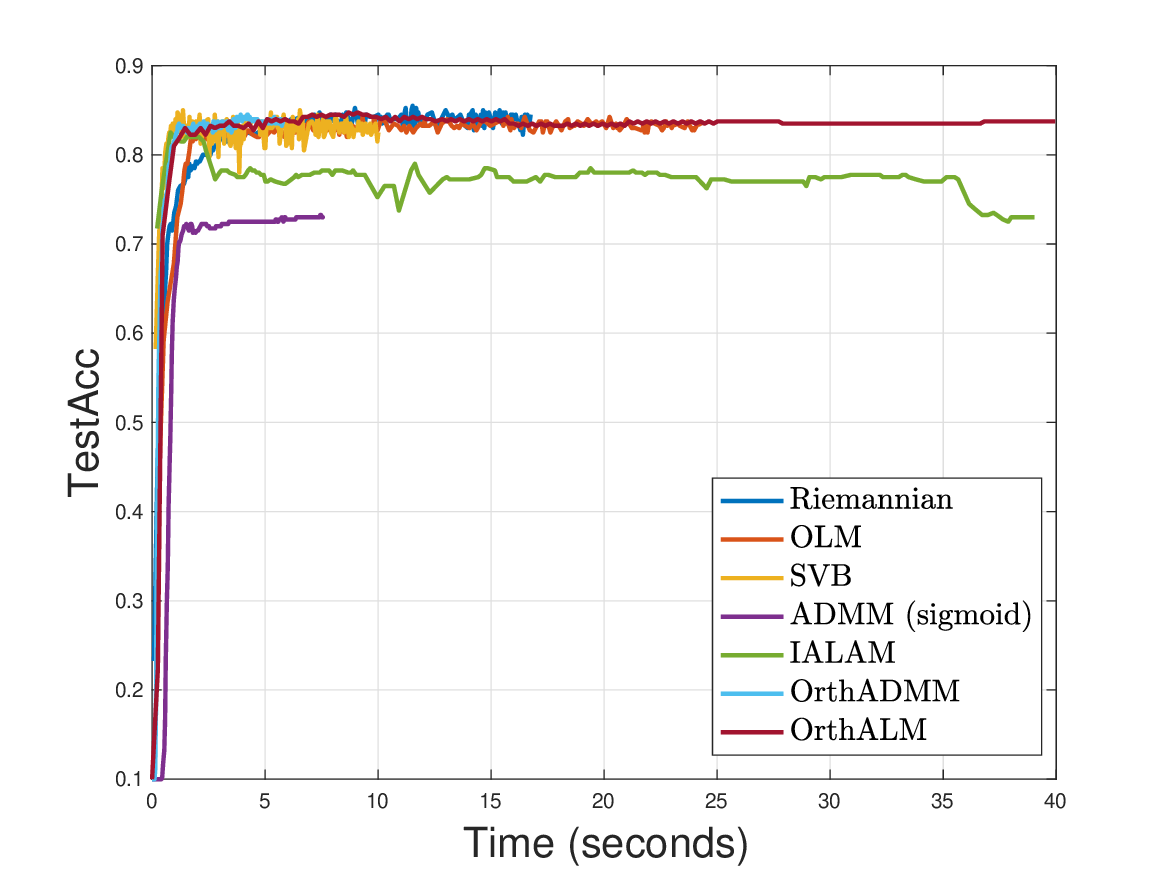}
         \caption{TestAcc}
         \label{9.3f}
     \end{subfigure}
     \hfill
     \begin{subfigure}[b]{0.325\textwidth}
         \centering
         \includegraphics[width=\textwidth]{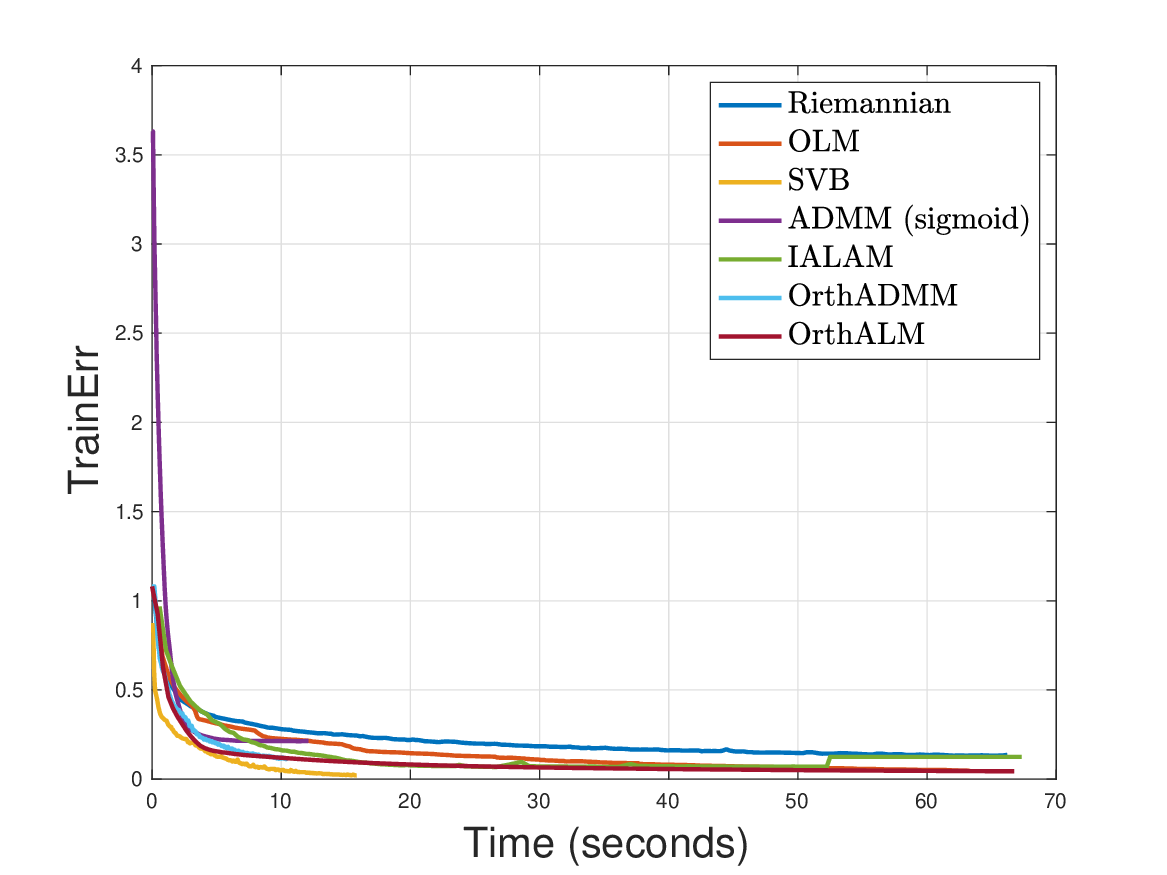}
         \caption{TrainErr}
         \label{9.3g}
     \end{subfigure}
     \hfill
     \begin{subfigure}[b]{0.325\textwidth}
         \centering
         \includegraphics[width=\textwidth]{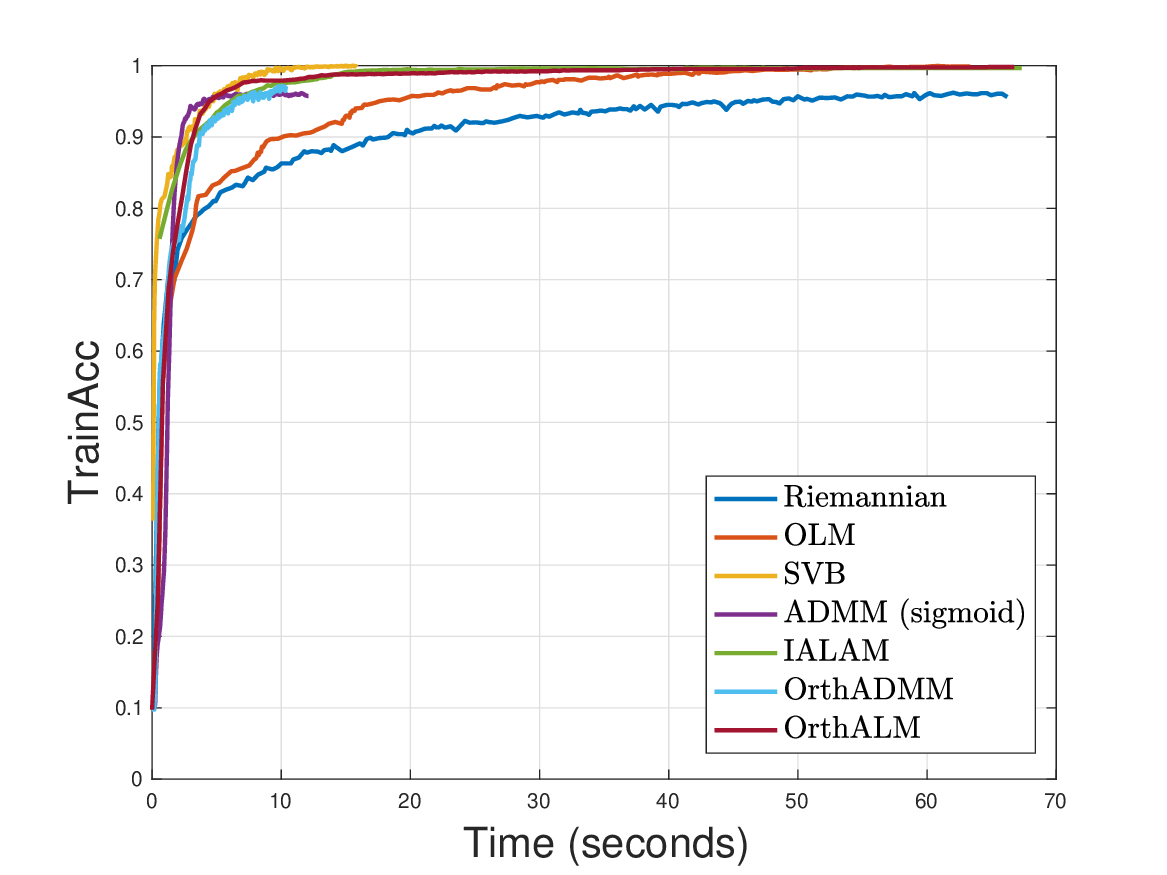}
         \caption{TrainAcc}
         \label{9.3h}
     \end{subfigure}
     \hfill
     \begin{subfigure}[b]{0.325\textwidth}
         \centering
         \includegraphics[width=\textwidth]{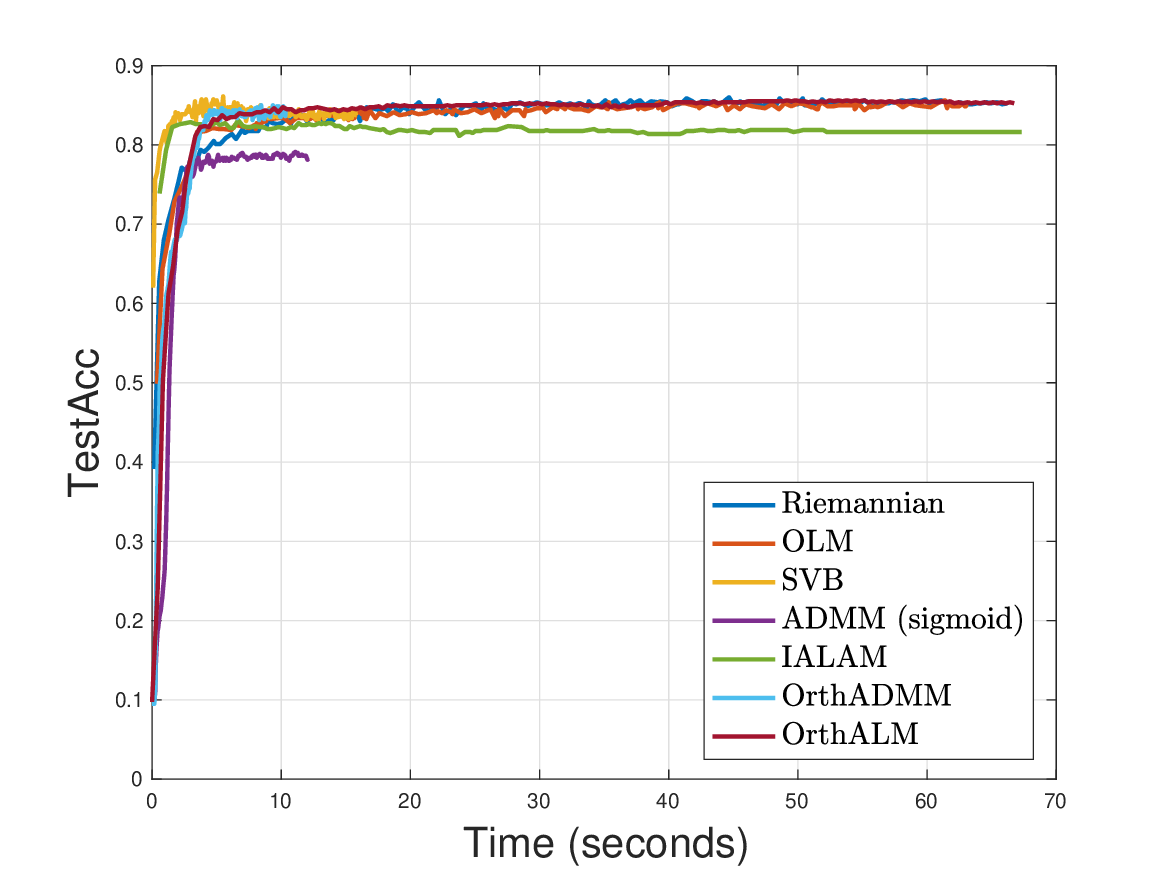}
         \caption{TestAcc}
         \label{9.3i}
     \end{subfigure}
     \hfill
     \begin{subfigure}[b]{0.325\textwidth}
         \centering
         \includegraphics[width=\textwidth]{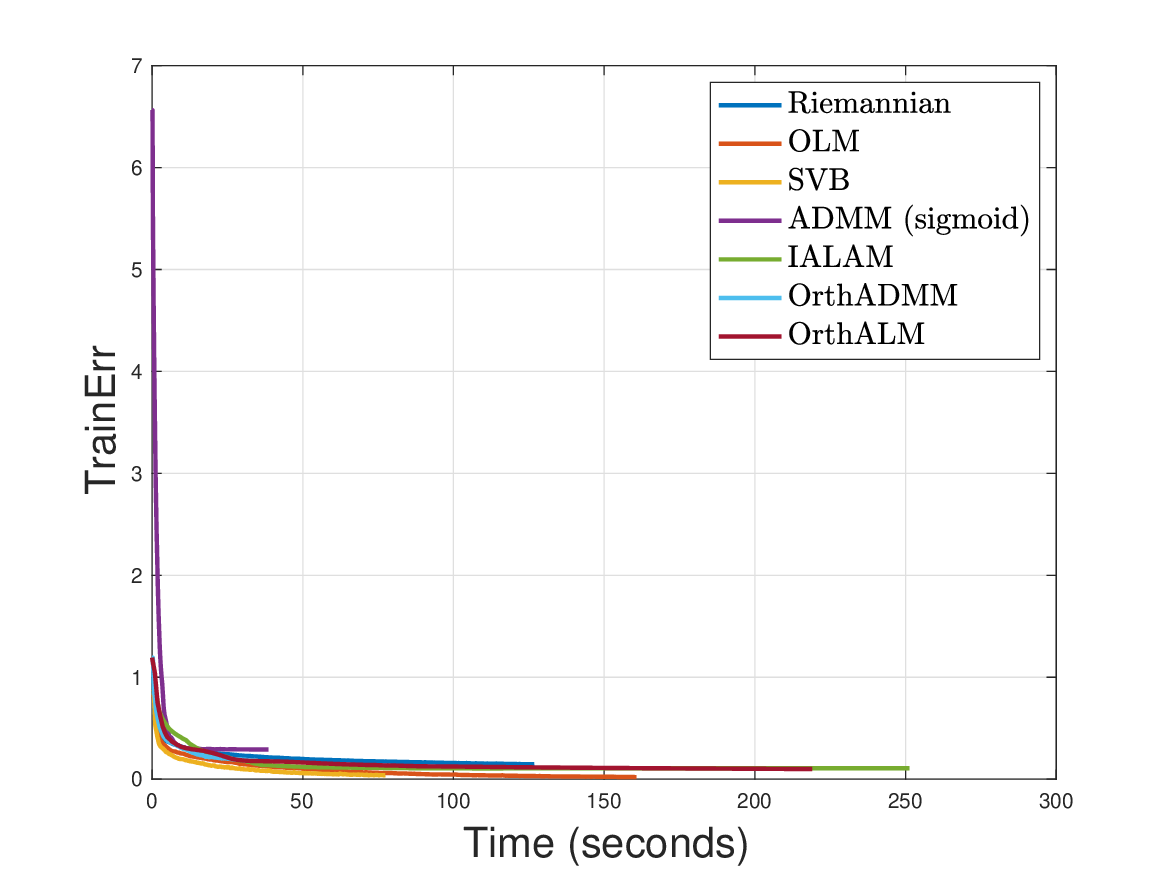}
         \caption{TrainErr}
         \label{9.3j}
     \end{subfigure}
     \hfill
     \begin{subfigure}[b]{0.325\textwidth}
         \centering
         \includegraphics[width=\textwidth]{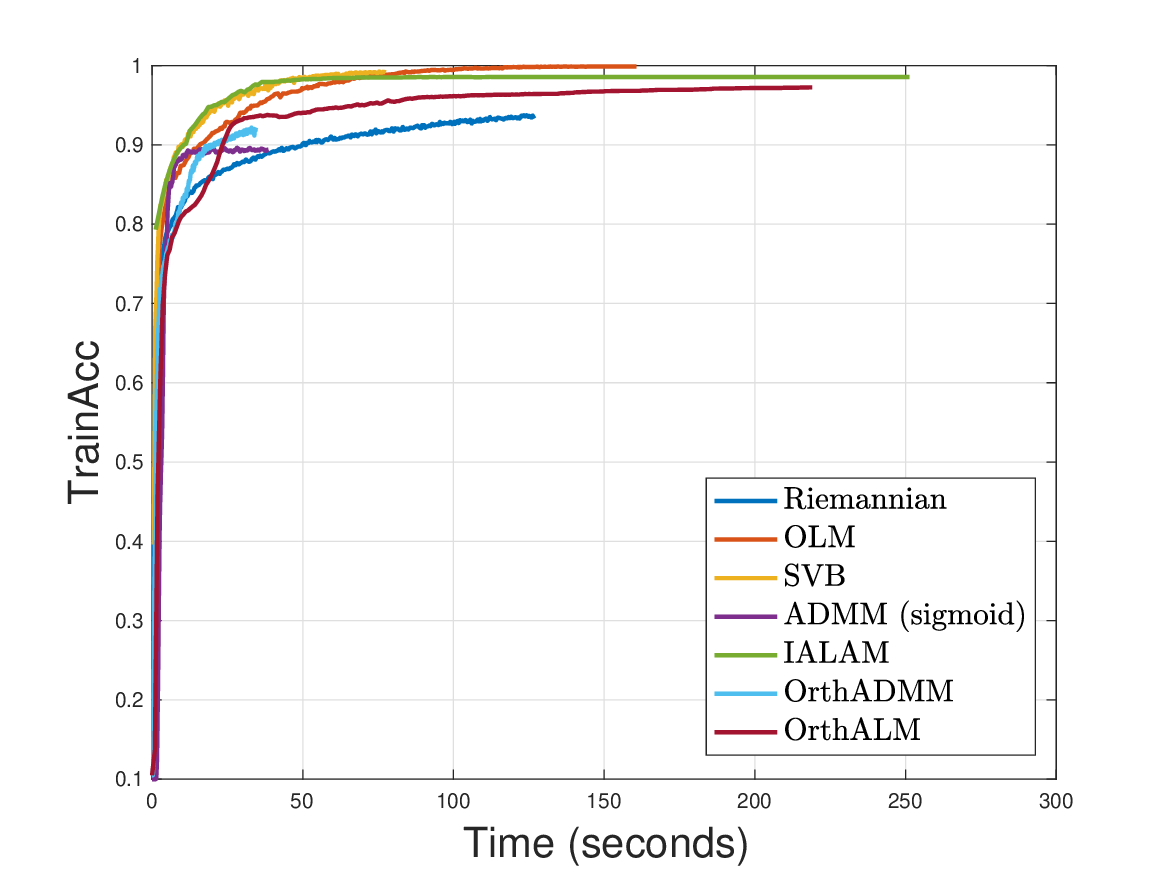}
         \caption{TrainAcc}
         \label{9.3k}
     \end{subfigure}
     \hfill
     \begin{subfigure}[b]{0.325\textwidth}
         \centering
         \includegraphics[width=\textwidth]{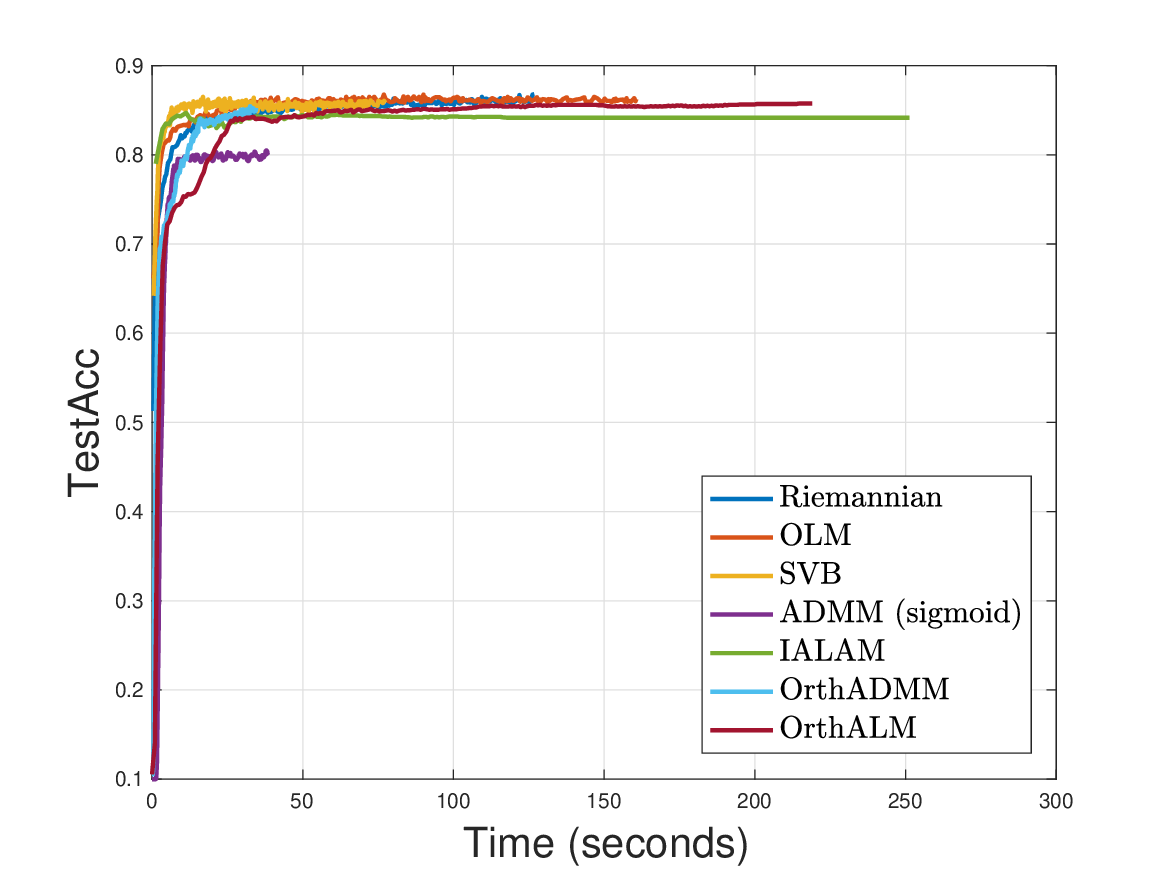}
         \caption{TestAcc}
         \label{9.3l}
     \end{subfigure}
        \caption{Comparison among OrthADMM, OrthALM, and five other existing approaches on the Fashion-MNIST dataset across different training sample sizes: $m = 500$ (first row, (a)--(c)), $m = 1000$ (second row, (d)--(f)), $m = 2000$ (third row, (g)--(i)), and $m = 5000$ (fourth row, (j)--(l)).}
        \label{fig:fmnist}
\end{figure}

We evaluate our OrthADMM and OrthALM alongside the five baseline algorithms using training sample sizes $m \in \{500, 1000, 2000, 5000\}$ and corresponding test sample sizes $m_{\text{test}} = 0.4m$. 
To ensure a balanced representation, we randomly pick up data entries from each class of the Fashion-MNIST under uniform distribution. 
We maintain the identical network architectures and maximum epoch budgets ($\text{epoch}_{\max}$) specified in \eqref{mnist.network.para}. 
The numerical results obtained at $\text{epoch}_{\max}$ are summarized in Table~\ref{tab:fmnist}, and the corresponding convergence trajectories of the training loss, training accuracy, and test accuracy with respect to CPU time are plotted in Figure~\ref{fig:fmnist}.
We draw several key observations from these results:
\begin{enumerate}

    \item As the training sample size $m$ increases from $500$ to $5000$, the test accuracy for all seven methods is ultimately higher at $m = 5000$. However, the overall gain in test accuracy on the Fashion-MNIST dataset across this range is notably smaller than the corresponding improvement observed on the MNIST dataset. Concurrently, as the classification task scales in both sample size and network width, the training accuracy of most algorithms experiences a mild decline while still remaining at a high level.

    \item Across all four evaluated sample sizes on the Fashion-MNIST dataset, both OrthADMM and OrthALM demonstrate highly competitive generalization capabilities that are comparable to the SGD-based approaches. 
    Furthermore, when compared to the two existing splitting and alternating-based methods for deep neural networks---IALAM and ADMM (sigmoid)---our proposed frameworks achieve higher test accuracies across every evaluated sample size $m$.

    \item A direct comparison between our two proposed algorithms reveals a clear trade-off between empirical fitting precision and computational efficiency. Across all sample sizes $m \in \{500, 1000, 2000, 5000\}$, OrthALM consistently achieves a lower training error and a training accuracy at least as high as that of OrthADMM; however, this performance advantage requires significantly greater CPU time.
    OrthADMM has the lowest reported CPU time for $m\in\{1000,2000,5000\}$. At $m=500$, it takes $2.4$ seconds, compared with $2.0$ seconds for OLM and $1.8$ seconds for SVB, while attaining higher test accuracy than both.

    \item As illustrated across all four rows of Figure~\ref{fig:fmnist}, both OrthADMM and OrthALM exhibit a sharp initial decline in training loss, as well as a rapid improvement in both training and test accuracies during the early stages of optimization.
    
\end{enumerate}

\begin{table}[htbp]
\caption{Performance comparison of our OrthADMM and OrthALM against five state-of-the-art baselines on the Fashion-MNIST dataset.}\label{tab:fmnist}
\begin{tabular}{llcccc}
\toprule
Method & Sample Size & TrainErr & TrainAcc & TestAcc & Time (s) \\
\midrule
OrthADMM & $\mathbf{m = 500}$  & 4.39e-2 & 1.00  & 0.835 & 2.4   \\
         & $\mathbf{m = 1000}$ & 6.18e-2 & 0.996 & 0.845 & 5.8   \\
         & $\mathbf{m = 2000}$ & 1.12e-1 & 0.972 & 0.850 & 10.4  \\
         & $\mathbf{m = 5000}$ & 1.79e-1 & 0.922 & 0.855 & 34.5  \\
\midrule
OrthALM  & $\mathbf{m = 500}$  & 1.87e-2 & 1.00  & 0.840 & 16.0  \\
         & $\mathbf{m = 1000}$ & 1.38e-2 & 0.999 & 0.838 & 39.9  \\
         & $\mathbf{m = 2000}$ & 4.41e-2 & 0.998 & 0.856 & 66.7  \\
         & $\mathbf{m = 5000}$ & 9.93e-2 & 0.973 & 0.858 & 219.0 \\
\midrule
OLM      & $\mathbf{m = 500}$  & 1.01e-2 & 1.00  & 0.795 & 2.0   \\
         & $\mathbf{m = 1000}$ & 2.99e-2 & 1.00  & 0.833 & 24.3  \\
         & $\mathbf{m = 2000}$ & 4.66e-2 & 0.999 & 0.854 & 63.3  \\
         & $\mathbf{m = 5000}$ & 2.03e-2 & 0.999 & 0.863 & 160.7 \\
\midrule
SVB      & $\mathbf{m = 500}$  & 2.55e-2 & 0.996 & 0.810 & 1.8   \\
         & $\mathbf{m = 1000}$ & 9.28e-3 & 1.00  & 0.825 & 10.0  \\
         & $\mathbf{m = 2000}$ & 2.09e-2 & 0.999 & 0.838 & 15.8  \\
         & $\mathbf{m = 5000}$ & 3.79e-2 & 0.991 & 0.859 & 77.4  \\
\midrule
Riemannian & $\mathbf{m = 500}$  & 9.07e-2 & 0.988 & 0.840 & 16.9  \\
           & $\mathbf{m = 1000}$ & 1.45e-1 & 0.964 & 0.845 & 16.7  \\
           & $\mathbf{m = 2000}$ & 1.32e-1 & 0.956 & 0.853 & 66.2  \\
           & $\mathbf{m = 5000}$ & 1.46e-1 & 0.933 & 0.868 & 126.8 \\
\midrule
IALAM    & $\mathbf{m = 500}$  & 1.32e-2 & 1.00  & 0.800 & 11.4  \\
         & $\mathbf{m = 1000}$ & 9.88e-3 & 1.00  & 0.783 & 39.0  \\
         & $\mathbf{m = 2000}$ & 3.39e-2 & 0.997 & 0.819 & 67.3  \\
         & $\mathbf{m = 5000}$ & 5.26e-2 & 0.986 & 0.844 & 251.3 \\
\midrule
ADMM (sigmoid) & $\mathbf{m = 500}$  & 6.45e-2 & 1.00  & 0.775 & 4.9   \\
               & $\mathbf{m = 1000}$ & 1.20e-1 & 0.992 & 0.733 & 7.6   \\
               & $\mathbf{m = 2000}$ & 2.13e-1 & 0.962 & 0.791 & 12.1  \\
               & $\mathbf{m = 5000}$ & 2.91e-1 & 0.896 & 0.805 & 38.6  \\
\botrule
\end{tabular}
\end{table}

\section{Conclusion}\label{sec:conclusion}

This paper investigated the design, theoretical convergence, and numerical performance of splitting and alternating-based algorithms for training OrthDFNs with nonsmooth activations. 
We reformulated network training into an equality-constrained minimization~\eqref{ori.prob.orth.three.sets.vari} and its $\ell_2$-regularized counterpart~\eqref{prob.orth.three.sets.vari}, proposing OrthALM to solve the latter. Furthermore, we introduced OrthADMM as a fast, proximal-free variant for the unregularized problem~\eqref{ori.prob.orth.three.sets.vari}. In both algorithms, every block coordinate subproblem admits a closed-form solution. 
Theoretically, we proved an $\mathcal{O}(1/\epsilon^2)$ iteration complexity for the inner algorithm, and guaranteed that any accumulation point generated by OrthALM is a KKT point of problem~\eqref{prob.orth.three.sets.vari}. 
Empirically, OrthALM achieved lower training errors and constraint violations in most tested configurations, whereas OrthADMM delivered competitive performance with substantial computational speedups over state-of-the-art methods.
Future research will pursue two primary directions. First, establishing the global convergence of OrthADMM remains an open challenge. Second, extending the orthogonality constraints and our algorithmic framework to broader architectures, such as convolutional and recurrent neural networks, offers potential for enhancing the generalization ability and training efficiency of deep learning models.

\backmatter

\section*{Acknowledgements}
The authors would like to thank Prof. Ren-Cang Li for his
valuable guidance and insightful discussions.

\section*{Statements and Declarations}

\begin{itemize}
\item \textbf{Funding:} This work was supported in part by the Hong Kong Research Grants Council General Research Fund under Grant 12300322.

\vspace{1em}

\item \textbf{Conflict of interest/Competing interests:}
The authors have no relevant financial or non-financial
interests to disclose.

\vspace{1em}

\item \textbf{Data availability:} The Extended Yale B, MNIST, and Fashion-MNIST datasets analyzed during the current study are publicly available benchmark datasets cited in the manuscript.

\vspace{1em}

\item \textbf{Code availability:} Our MATLAB source code is publicly available on GitHub at \url{https://github.com/WANGJiawei000/OrthADMM_and_OrthALM_for_OrthDFNs}.

\end{itemize}

\bibliography{reference,refLi}

@book{absil2008optimization,
  title={Optimization {A}lgorithms on {M}atrix {M}anifolds},
  author={Absil, P.-A. and Mahony, Robert and Sepulchre, Rodolphe},
  year={2008},
  publisher={Princeton University Press},
  address={Princeton, NJ}
}

@article{li2024theory,
  title={A {T}heory of the {NEP}v {A}pproach for {O}ptimization {O}n the {S}tiefel {M}anifold},
  author={Li, Ren-Cang},
  journal={Foundations of Computational Mathematics},
  volume={26},
  number={1},
  pages={179--244},
  year={2026},
  publisher={Springer}
}

@book{rockafellar2009variational,
  title={Variational Analysis},
  author={Rockafellar, R. Tyrrell and Wets, Roger J.-B.},
  year={1998},
  publisher={Springer Science \& Business Media},
  address={Berlin, Heidelberg}
}

@article{georghiades2001from,
  title={From few to many: {I}llumination cone models for face recognition under variable lighting and pose},
  author={Georghiades, Athinodoros S. and Belhumeur, Peter N. and Kriegman, David J.},
  journal={IEEE Transactions on Pattern Analysis and Machine Intelligence},
  volume={23},
  number={6},
  pages={643--660},
  year={2001},
  publisher={IEEE}
}

@article{lee2005acquiring,
  title={Acquiring linear subspaces for face recognition under variable lighting},
  author={Lee, Kuang-Chih and Ho, Jeffrey and Kriegman, David J},
  journal={IEEE Transactions on Pattern Analysis and Machine Intelligence},
  volume={27},
  number={5},
  pages={684--698},
  year={2005},
  publisher={IEEE}
}

@article{li2019orthogonal,
  title={Orthogonal deep neural networks},
  author={Li, Shuai and Jia, Kui and Wen, Yuxin and Liu, Tongliang and Tao, Dacheng},
  journal={IEEE Transactions on Pattern Analysis and Machine Intelligence},
  volume={43},
  number={4},
  pages={1352--1368},
  year={2021},
  publisher={IEEE}
}

@article{huang2018orthogonal, title={Orthogonal {W}eight {N}ormalization: {S}olution to {O}ptimization {O}ver {M}ultiple {D}ependent {S}tiefel {M}anifolds in {D}eep {N}eural {N}etworks}, volume={32}, url={https://ojs.aaai.org/index.php/AAAI/article/view/11768}, DOI={10.1609/aaai.v32i1.11768}, number={1}, journal={Proceedings of the AAAI Conference on Artificial Intelligence}, author={Huang, Lei and Liu, Xianglong and Lang, Bo and Yu, Adams and Wang, Yongliang and Li, Bo}, year={2018}, month={Apr.} }

@inproceedings{carreira2014distributed,
  title={Distributed optimization of deeply nested systems},
  author={Carreira-Perpi\~{n}\'{a}n, Miguel and Wang, Weiran},
  booktitle={Artificial intelligence and statistics},
  pages={10--19},
  year={2014},
  organization={PMLR}
}

@article{zeng2021admm,
  title={On {A}{D}{M}{M} in deep learning: {C}onvergence and saturation-avoidance},
  author={Zeng, Jinshan and Lin, Shao-Bo and Yao, Yuan and Zhou, Ding-Xuan},
  journal={Journal of Machine Learning Research},
  volume={22},
  number={199},
  pages={1--67},
  year={2021}
}

@article{wang2025augmented,
  title={An augmented {L}agrangian method for training recurrent neural networks},
  author={Wang, Yue and Zhang, Chao and Chen, Xiaojun},
  journal={SIAM Journal on Scientific Computing},
  volume={47},
  number={1},
  pages={C22--C51},
  year={2025},
  publisher={SIAM}
}

@article{liu2023inexact,
  title={An inexact augmented {L}agrangian algorithm for training leaky {R}e{L}{U} neural network with group sparsity},
  author={Liu, Wei and Liu, Xin and Chen, Xiaojun},
  journal={Journal of Machine Learning Research},
  volume={24},
  number={212},
  pages={1--43},
  year={2023}
}

@article{elman1990finding,
  title={Finding structure in time},
  author={Elman, Jeffrey L},
  journal={Cognitive science},
  volume={14},
  number={2},
  pages={179--211},
  year={1990},
  publisher={Wiley Online Library}
}

@misc{lecun1998mnist,
  title={The {M}{N}{I}{S}{T} database of handwritten digits},
  author={LeCun, Yann},
  url={http://yann.lecun.com/exdb/mnist/},
  year={1998}
}

@article{xiao2017fashion,
  title={{Fashion-MNIST}: a novel image dataset for benchmarking machine learning algorithms},
  author={Xiao, Han and Rasul, Kashif and Vollgraf, Roland},
  journal={arXiv preprint arXiv:1708.07747},
  year={2017}
}

@article{saxe2013exact,
  title={Exact solutions to the nonlinear dynamics of learning in deep linear neural networks},
  author={Saxe, Andrew M and McClelland, James L and Ganguli, Surya},
  journal={arXiv preprint arXiv:1312.6120},
  year={2013}
}

@article{lecun1998gradient,
  title={Gradient-based learning applied to document recognition},
  author={LeCun, Yann and Bottou, L{\'e}on and Bengio, Yoshua and Haffner, Patrick},
  journal={Proceedings of the IEEE},
  volume={86},
  number={11},
  pages={2278--2324},
  year={1998},
  publisher={IEEE}
}

\end{document}